\documentclass{article}

\usepackage{arxiv}

\usepackage[normalem]{ulem} 

\usepackage[utf8]{inputenc} 
\usepackage[T1]{fontenc}    
\usepackage{hyperref}       
\usepackage{url}            
\usepackage{booktabs}       
\usepackage{amsfonts}       
\usepackage{microtype}      
\usepackage{color}
\usepackage{lipsum}
\usepackage{amsthm}
\usepackage{indentfirst}
\usepackage{graphicx}
\usepackage{epstopdf}
\usepackage{fourier}
\usepackage{bm}
\usepackage{mathtools}
\usepackage{latexsym,enumerate}
\usepackage{ulem}

\def\bm{\boldsymbol}
\newcommand{\BEA}{\begin{eqnarray}}
\newcommand{\EEA}{\end{eqnarray}}

\newcommand{\peq}{p_{eq}}
\newcommand{\td}{\operatorname{d}}

\newcommand{\BR}{\mathbb{R}}

\newcommand{\id}{\operatorname{id}}

\newtheorem{theorem}{Theorem}[section]
\newtheorem{corollary}{Corollary}

\newtheorem{lemma}[theorem]{Lemma}
\newtheorem{proposition}{Proposition}

\theoremstyle{definition}
\newtheorem{definition}[theorem]{Definition}
\newtheorem{remark}{Remark}

\title{Estimating linear response statistics using orthogonal polynomials: An RKHS formulation}

\author{
  He Zhang \\
  Department of Mathematics \\
  The Pennsylvania State University, University Park, PA 16802, USA\\
  \texttt{hqz5159@psu.edu} \\
  \And
  John Harlim \\
  Department of Mathematics, Department of Meteorology and Atmospheric Science, \\ Institute for Computational and Data Sciences \\
  The Pennsylvania State University, University Park, PA 16802, USA\\
  \texttt{jharlim@psu.edu} \\
  \And
  Xiantao Li \\
  Department of Mathematics\\
  The Pennsylvania State University, University Park, PA 16802, USA\\
  \texttt{xxl12@psu.edu} \\
}

\begin{document}

\maketitle

\begin{abstract}
In the paper, we study the problem of estimating linear response statistics under external perturbations using time series of unperturbed dynamics. A standard approach to this estimation problem is to employ the Fluctuation-Dissipation Theory, which, in turn, requires the knowledge of the functional form of the underlying unperturbed density that is not available in general. To overcome this issue, we consider a nonparametric density estimator formulated by the kernel embedding of distributions. To avoid the computational expense associated with using radial type kernels, we consider the ``Mercer-type'' kernels constructed based on the classical orthogonal bases defined on non-compact domains, such as the Hermite and Laguerre polynomials. While the resulting representation is analogous to Polynomial Chaos Expansion (PCE), by studying the orthogonal polynomial approximation in the reproducing kernel Hilbert space (RKHS) setting, we establish the uniform convergence of the estimator. More importantly,  the RKHS formulation allows one to systematically address a practical question of identifying the PCE basis for a consistent estimation through the decay property of the target functions that can be quantified using the available data.
In terms of the linear response estimation, our study provides practical conditions for the well-posedness of not only the estimator but also the well-posedness of the underlying response statistics. Given a well-posed estimator, we provide a theoretical guarantee for the convergence of the estimator to the underlying linear response statistics. Finally, we provide a statistical error bound for the density estimation that accounts for the Monte-Carlo averaging over non-i.i.d time series and the biases due to a finite basis truncation. This error bound provides a means to understand the feasibility as well as limitation of the kernel embedding with Mercer-type kernels. Numerically, we verify the effectiveness of the kernel embedding linear response estimator on two stochastic dynamics with known, yet, non-trivial equilibrium densities.
\end{abstract}


\section{Introduction} \label{sec:intro}

Estimating linear response statistics of dynamical systems under external forces is a problem of broad interest. This forward Uncertainty Quantification (UQ) problem has many important applications. For example, in climate dynamics, the linear response can be used as a proxy that quantifies the climate change statistics corresponding to exogenous forcing such as the volcanic eruptions or even the anthropogenic factor such as the human activities \cite{leith:75,mag:05}.  In statistical mechanics, the linear response provides a simple route to determine transport coefficients from microscopic fluctuations \cite{Toda-Kubo-2}, such as viscosity, diffusion coefficients, and heat conductivity via Green-Kubo type of formulas \cite{green1954markoff,kubo1957statistical}.  One of the popular approaches to quantify the linear response statistics is using the Fluctuation-Dissipation theory (FDT), which in rough terms,  states that the leading order (linear) statistical response to small perturbations can be approximated by a two-point statistic of the unperturbed dynamics. In this paper, we consider computing linear response in the context of ergodic stochastic differential equations, the validity of which was studied in \cite{hairer2010simple}.

In practice, while the relevant two-point statistics can be numerically estimated by a Monte-Carlo average over samples of unperturbed dynamics at equilibrium state, the integrand (or the function to be averaged) depends on the explicit expression of the equilibrium density of the unperturbed dynamical system which is unknown in general. Examples include PDEs that exhibit spatial-temporal chaos \cite{jolly1990approximate}, 
non-equilibrium steady states in statistical mechanics, \cite{evans1984nonlinear,baiesi2009fluctuations}, coarse-graining molecular dynamics \cite{roux1995calculation}, where the potential of mean force has to be estimated from data, etc.

The main problem that arises in such applications is the density estimation of the unknown equilibrium distribution of the unperturbed dynamics, which in turn, allows us to compute the linear response statistics via the FDT theory. Density estimation is a long-standing problem in computational statistics and machine learning. There are many approaches developed to tackle this fundamental problem and one can classify them based on the form of the models such as parametric, non-parametric, or semi-parametric, or based on the training methodology such as linear or nonlinear estimators. The parametric or semi-parametric type models are suitable for problems where some physical knowledge is known. In the absence of physical knowledge, a popular parametric density estimator is the Gaussian Mixture Models (which is also known as the Radial Basis Models in some literature) \cite{hwang1994nonparametric}. Regardless of whether physical knowledge is known or not, most of them are nonlinear estimators since their training phase involves a minimization procedure to estimate the latent parameters. Two practical issues of such nonlinear estimators are: (1) The difficulty in finding the global minimizer using numerical methods that provably converge to local minima, especially when the parameter space is high-dimensional; (2) The identifiability of these parameters when the form of the underlying density function is not known. While these issues are not solved, recent successes of Deep Neural Network (DNN) suggest that high-dimensional functions can be approximated effectively using compositions of ReLU or sigmoid-type functions, even when the global minimizer is not found \cite{caogu:2019}. In this direction, there are several recently introduced density estimators that have adopted DNN, such as the Neural Autoregressive Distribution Estimation \cite{uria2016neural} and its variant, the Masked Autoregressive Flow \cite{papamakarios2017masked}. We should point out that with this type of estimators, the complexity of each function evaluation is on the order of $M$, where $M$ is the number of parameters in the neural network model \cite{uria2016neural,wang2019nonparametric}. For our application, as we shall see, since we need to evaluate the derivatives of the density on each training data point of a set at least on the order of $N=\mathcal{O}(10^7)$, a small $M$ will be desirable. While nonlinear estimator is a promising direction that warrants a thorough investigation, we will not pursue it in the current paper.

In contrast to these approaches, the nonparametric approach requires the least modeling assumptions and our goal is to construct such an estimator that is consistent in the limit of large data. Among the existing approaches, it is widely accepted that the classical Kernel Density Estimation (KDE) \cite{rosenblatt1956remarks} is not effective for problems with dimension higher than three (see e.g., \cite{hwang1994nonparametric, bsp:2013,wang2019nonparametric}). In addition to the bandwidth specification issue, the KDE implemented with radial-type kernels, is not practical for our application. In particular, given a training data set of size $N$, which can be of order $10^7$ or larger depending on applications, the estimator will be defined as an average of the radial functions $\left\{k(x,x_i) = h(\|x-x_i\|)\right\}_{i=1}^N$, where $\|\cdot\|$ denotes, e.g, the $d$-dimensional Euclidean norm for some $h>0$. With such an estimator, computing the linear response statistics will require evaluating the norm of the distances between pairs of training data points roughly $N^2/2$ times in addition to evaluating $h$ and its derivatives on each of these $N^2/2$ distances. Alternatively, a theoretically consistent nonparametric approach that also provably avoids the curse of dimensionality is the Bayesian Sequential Partitioning (BSP) \cite{bsp:2013,NIPS2017_7059}, which involves the Sequential Important Sampling algorithm. While the dimensionality aspect of BSP is appealing, the fact that the estimator is in the form of piecewise constant functions supported on binary partitions makes it not suitable for our application that requires the functional form of the derivatives of the density. While nonparametric methods that directly estimate the derivative of log density exists (e.g.~\cite{sriperumbudur2017density}), such techniques will restrict our problem to exponential type densities. In addition, the method proposed in \cite{sriperumbudur2017density} requires an inversion of a matrix of size $Nd\times Nd$, which is not feasible for problems with large $N$. 
 
In this paper, we consider the kernel embedding of distributions \cite{muandet2017kernel}, which is a linear density estimator. We should point out that the concept of kernel embedding of distribution was introduced in \cite{muandet2017kernel} to characterize probability distributions with the {\it kernel mean embedding}, which are nothing but statistical quantities of relevant feature maps corresponding to a reproducing kernel Hilbert space (RKHS). Here, we use kernel embedding to estimate probability density functions. In particular, the estimator will be represented as a linear superposition of the RKHS basis functions with the kernel mean embedding as the expansion coefficients. For this application, the direct use of the kernel embedding estimator with radial-type kernels is not practical. Beyond the same computational issue that hampers the KDE method, the kernel embedding implemented with radial-type kernels requires an inversion of an $N\times N$ matrix. To avoid the complexity of the function evaluations and the large matrix inversion, we consider the ``Mercer-type'' kernels constructed using the classical orthogonal polynomials of weighted $L^2$-space. Practically, the orthogonality replaces the large matrix inversion problem with a transpose operation. Such representation allows us to conveniently construct a hypothesis model with $M$ expansion coefficients, with $M\ll N$
such that the calculation of linear response statistics amounts to evaluating $M$ polynomial basis functions (instead of $N/2$ radial basis functions) on $N$ training data points.

We should point out that the resulting estimator (kernel embedding with polynomial basis) is essentially the polynomial chaos expansion (when classical orthogonal polynomials are used, e.g., Hermite polynomials), whose convergence is often understood in $L^2$-sense, relying on the Cameron-Martin theorem. See e.g., \cite{xiu:2010,ernst2012convergence} for specific conditions for the convergence. In this paper, we will show that the uniform convergence of the polynomial chaos expansion can be naturally achieved using the theory of Reproducing Kernel Hilbert Space. To achieve this goal, we develop an RKHS using the ``Mercer-type" kernels induced by orthogonal polynomials, which then allows its components to inherit the ``nice'' properties of the kernel, such as boundedness and smoothness. As opposed to the compact domain setting for which the Mercer theorem is valid, these desirable properties, specifically, boundedness and smoothness, are not automatically inherited by the kernel definition alone when its domain is not a compact metric space, which leads us to call the constructed kernel as the ``Mercer-type''.  One of the practical issues in polynomial chaos expansion is how to choose the appropriate basis when the underlying distribution is different than the basic distributions as listed in \cite{xiu:2010}. This issue can be naturally understood in the RKHS setting as a problem of determining how ``large''  the resulting RKHS space needs to be to guarantee a consistent estimator. To understand this, we generalize the notion of $c_0$-universality of the kernel that was introduced in \cite{sriperumbudur2010relation,sriperumbudur2011universality} on space of continuous functions on $\mathbb{R}^d$ with appropriate decay rate. Our study provides a systematic way to choose the basis (or through its corresponding weight function) for approximating density functions with a decaying rate of Gaussian or faster. 

While the kernel embedding density estimator is consistent in the limit of large $M$ and $N$, a finite truncation of these parameters poses an undesirable issue, namely, admitting negative values in estimating positive functions defined on non-compact domains. In principle, we have traded the positivity of the estimators (that is guaranteed using the radial-type kernels) with orthogonality for computational efficiency. In this study, we will provide a detailed characterization of a compact subset for which the density estimate admits only positive values and, most importantly, the linear response statistics can be estimated by a well-defined estimator restricted to this subset up to any desirable precision. Practically, the well-posed estimator on this restricted domain can be numerically realized with a straightforward parameter selection criterion. We will use the designed criterion to identify $M$ such that the linear response estimate based on the training data in the restricted domain is an accurate estimator.

Finally, we also provide a statistical error bound for the density estimation that accounts for the Monte-Carlo discretization, averaging over non-i.i.d time series with $\alpha$-mixing property \cite{Davydov:68,Hang:14}. This error bound provides a means to understand the feasibility as well as limitation of the kernel embedding with Mercer-type kernels in general problems. This error bound serves as a definitive confirmation to the limitation of the polynomial chaos expansion that were reported in many numerical studies, such as \cite{branicki2013fundamental} in the context of forward UQ and \cite{lu2015limitations} in the context of Bayesian inference. It is, in fact, well-known in statistics that the optimal rate of any linear estimators is of order-$N^{-\frac{2\ell}{2\ell+d}}$ \cite{stone1982optimal}, where the parameter $\ell$ denotes the smoothness of the function. This means that only when the function is very smooth, such as $\ell=d$, the error bound becomes independent of dimension, $N^{-\frac{2}{3}}$. Indeed, for data that lie on a smooth $d$-dimensional manifold embedded in $\mathbb{R}^n$, one can construct a Mercer-type kernel based on the orthogonal basis functions obtained via the diffusion maps algorithm \cite{cl:06} and achieve the optimal rate as reported in \cite{stone1982optimal}. While this approach has been explored in \cite{JH:19}, the errors in the estimation of these eigenbasis, $\mathcal{O}\left(\left(\frac{\log N}{N}\right)^{\frac{1}{2d}}\right)$, do not escape the curse of dimension \cite{trillos2019error}. Despite this fundamental limitation, our study sheds light on what one can achieve from implementing a linear estimator in approximating linear response statistics.

\paragraph{Outline and Main Contributions}
We close this introduction with an outline of the rest of the paper, summarizing the main contributions of each section. Main results:

\begin{itemize}

\item In Section~\ref{sec:review}, we provide a quick review of the FDT linear response theory and relevant results on kernels and RKHS. To have a well-defined estimation problem, we deduce sufficient conditions to guarantee that the underlying FDT linear response operator is bounded uniformly in time (see Lemma~\ref{lem:bound_kA}). The main novel contribution in this section is the generalization of the $c_0$-universality of kernels to weighted $c_0$-universality, which will be used to characterize the denseness of a given RKHS in the space of continuous functions with  certain decay rate. We derive the connection between $c_0$-universal kernels and weighted $c_0$-universal kernels (see Lemma~\ref{lem:weight_c0}) so that the principles in $c_0$-universality can be shifted to the weighted scenario.

\item In Section~\ref{sec:orthogonal}, we discuss a framework for constructing RKHS from classical orthogonal polynomials of weighted $L^2$-space. The main contribution is summarized in Proposition~\ref{prop:RKHS_d}. In Section~\ref{sec:Her}, we study the RKHS constructed using the Hermite polynomials. In this case, the resulting kernel is the well-known Mehler kernel \cite{mehler1866ueber}. We specify the regularity of the resulting RKHS in Corollary~\ref{lem:smooth}. In Section~\ref{sec:Lag}, we also study the RKHS constructed using the Laguerre polynomials. We show the boundedness of the resulting kernel, known as the Hille-Hardy kernel, in Lemma~\ref{lem:HHF}. In Section~\ref{sec:richness}, we show that the RKHS associated with the Mehler kernel is ``rich enough'' to approximate any continuous density function with Gaussian (or faster) decay-rate of arbitrary variance (see Corollary~\ref{corollary:rich} and Remark~\ref{remark_rich}). 

\item In Section~\ref{sec:KELR}, we introduce the kernel embedding approximation to the linear response statistics. To simplify the discussion and provide a concrete error bound, we only present results based on the Mehler kernel. The same overall conclusion holds with different constants in error bound and different weighted $L^2$-space, if the same technique is applied on the Hille-Hardy kernel. Based on the RKHS induced by the Mehler kernel, we consider the kernel embedding estimates of the target equilibrium distribution function of an ergodic It\^o diffusion. Consequently, we define the kernel embedding linear response operator as the FDT response operator associated with the kernel embedding density estimate. We discuss to which extent the well-posedness of the estimator can be conceived, excluding the data points that admit negative-value density estimates. Then, we summarize the consistency of the proposed estimator in the limit of $M\to\infty$ in Proposition~\ref{prop:consis}. Using the regularity of the functions in RKHS induced by the Mehler kernel, we specify sufficient conditions for the external forces to admit a well-defined FDT response operator for the class of target function in the RKHS (see Remark~\ref{rem:cond_forcing}). {   We also provide a simple criterion for choosing the parameter $M$ and other parameters for a well-defined estimator that avoids negative values density estimates.  In Proposition~\ref{prop:mixing}, we provide an error bound for the Monte-Carlo approximation under non-i.i.d data and discuss the implication of this result.}

\item In Section~\ref{sec:num}, we numerically examine the kernel embedding linear response. We choose two ergodic SDEs with known analytical equilibrium densities so that our approach can be directly validated. We will demonstrate the effectiveness of the Hermite and Laguerre polynomials in approximating densities with symmetric and non-symmetric decaying tails, respectively. Compared to the linear response estimator obtained via the classical KDE, we numerically find that the proposed estimator is computationally more efficient and more accurate.

\item In Section~\ref{sec:sum}, we close this paper with a summary and discussion on open problems and future research plans that stem from this study.
\end{itemize}

Some proofs are reported in the Appendices to improve the readability.

\section{Preliminary on the theory of linear response and RKHS} \label{sec:review}

In this section, we first review the linear response theory, which motivates the estimation of the equilibrium density of SDEs from  their time series. Then we include a short survey on RKHS and the relevant results. We refer readers to \cite{mag:05, muandet2017kernel, christmann2008support} for a more comprehensive discussion.

\subsection{Linear response theory}\label{sec:FDT}

Fluctuation-Dissipation Theory (FDT) is a mathematical framework for quantifying the linear response of a dynamical system subjected to small external forcing \cite{leith:75}. The linear response statistics, determined based on two-point equilibrium statistics of the unperturbed dynamics, provide estimates for the non-equilibrium properties. In statistical mechanics literature, FDT is a linear response approach \cite{evans2008statistical} which serves as a foundation for defining transport coefficients, e.g., viscosity, diffusion constant, heat conductivity, etc.

The review will be presented in the context of the $d$-dimensional (time-homogeneous) SDEs, also known as the It\^o diffusions \cite{Pav_book:14}. The unperturbed and perturbed systems of SDEs are written, respectively, as follows,
\BEA
\dot X &=& b(X) + \sigma(X)\dot W_t,  \label{eq:Ito}\\
\dot X^{\delta} &=& \left[b(X^{\delta}) + c(X^{\delta})\delta f(t)\right]+ \sigma(X^{\delta}) \dot U_{t}, \label{eq:Ito_per}
\EEA
where $W_{t}$ and $U_t$ are standard Wiener processes. In both the unperturbed and perturbed  systems in \eqref{eq:Ito}-\eqref{eq:Ito_per}, the vector field $b:\mathbb{R}^d\to\mathbb{R}^d$ denotes the drift and $\sigma:\mathbb{R}^d\to \mathbb{R}^d \times \mathbb{R}^d$ denotes the diffusion tensor.  In \eqref{eq:Ito_per}, an order-$\delta$ ($\delta\ll 1$) external perturbation is introduced, in the form of $f(x,t) = c(x)\delta f(t)$. 

We assume the unperturbed system governed by Eq. \eqref{eq:Ito} is ergodic with a positive equilibrium density $p_{eq}(\bm{x})$, where $\bm{x}\in \BR^{d}$, as the unique classical solution of the stationary Fokker-Planck equation (see Theorem 4.1 of \cite{Pav_book:14} for detailed conditions). We should clarify that while we refer to $p_{eq}$ as the equilibrium density (following the classical nomenclature of linear response theory in statistical mechanics \cite{Toda-Kubo-2}), the formulation in this paper is not restrictive to reversible processes. We also assume that the statistical properties associated with Eq. \eqref{eq:Ito_per} can be characterized by a perturbed density $p^{\delta}(\bm{x},t)$, which solves the corresponding Fokker-Planck equation with initial condition $p^{\delta}(\bm{x},0) = p_{eq}(\bm{x})$, that is, we initiate the perturbed dynamics in \eqref{eq:Ito_per} at the equilibrium state of the unperturbed dynamics in \eqref{eq:Ito}. In the following formulation, we do not assume that \eqref{eq:Ito_per} possesses a stationary distribution.

For any integrable observable $A(\cdot)$, we define the difference of two expectations
\begin{equation}\label{eq:full_resp}
\Delta \mathbb{E}[A](t): = \mathbb{E}_{p^{\delta}}[A(X^{\delta})](t) - \mathbb{E}_{p_{eq}}[A(X)],
\end{equation}
as the full response statistics. In Appendix A of \cite{HLZ:19b}, we have showed that computing the non-equilibrium statistics $ \mathbb{E}_{p^{\delta}}[A(X^{\delta})](t)$ in \eqref{eq:full_resp} requires extensive numerical simulations of the perturbed dynamics in \eqref{eq:Ito_per}. The linear response theory allows us to estimate the order-$\delta$ term of \eqref{eq:full_resp} by a convolution integral,
\begin{equation}\label{eq:lin_resp}
\Delta \mathbb{E}[A](t) =\int^{t}_{0}k_{A}(t-s)\delta f(s) \td s + \mathcal{O}(\delta^2),
\end{equation}
avoiding simulations of the unperturbed dynamics in \eqref{eq:Ito_per}. Specifically, FDT formulates the linear response operator, $k_{A}(t)$ in \eqref{eq:lin_resp}, as the following two-point statistics:
\begin{equation}\label{eq:lin_oper}
k_{A}(t): = \mathbb{E}_{p_{eq}} \left [A(X(t)) \otimes B (X(0)) \right], \quad B_{i}(X):= -\frac{\partial_{X_i}\left[c_{i}(X)p_{eq}(X)\right]}{p_{eq}(X)},
\end{equation}
where $B_{i}$ and $c_i$ denote the $i^{\text{th}}$ components of $B$ and $c$, respectively. In \cite{evans2008statistical}, the variable $B$ is called the conjugate variable to the external forcing. In general, if $A$ in Eq. \eqref{eq:full_resp} is an $m$-dimensional vector, then the linear response operator $k_{A}(t)$ is an $m$-by-$d$ matrix. The significance of FDT is that the response operator is defined without involving the perturbed density $p^{\delta}(x,t)$. Rather, it can be evaluated at equilibrium of the unperturbed dynamics. For a given $t \ge 0$, the value of $k_{A}(t)$ can be computed using a Monte-Carlo sum based on the time series of the unperturbed system \eqref{eq:Ito} at $p_{eq}$. For example, let $\left\{X_{n}=X(t_{n})\right\}_{n=1}^{N}$ be the time series generated at $p_{eq}$ with step length $\Delta t = t_{n+1}-t_{n}$, then for $t = s \Delta t$, the Monte-Carlo approximation can be written as 
\begin{equation}\label{eq:MC_approx}
k_{A}(t) \approx \frac{1}{N-s} \sum_{n= 1}^{N-s} A\left(X_{n+s}\right)\otimes B\left(X_n\right).
\end{equation}
In practice, the computation of \eqref{eq:MC_approx} can be done efficiently using the block averaging algorithm \cite{flyvbjerg1989error}.

In applications, the major issue comes from the conjugate variable $B$ in the linear response operator \eqref{eq:lin_oper}. Since $B$ depends on the explicit formula of $p_{eq}$, which may not be available, one cannot directly apply the Monte-Carlo approximation in \eqref{eq:MC_approx} given only the time series of $\{X_{n}\}$ at $p_{eq}$. Thus, it is natural to ask how to learn the density function $p_{eq}$ from the observed time series such that the conjugate variable $B$ can be determined via the estimated density, $\hat{p}_{eq}$. To guarantee a well-posed estimation problem, the following lemma provides conditions on observables $A$ and $B$ such that the two-point statistics $k_{A}(t)$ in \eqref{eq:lin_oper} is bounded $\forall t \geq 0$.

\begin{lemma}  \label{lem:bound_kA}
Let $X$ be the solution of \eqref{eq:Ito} initially at the equilibrium. Assume that $A$ and $B$ in \eqref{eq:lin_oper} have finite second moments with respect to $p_{eq}$, then the linear response operator $k_A(t)$ in \eqref{eq:lin_oper} is well-defined. In particular, we have
\begin{equation*} 
k_{A}(t) \leq \mathbb{E}_{p_{eq}} \left[A^{2}(X) \right]^{\frac{1}{2}} \otimes \mathbb{E}_{p_{eq}}\left[B^{2}(X)\right]^{\frac{1}{2}}, \quad \forall t \geq 0,
\end{equation*}
\end{lemma}
where the inequality and the square/square root operation are defined componentwise.
\begin{proof}
Let $\mathcal{L}$ be the generator of the It\^{o} diffusion \eqref{eq:Ito} \cite{Pav_book:14}, and $e^{t\mathcal{L}}$ be the corresponding semigroup. Let us denote the transition kernel of \eqref{eq:Ito} as,
\begin{equation*}
  p(\bm{x},t |  \bm{y},0) = e^{t \mathcal{L}^{*}} \delta(\bm{x}- \bm{y}), \quad t \geq 0,
\end{equation*}
where $\mathcal{L}^{*}$, acting on $\bm{x}$, is the adjoint operator of $\mathcal{L}$ in the standard $L^{2}(\BR^{d})$ space. The transition kernel $ p(\bm{x},t |  \bm{y},0)$, as the solution of the Fokker-Planck equation of \eqref{eq:Ito} with initial condition $\delta(\bm{x} - \bm{y})$ can be interpreted as a density function of $\bm{x}$. With these definitions, $k_{A}(t)$ in \eqref{eq:lin_oper} can be specified as a double integral \cite{Pav_book:14}
\begin{equation*}
  k_{A}(t) = \int_{\BR^{d}} \int_{\BR^{d}} A(\bm{x}) \otimes B(\bm{y}) p(\bm{x},t  |  \bm{y},0)  p_{eq}(\bm{y}) \td \bm{x} \td \bm{y}.
\end{equation*}
Using the Cauchy-Schwarz inequality, we have
\begin{equation}\label{eq:CS_1}
\begin{split} 
 k_{A}(t) &\leq  \left(\int_{\BR^{2d}} A^{2}(\bm{x}) p(\bm{x},t  |  \bm{y},0)  p_{eq}(\bm{y}) \td \bm{x} \td \bm{y} \right)^{\frac{1}{2}} \otimes \left(\int_{\BR^{2d}}  B^{2}(\bm{y}) p(\bm{x},t  |  \bm{y},0)  p_{eq}(\bm{y}) \td \bm{x} \td \bm{y} \right)^{\frac{1}{2}} \\
  & =  \left(\int_{\BR^{2d}}  \left(e^{tL}A^{2}(\bm{x}) \right) \delta(\bm{x} - \bm{y}) p_{eq}(\bm{y}) \td \bm{y}  \td \bm{x}\right)^{\frac{1}{2}} \otimes  \left( \int_{\BR^{2d}}   B^{2}(\bm{y}) p_{eq}(\bm{y})p(\bm{x},t  |  \bm{y},0)   \td \bm{x} \td \bm{y} \right)^{\frac{1}{2}} \\
  & =  \left(\int_{\BR^{d}}  \left(e^{tL}A^{2}(\bm{x}) \right) p_{eq}(\bm{x}) \td \bm{x}\right)^{\frac{1}{2}} \otimes  \left( \int_{\BR^{d}}   B^{2}(\bm{y}) p_{eq}(\bm{y}) \td \bm{y} \right)^{\frac{1}{2}} \\ 
  &=\mathbb{E}_{p_{eq}} \left[A^{2}(X) \right]^{\frac{1}{2}} \otimes \mathbb{E}_{p_{eq}}\left[B^{2}(X)\right]^{\frac{1}{2}}.
\end{split}
\end{equation}
In Eq. \eqref{eq:CS_1}, the inequality and identity are all defined componentwise.
\end{proof}

We should point out that the validity of the linear response theory for SDE has been discussed in \cite{hairer2010simple} under a very general setting. Here, the purpose of Lemma~\ref{lem:bound_kA} is to guarantee that the linear response operator in \eqref{eq:lin_oper}, which is the central object that we wish to approximate in this article, is bounded uniformly in time. The boundedness of the second moment of the observable $A$ with respect to $\peq$ is fulfilled in many applications. As for the conjugate function $B$, since it is related to the function $c(\cdot)$ in the external forcing through the formula in \eqref{eq:lin_oper}, the Lemma provides a condition for admissible external forcings. In Section~\ref{sec:KELR}, for a specific class of $p_{eq}$, we will provide a more concrete condition on $c(\cdot)$ such that $B$ has a bounded second moment with respect to $\peq$ (see Remark~\ref{rem:cond_forcing}).

It is worthwhile to mention that learning the distribution function from the observations is a classical problem in statistics and machine learning. In general, there are two types of approaches: parametric and nonparametric. The kernel embedding formulation, which is the focus of this paper, belongs to the nonparametric category. However, direct use of the kernel embedding formulation with radial-type kernels, e.g., $k(x,y)= h(\|x-y\|)$) to approximate $\peq$ is computationally expensive for this application. In particular, if the length of the training data is $N$, the Monte-Carlo integral in \eqref{eq:MC_approx} will require  computing the $\|\cdot\|$ norm in the kernel function $N(N-s)$-times, since the evaluation of $B$ (as defined in \eqref{eq:lin_oper}) on each sample point $X_i$ requires the computation of $h(\|X_i-X_j\|)$, for all $j=1,\ldots, N$. This computational cost is similar to the complexity of using the kernel density estimation with radial-type kernels. In addition, the kernel embedding expansion implemented with radial type kernels will require an inversion of an $N\times N$ matrix (i.e., solving a large linear system) that determines the expansion coefficients. To overcome these practical issues, we will consider the Mercer-type kernels, constructed using a set of orthogonal polynomials. With such kernels, the resulting kernel embedding approximation on $B$ (or effectively $\peq$) becomes a parametric model since the number of parameters is smaller than the size of the data. Furthermore, the expansion coefficients can be determined using only transpose operations (thus, avoiding a large matrix inversion), thanks to the orthogonality.

To facilitate a self-contained discussion, we now provide a quick review of the relevant background material on RKHS.

\subsection{Kernels and RKHS}

In this subsection, we review the notions of \textit{kernel}, \textit{feature space}, \textit{feature map}, and RKHS. Then we summarize a few useful results which will be used in the later proofs. All unlisted proofs can be found in Chapter 4 of \cite{christmann2008support}.  We will restrict our discussions to real-valued functions since it provides simpler notations and is adequate for our application.

\begin{definition}\label{def:kernel}
Let $X$ be a non-empty set. A function $k: X\times X \rightarrow \mathbb{R}$ is called a kernel on $X$ if there exists an $\mathbb{R}$-Hilbert space $H$ and a map $\Phi: X \rightarrow H$ such that $\forall x,y \in X$ we have
\begin{equation}\label{eq:def_kernel}
       k(x,y) = \left\langle \Phi(x), \Phi(y) \right\rangle_{H},
\end{equation}
where $\langle \cdot, \cdot \rangle_{H}$ denotes the inner product of $H$. We call $\Phi$ a feature map and $H$ a feature space of $k$.
\end{definition}

By Eq. \eqref{eq:def_kernel}, for any fixed $x_{1}, x_{2}, \dots,  x_{n} \in X$, the $n\times n$ Gram matrix 
\begin{equation}\label{eq:Gram}
K_{n}: = \left(k(x_{i}, x_{j}) \right)_{1\leq i,j\leq n} 
\end{equation}
is symmetric positive definite (SPD), that is, $\forall \bm{a} = (a_{1},a_{2}, \dots, a_{n})^{\top} \in \BR^{n}$, the bilinear form
\begin{equation*}
\begin{split}
\bm{a}^{\top} K_{n} \bm{a} &= \sum_{i =1}^{n} \sum_{j=1}^{n} a_{i}a_{j} k(x_{i},x_{j}) = \sum_{i =1}^{n} \sum_{j=1}^{n} a_{i}a_{j} \left \langle \Phi(x_{i}) , \Phi(x_{j}) \right\rangle_{H} \\
&= \left  \langle \sum_{i=1}^{n}a_{i}\Phi(x_{i}) ,\sum_{j=1}^{n} a_{j}\Phi(x_{j}) \right \rangle_{H} = \left \| \sum_{i=1}^{n}a_{i}\Phi(x_{i})\right\|_{H}^{2} \geq 0.
\end{split}
\end{equation*}
We say that a function $k: X\times X \rightarrow \BR$ is SPD if $\forall n\geq 1$, and $(x_{1}, x_{2}, \dots, x_{n}) \in X^{n}$ the corresponding Gram matrix \eqref{eq:Gram} is SPD. Here, we follow the convention in \cite{christmann2008support}. Our next lemma states that symmetric positive definiteness is not only a necessary, but also a sufficient condition for $k$ to be a kernel.
\begin{lemma}\label{lem:char_kernel}
A function $k:X\times X \rightarrow \BR$ is a kernel if and only if it is SPD.
\end{lemma}
The lemma above is useful for checking whether a given function is a kernel. With the concept of kernel, we now define the RKHS.

\begin{definition}\label{def:RKHS}
Let $X$ be a non-empty set and $\mathcal{H}$ be a $\BR$-Hilbert function space over $X$, i.e., a $\mathbb{R}$-Hilbert space of functions that maps  $X$ to $\mathbb{R}$. Then $\mathcal{H}$ is called an RKHS with kernel $k$, if $k(\cdot,x)\in \mathcal{H}$, $\forall x\in X$, and we have the reproducing property
\begin{equation}\label{eq:rep_prop}
f(x) = \left\langle f(\cdot), k(\cdot,x) \right\rangle_{\mathcal{H}}
\end{equation}
holds for all $f\in \mathcal{H}$ and all $x\in X$. In particular, we call such $k(\cdot, \cdot)$ a reproducing kernel of $\mathcal{H}$.
\end{definition}

From Definition~\ref{def:RKHS}, it seems that the reproducing kernel is a result of RKHS. However, there is a one-to-one correspondence between the RKHS and kernel (see Theorem 4.20 and 4.21 of \cite{christmann2008support}). In Section~\ref{sec:orthogonal}, where we construct the RKHS  for the estimation of linear response, we will first define a kernel, then build an RKHS to ``promote'' such kernel to a reproducing kernel.

In the rest of the Section, $\mathcal{H}$ always denotes as the RKHS with kernel $k$. The RKHS has the remarkable property that the norm convergence implies the pointwise convergence. More precisely, consider $f_{n} \rightarrow f$ in $\mathcal{H}$, that is, $\left\| f_{n} - f \right\|_{\mathcal{H}}\rightarrow 0$ as $n\rightarrow \infty$. Then, $\forall x\in X$, we have
\begin{equation}\label{eq:pointwise}
\left|(f_{n}-f)(x)\right| = \left|\left \langle f_{n}-f, k(\cdot, x) \right\rangle_{\mathcal{H}}\right| \leq \left \|f_{n}- f \right\|_{\mathcal{H}} \left\|k(\cdot, x) \right\|_{\mathcal{H}} \rightarrow 0,
\end{equation}
as $n\rightarrow \infty$. Eq. \eqref{eq:pointwise} also suggests that if $\|k(\cdot, x)\|_{\mathcal{H}}$ is bounded uniformly in $x \in X$, we will have the uniform convergence of $f_{n}$ to $f$. We arrive at the following lemma.

\begin{lemma} \label{lem:uniform}
Let $X$ be a topological space and $k$ be a kernel on $X$ with RKHS $\mathcal{H}$. If $k$ is bounded in the sense that
\begin{equation*}
       \| k\|_{\infty}: = \sup_{x\in X} \sqrt{k(x,x)} < \infty.
\end{equation*}
and $k(\cdot,x): X \rightarrow \BR$ is continuous $\forall x \in X$, then $\mathcal{H}\subset C_{b}(X)$ (space of bounded and continuous functions on $X$), and the inclusion $\id: \mathcal{H} \rightarrow C_{b}(X)$ is continuous with $\left \|\id: \mathcal{H} \rightarrow C_{b}(X)\right\| = \|k\|_{\infty}$.
\end{lemma}
Here, to see the connection between $\|k\|_{\infty}$ and $\|k(\cdot, x)\|_{\mathcal{H}}$, simply notice that $k(\cdot, x) \in \mathcal{H}$, and with the reproducing property \eqref{eq:rep_prop} we have
\begin{equation} \nonumber 
\left\|k(\cdot, x) \right\|_{\mathcal{H}}^{2} = \left \langle k(\cdot, x), k(\cdot, x) \right \rangle_{\mathcal{H}} = k(x,x).
\end{equation}
Thus, $\|k\|_{\infty}$ is the upper bound of $\|k(\cdot, x)\|_{\mathcal{H}}$, and $\forall f\in \mathcal{H}$,
\begin{equation}\label{eq:decay_rate}
 \left| f(x) \right| = \left| \left \langle f , k(\cdot, x) \right\rangle_{\mathcal{H}} \right|  \leq \left\|f\right\|_{\mathcal{H}} \left\| k(\cdot, x) \right\|_{\mathcal{H}} = \left\|f \right\|_{\mathcal{H}} k^{\frac{1}{2}}(x,x), \quad \forall x \in X,
\end{equation}
that is, $f(x)$ has the same decay rate as $k^{\frac{1}{2}}(x,x)$. In this paper, we focus on the RKHS with a bounded kernel.

As a subspace of $C_{b}(X)$, it is natural to ask whether the RKHS $\mathcal{H}$ is dense in the Banach space $C_{b}(X)$ equipped with the uniform norm. The density of $\mathcal{H}$ in $C_{b}(X)$ is equivalent to the $c$-universality \cite{sriperumbudur2011universality} of the corresponding kernel $k$, which has been discussed by Steinwart \cite{steinwart2001influence} for compact $X$. In this paper, we are interested in the case where $X$ is non-compact, e.g., $X = \BR^{d}$, and the target $f$ is a continuous density function vanishing at infinity. For a locally compact Hausdorff (LCH) space $X$,  let $C_{0}(X)$ denote the space of all continuous functions on $X$ that vanish at infinity, that is, $\forall \delta >0$, the set $ \left\{ \bm{x} \in X \;|\; |f(\bm{x})| \geq \delta \right\} $ is compact $ \forall f\in C_{0}(X)$. $C_{0}(X)$, like $C_{b}(X)$, is a Banach space with respect to the infinite-norm $\|\cdot\|_{\infty}$. As an analog of the $c$-universal kernel, the concept of $c_{0}$-universal was introduced by Sriperumbudur et al. in \cite{sriperumbudur2011universality}.

\begin{definition} \label{def:c_0}
($c_{0}$-universal) Let $X$ be an LCH space and let $k$ be a bounded kernel on $X\times X$ and $k(\cdot, x) \in C_{0}(X)$, $\forall x\in X$. The kernel $k$ is said to be $c_{0}$-universal if the RKHS, $\mathcal{H}$, induced by $k$ is dense in $C_{0}(X)$ with respect to the uniform norm, that is, $\forall g \in C_{0}(X)$ and $\forall \epsilon >0$, there exists a function $\hat{g} \in \mathcal{H}$ such that $\|g - \hat{g}\|_{\infty} < \epsilon$.
\end{definition}

A series of characterizations of $c_{0}$-universality for different types of kernels has been developed in \cite{sriperumbudur2010relation,sriperumbudur2011universality} based on the Hahn-Banach theorem and Riesz representation theorem. When $X = \BR^{d}$, a typical example of $c_{0}$-universal kernel is the Gaussian kernel $k(\bm{x},\bm{y}) = \exp(-\theta \|\bm{x} - \bm{y}\|^{2})$, $\bm{x}, \bm{y} \in \BR^{d}$, for some $\theta>0$. To facilitate the applications in this paper (see Section~\ref{sec:orthogonal}), we generalize the concept of $c_{0}$-universality to a weighted $C_{0}$-space.

\begin{lemma} \label{lem:weighted}
Let $X$ be an LCH space and $q$ be a bounded positive continuous function on $X$. Then we have the following results.
\begin{enumerate}
\item  The set of functions
\begin{equation*}
C_{0}(X, q^{-1}): = \left\{ f \in C(X) \; \big| \;  fq^{-1} \in C_{0}(X)    \right\}
\end{equation*}
defines a vector space.
\item The map $\|\cdot \|_{C_{0}(q^{-1})}: C_{0}(X, q^{-1}) \rightarrow \BR$ defined as
\begin{equation*}
\left \|f  \right\|_{C_{0}(q^{-1})} : = \left\| fq^{-1} \right\|_{\infty}, \quad f\in C_{0}(X, q^{-1}),
\end{equation*}
is a norm. Moreover, $C_{0}(X, q^{-1})$, equipped with the norm $\|\cdot \|_{C_{0}(q^{-1})}$, is a Banach space.
\item The Banach spaces $C_{0}(X, q^{-1})$ and $C_{0}(X)$ are isometrically isomorphic.
\end{enumerate}
\end{lemma}

\begin{proof}
It is straightforward to check that $C_{0}(X, q^{-1})$ is a normed vector space with respect to the norm $\|\cdot\|_{C_{0}(q^{-1})}$. To see $C_{0}(X, q^{-1})$ is closed under the topology induced by $\|\cdot\|_{C_{0}(q^{-1})}$, we introduce the linear bijection $ Q: C_{0}(X) \rightarrow C_{0}(X, q^{-1})$ defined as  $Q(g):=gq$ for every $g\in C_0(X)$. The bijection $Q$ is norm-preserving in the sense that $\forall g \in C_{0}(X)$,
\begin{equation*}
\| g\|_{\infty} = \| gq \|_{C_{0}(q^{-1})} = \| Qg \|_{C_{0}(q^{-1})}.
\end{equation*}
Therefore, for any Cauchy sequence $\left\{f_{n}\right\}$ in $C_{0}(X, q^{-1})$, $\left\{g_{n}: = Q^{-1}{f_{n}}\right\}$ defines a Cauchy sequence in Banach space $C_{0}(X)$. Let $g_{n} \rightarrow g^{*}$ in $C_{0}(X)$, then $f^{*}: = Qg^{*} \in C_{0}(X, q^{-1})$ and $f_{n}\rightarrow f^{*}$ in $C_{0}(X, q^{-1})$. With $C_{0}(X, q^{-1})$ being a Banach space, the operator $Q$ becomes an isometrical isomorphism.
\end{proof}
In practice, we use the weight function $q$ to characterize the decay rate of continuous functions. For example, take $X = \BR^{d}$, and $q \propto \exp(-\theta \|\bm{x} \|^{2})$ for some $\theta >0$, then the functions in $C_{0}(\BR^{d}, q^{-1})$ are continuous with a Gaussian decay rate. Motivated by the decay rate \eqref{eq:decay_rate} of functions in $\mathcal{H}$, the following lemma provides conditions for $\mathcal{H}$ to be dense in $C_{0}(X, q^{-1})$.

\begin{lemma} \label{lem:weight_c0}
(weighted $c_{0}$-universal) 
Let $X$ and $q$ be the same as in Lemma~\ref{lem:weighted}, and the kernel $k$, satisfying $k(\cdot, x) \in C_{0}(X, q^{-1})$, $\forall x\in X$. Then, $\tilde{k}(x,y):= q^{-1}(x)k(x,y)q^{-1}(y)$ defines a kernel on $X$, and the RKHS $\mathcal{H}$ induced by $k$ is dense in $C_{0}(X, q^{-1})$ if and only if the kernel $\tilde{k}$ is $c_{0}$-universal.
\end{lemma}

\begin{proof}
To begin with, by Lemma~\ref{lem:char_kernel}, $\tilde{k}$ defines a kernel, and $\tilde{k}(\cdot, x) = q^{-1}(\cdot)k(\cdot, x)q^{-1}(x) \in C_{0}(X)$, $\forall x\in X$ since $k(\cdot, x) \in C_{0}(X, q^{-1})$, $\forall x\in X$.

By Definition~\ref{def:c_0}, it is enough to show that $\mathcal{H}$ being dense in $C_{0}(X, q^{-1})$ is equivalent to $\tilde{\mathcal{H}}$, the RKHS induced by $\tilde{k}$, being dense in $C_{0}(X)$. Recall that, by Lemma~\ref{lem:weighted}, the Banach spaces $C_{0}(X, q^{-1})$ and $C_{0}(X)$ are isometrically isomorphic with $ Q: C_{0}(X) \rightarrow C_{0}(X, q^{-1})$, defined as $Q(g)=gq$ for every $g\in C_0(X)$, being the isomorphism. Following the same idea, we take
\begin{equation*}
\tilde{\mathcal{H}}:= \left\{fq^{-1} \; \big| \; f \in \mathcal{H} \right\}
\end{equation*}
with the inner product
\begin{equation*}
\langle g_{1}, g_{2} \rangle_{\tilde{\mathcal{H}}} := \langle g_{1}q, g_{2}q\rangle_{\mathcal{H}}, \quad \forall g_{1},g_{2} \in \tilde{\mathcal{H}}.
\end{equation*}
Since $\mathcal{H}$ is a Hilbert space, $\tilde{\mathcal{H}}$ equipped with the inner product $\langle \cdot, \cdot \rangle_{\tilde{\mathcal{H}}}$ is also a Hilbert space. With $\tilde{k}(\cdot, x) = q^{-1}(\cdot)k(\cdot, x)q^{-1}(x)$, and $k(\cdot, x) \in \mathcal{H}$, $\forall x\in X$, we have $\tilde{k}(\cdot, x) \in \tilde{\mathcal{H}}$, $\forall x\in X$.

In terms of the reproducing property, $\forall g = fq^{-1} \in \tilde{\mathcal{H}}$ with $f \in \mathcal{H}$, we have
\begin{equation*}
\langle g, \tilde{k}(\cdot, x) \rangle_{\tilde{\mathcal{H}}} = \langle f, k(\cdot, x) q^{-1}(x) \rangle_{\mathcal{H}} = \langle f, k(\cdot, x)  \rangle_{\mathcal{H}}q^{-1}(x) = f(x)q^{-1}(x) = g(x), \quad \forall x\in X.
\end{equation*}
Thus, $\tilde{\mathcal{H}}$ indeed is the RKHS induced by $\tilde{k}$ satisfying $Q \left(\tilde{\mathcal{H}} \cap C_{0}(X)\right)$ = $\mathcal{H} \cap C_{0}(X, q^{-1})$. Finally, by the fact that $Q$ defines an isometrical isomorphism between $C_{0}(X)$ and $C_{0}(X, q^{-1})$, we reach the equivalence. 
\end{proof}

Our next lemma characterizes how the differentiability of a kernel is inherited by the functions of its RKHS. In particular, we take $X \subset \mathbb{R}^{d}$ to be an open subset, and introduce the multi-index notation $\vec{m} = (m_{1}, m_{2}, \dots, m_{d})$ with $m_{i}$ being nonnegative integers. Then, we say that the kernel $k$ is $M$-times continuously differentiable if $\partial^{\vec{m}, \vec{m}}k: X\times X \rightarrow \mathbb{R}$ exist and are continuous for all multi-indices $\vec{m}$ with $\|\vec{m}\|_{1}:= \sum m_{i} \leq M$. Recall that
\begin{equation*}
\partial^{\vec{m}, \vec{m}}k(\bm{x}, \bm{y}): = \frac{\partial^{\vec{m}} }{\partial \bm{x}^{\vec{m}}} \frac{\partial^{\vec{m}} }{\partial \bm{y}^{\vec{m}}} k(\bm{x}, \bm{y}), \quad \bm{x}, \bm{y} \in \BR^{d}.
\end{equation*}

\begin{lemma} \label{lem:diff}
Let $X$ be an open subset of $\mathbb{R}^{d}$, and kernel $k$ be an $M$-times continuously differentiable kernel. Then every $f\in \mathcal{H}$ is $M$-times continuously differentiable in $X$, and $\forall \|\vec{m}\|_{1}\leq M$, we have
\begin{equation*}
\left | \partial^{\vec{m}}f(\bm{x})  \right| \leq \left\|f \right\|_{\mathcal{H}} \cdot \left ( \partial^{\vec{m}, \vec{m}}k(\bm{x},\bm{x})\right)^{\frac{1}{2}}, \quad \forall \bm{x}\in \mathbb{R}^{d}.
\end{equation*}
\end{lemma}

From Lemma~\ref{lem:uniform} and \ref{lem:diff} we have learned the importance of constructing reproducing kernel $k$ with certain ``nice'' properties. In practice, we can construct a bounded kernel using \eqref{eq:def_kernel} based on a feature map $\Phi: X\rightarrow H$, where $H$ is a  Hilbert space (e.g., $\ell_{2}$-space). Then, we define the corresponding RKHS $\mathcal{H}$ as a subspace of $C_{0}(X)$ such that $\mathcal{H}$ yields the reproducing property. The coming section will follow this procedure of constructing RKHS, using a feature map induced by orthogonal polynomials.

\section{From orthogonal polynomials to RKHS} \label{sec:orthogonal}

In Section \ref{sec:review}, we have discussed the density estimation problem arising from the linear response theory. Given observed time series, our goal will be to approximate the target equilibrium density function $\peq$ with appropriate RKHS functions. In practice, such RKHS should be well selected such that:
\begin{enumerate}
\item The corresponding reproducing kernel $k$ has all good properties in Lemma~\ref{lem:uniform} and \ref{lem:diff} so that we can derive convergence results of the estimates. We will show examples of RKHSs constructed from Hermite and Laguerre polynomials (see Section~\ref{sec:general}).
\item The RKHS is rich enough in the sense of Lemma~\ref{lem:weight_c0} { such that one can estimate}  $p_{eq}$ that has a certain decaying rate in arbitrary precision (see Section~\ref{sec:richness}).
\end{enumerate}

\subsection{Constructing RKHS via orthogonal polynomials} \label{sec:general}

Inspired by the Mercer's theorem \cite{christmann2008support} and reproducing kernel weighted Hilbert space used in \cite{berry2017correcting, JH:18,JH:19,HLZ:19b}, we consider the orthonormal polynomials with respect to  $L^2(\mathbb{R}^{d}, \bm{W})$, to construct our kernel and RKHS. Here, $L^2(\mathbb{R}^{d}, \bm{W})$ denotes the product space $\prod_{i=1}^{d} L^2(\mathbb{R}, {W})$, and $\bm{W}(\bm{x})= \prod_{i=1}^{d}W(x_{i})$ for $\bm{x} = (x_{1},x_{2},\dots, x_{d}) \in \BR^{d}$. 

For a more concise discussion, we begin our analyses with the class of one-dimensional weight functions, satisfying conditions specified in the following lemma. In the lemma, the orthonormal polynomials are obtained via applying Gram-Schmidt process to $\{1,x,$ $x^2,\dots\}$ based on the inner product in $L^{2}(\mathbb{R}, W)$ and $p_{n}$ always denotes a polynomial of degree $n$.

\begin{lemma}\label{lem:bound}
Let $W = e^{-2Q}$, where $Q: \mathbb{R} \rightarrow \mathbb{R}$ is an even $C^{2}$-function. We assume that $Q>0$  on $(0, \infty)$, and
 there exist real numbers $B\geq A>1$, such that
\begin{equation}
A \leq \frac{\left(xQ'(x)\right)'}{Q'(x)}\leq B, \quad \forall x \in (0, \infty).\nonumber
\end{equation}
Then the orthonormal polynomials $\{p_{n}\}$  in $L^{2}(\mathbb{R}, W)$ satisfy
\begin{equation}\label{eq:bound}
C_{1}n^{\frac{1}{6}-\frac{1}{2A}}\leq \sup_{x\in \mathbb{R}} \left| p_{n}(x) \right| W^{\frac{1}{2}}(x)  \leq C_{2}n^{\frac{1}{6}-\frac{1}{2B}}, \quad  n =0,1,\dots,
\end{equation}
where $C_{1}$ and $C_{2}$ are positive constants independent of $n$.
\end{lemma}

\begin{proof}
Eq. \eqref{eq:bound} is a combination of Corollary 1.4 and Lemma 5.2 in \cite{levin:92}.
\end{proof}

A typical type of functions satisfying the conditions in Lemma~\ref{lem:bound} is $Q_{m}(x) = |x|^{m}$ for $m>1$, and we have
\begin{equation*}
\frac{\left(xQ_{m}'(x)\right)'}{Q_{m}'(x)} = m.
\end{equation*}
Let $W_{m} = e^{-2Q_{m}}$, and \eqref{eq:bound} leads to
\begin{equation*}
\sup_{x\in \mathbb{R}} \left| p_{n}(x) \right| W^{\frac{1}{2}}_{m}(x) \sim n^{\frac{1}{6}-\frac{1}{2m}}.
\end{equation*}

With the control of the $L^{\infty}$-norm of the orthonormal polynomials in \eqref{eq:bound}, the following lemma defines the bounded kernel we need to build our RKHS.

\begin{lemma}\label{lem:kernel}
Let $W$ and $p_{n}(x)$ be as in Lemma~\ref{lem:bound}. Given a sequence of monotonically decreasing positive real numbers $\{\lambda_{n}\}_{n=0}^{\infty}$ satisfying
\begin{equation}\label{eq:condlam}
\sum_{n=0}^{\infty} \lambda_{n} n^{\frac{1}{3} - \frac{1}{B}} < \infty,
\end{equation}
then, for $\beta \geq \frac{1}{2}$, the bivariate function $k_{\beta}(\cdot,\cdot): \mathbb{R} \times \mathbb{R} \rightarrow \mathbb{R}$ defined by,
\begin{equation}\label{eq:kernel_1}
k_{\beta}(x,y): = \sum_{n=0}^{\infty} \lambda_{n} p_{n}(x)p_{n}(y) W^{\beta}(x)W^{\beta}(y),
\end{equation}
is a well-defined bounded continuous function. Moreover, $k_{\beta}$ is a kernel.
\end{lemma}

\begin{proof}                                             Notice that by the uniform bound \eqref{eq:bound} in Lemma~\ref{lem:bound}, we have
\begin{equation*}
   \left| \lambda_{n} p_{n}(x)p_{n}(y) W^{\beta}(x)W^{\beta}(y) \right| \leq C^{2}_{2} \lambda_{n} n^{\frac{1}{3}-\frac{1}{B}} W^{\beta-\frac{1}{2}}(x)W^{\beta-\frac{1}{2}}(y),
\end{equation*}
and combining with the condition \eqref{eq:condlam}, we arrive at the uniform convergence of the summation in \eqref{eq:kernel_1}. Thus, $k_{\beta}(x,y)$ in \eqref{eq:kernel_1} is a well-defined continuous function on $\mathbb{R}^{2}$ with a decay rate,
\begin{equation} \label{eq:decay}
\left|k_{\beta}(x,y) \right| \leq C_{3} W^{\beta-\frac{1}{2}}(x)W^{\beta-\frac{1}{2}}(y), \quad C_{3}:= C_{2}^{2} \sum_{n=0}^{\infty} \lambda_{n} n^{\frac{1}{3}-\frac{1}{B}}.
\end{equation}
To show that $k_\beta$ is a kernel, we define the \textit{feature map} $\Phi_{\beta}: \mathbb{R} \rightarrow \ell_{2}$ as
\begin{equation*}
\Phi_{\beta}(x):= \left( \sqrt{\lambda_{0}}  p_{0}(x) W^{\beta}(x),  \sqrt{\lambda_{1}}  p_{1}(x) W^{\beta}(x), \dots, \sqrt{\lambda_{n}}  p_{n}(x) W^{\beta}(x),\dots \right), \quad x\in \BR.
\end{equation*}
With $k_{\beta}(x,y) = \langle \Phi_{\beta}(x) , \Phi_{\beta}(y) \rangle_{\ell_{2}}$, and by Definition~\ref{def:kernel}, $k_\beta$ is a kernel.
\end{proof}

To generalize Lemma~\ref{lem:kernel} to the $d$-dimensional case, consider the orthonormal polynomial in $L^{2}(\BR^{d}, \bm{W})$ of the form
\begin{equation*}
p_{\vec{m}}(\bm{x}) := \prod_{i=1}^{d} p_{m_{i}}(x_{i}), \quad \bm{x}\in \BR^{d},
\end{equation*}
where $\vec{m} = (m_{1}, m_{2}, \dots, m_{d})$ is a multi-index and $\{p_{n}\}$ are the orthonormal polynomials  in $L^{2}(\BR, W)$. Following Lemma~\ref{lem:kernel}, we define
\begin{equation}\label{eq:kernel_d}
k_{\beta}(\bm{x}, \bm{y}): =  \sum_{\vec{m}\geq 0} \lambda_{\vec{m}} p_{\vec{m}}(\bm{x}) p_{\vec{m}}(\bm{y}) \bm{W}^{\beta}(\bm{x})\bm{W}^{\beta}(\bm{y}),\quad \bm{x}, \bm{y} \in \BR^{d}, \; \lambda_{\vec{m}}: = \prod_{i=1}^{d} \lambda_{m_{i}}, 
\end{equation}
for $\beta \geq \frac{1}{2}$, as the $d$-dimensional generalization of the kernel $k_{\beta}$ in Eq.~\eqref{eq:kernel_1}. Here, $\lambda_{n}$ satisfies the condition in Lemma~\ref{lem:kernel}. The function $k_{\beta}(\bm{x}, \bm{y})$ in \eqref{eq:kernel_d} is well-defined since
\begin{equation*}
\sum_{\vec{m}\geq 0} \lambda_{\vec{m}} p_{\vec{m}}(\bm{x})p_{\vec{m}}(\bm{y}) \bm{W}^{\beta}(\bm{x})\bm{W}^{\beta}(\bm{y}) = \prod_{i=1}^{d} \left( \sum_{m_{i}=0}^{\infty} \lambda_{m_{i}} p_{m_{i}}(x_{i})p_{m_{i}}(y_{i}) W^{\beta}(x_{i})W^{\beta}(y_{i}) \right)< \infty.
\end{equation*}
Moreover, similar to \eqref{eq:decay}, $k_{\beta}(\bm{x}, \bm{y})$ yields the following decay rate
\begin{equation}\label{eq:decay_d}
 \left |   k_{\beta}(\bm{x}, \bm{y}) \right | \leq C_{3}^{d} \bm{W}^{\beta-\frac{1}{2}}(\bm x) \bm{W}^{\beta-\frac{1}{2}}(\bm{y}),
\end{equation} 
where $C_{3}$ is the same constant as in \eqref{eq:decay}. 

Formally, one can always define a kernel by the infinite sum in \eqref{eq:kernel_1}. In Lemma~\ref{lem:kernel}, the uniform convergence of the series in \eqref{eq:kernel_1} is proved based on the asymptotic behavior of a class of orthogonal polynomials (see Lemma~\ref{lem:bound}). In practice, one can also rely on existing identities to show the boundedness of the kernel defined in \eqref{eq:kernel_1} (see Section~\ref{sec:Lag} below).

Our next task is to define the RKHS, denoted by $\mathcal{H}_{\beta}$, such that the bounded kernel $k_{\beta}$ in \eqref{eq:kernel_d} is the reproducing kernel of $\mathcal{H}_{\beta}$. A crucial observation is that, by Definition~\ref{def:RKHS}, a necessary condition for $\mathcal{H}_{\beta}$ is $k_{\beta}(\cdot, \bm{x}) \in \mathcal{H}_{\beta}$, $\forall \bm{x}\in \BR^{d}$. Thus, we shall first study the function space containing $\left\{k_{\beta}(\cdot, \bm{x}), \; \bm{x} \in \BR^{d}\right\}$. We have the following lemma.

\begin{lemma} 
Let $W(x)$, $p_{n}(x)$, and $\{\lambda_{n}\}$ be as in Lemma~\ref{lem:kernel}. The kernel $k_{\beta}$ defined by \eqref{eq:kernel_d} satisfies $\forall \bm{x} \in \BR^{d}$, $k_{\beta}(\cdot, \bm{x}) \in L^{2}(\mathbb{R}^{d}, \bm{W}^{1-2\beta})\cap C_{0}(\mathbb{R}^{d})$. 
\end{lemma}

\begin{proof}
First notice that $\forall \bm{x}\in \BR^{d}$, $k_{\beta}(\cdot, \bm{x}) \in C_{0}(\BR^{d})$ is a direct result of the decay rate estimate  \eqref{eq:decay_d}. Introduce
\begin{equation*}
    \Psi_{\beta,\vec{m}}(\bm{x}):=  p_{\vec{m}}(\bm{x}) \bm{W} ^{\beta}(\bm{x}),
\end{equation*}
such that $\left\{\Psi_{\beta,\vec{m}}\right\}$ forms an orthonormal basis on $L^{2}(\BR^{d}, \bm{W}^{1-2\beta})$. We rewrite the kernel $k_{\beta}(\bm{x}, \bm{y})$ in \eqref{eq:kernel_d} as
\begin{equation}\label{eq:kernel_Mercer}
      k_{\beta}(\bm{x}, \bm{y}) = \sum_{\vec{m} \geq 0} \lambda_{\vec{m}} \Psi_{\beta, \vec{m}}(\bm{x}) \Psi_{\beta,\vec{m}}(\bm{y}).
\end{equation}
In particular, by the orthogonality of $\left\{ \Psi_{\beta, \vec{m}} \right\}$, we have 
\begin{equation*}
\int_{\BR^{d}} k_{\beta}^{2}(\bm{y}, \bm{x}) \bm{W}^{1-2\beta}(\bm{y}) \td \bm{y} = \sum_{\vec{m} \geq 0} \lambda^{2}_{\vec{m}} \Psi^{2}_{\beta,\vec{m}}(\bm{x})  \leq \lambda_{\vec{0}} k_{\beta}(\bm{x}, \bm{x}) \leq  \lambda_{0}^{d}\left \|k_{\beta} \right\|_{\infty}^{2},  \quad \forall \bm{x} \in \BR^{d}.
\end{equation*}
Thus, $\left\{k_{\beta}(\cdot, \bm{x}), \; \bm{x} \in \BR^{d}\right\} \subset L^{2}(\mathbb{R}^{d}, \bm{W}^{1-2\beta})\cap C_{0}(\mathbb{R}^{d})$.
\end{proof}

We shall emphasize that $L^{2}(\mathbb{R}^{d}, \bm{W}^{1-2\beta})\cap C_{0}(\mathbb{R}^{d})$ is the set consisting of all continuous \textit{functions} vanishing at the infinity which are $\bm{W}^{1-2\beta}$-weighted $L^{2}$-integrable, while each element in $L^{2}(\mathbb{R}^{d}, \bm{W}^{1-2\beta})$ represents an equivalent class of functions due to the difference of their topologies. With the topology induced by the weighted $L^{2}$-norm $\| \cdot \|_{L^{2}(\bm{W}^{1-2\beta})}$, the expansion formula for $ f \in L^{2}(\mathbb{R}^{d}, \bm{W}^{1-2\beta})$ given by 
\begin{equation} \label{eq:expansion}
 f  = \sum_{\vec{m}\geq 0} \hat{f}_{\vec m} \Psi_{\beta,\vec{m}}, \quad \hat{f}_{\vec{m}} := \int_{\mathbb{R}^{d}} f(\bm x) \Psi_{\beta,\vec{m}}(\bm x) \bm{W}^{1-2\beta}(\bm x) \td \bm x,\end{equation}
is valid in the sense that
\begin{equation}
\lim_{M \rightarrow \infty}\left \|  f  -  f_{M}\right\|_{L^{2}(\bm{W}^{1-2\beta})} = 0, \quad f_{M}:= \sum_{\|\vec{m}\|_{1} \leq M} \hat{f}_{\vec m} \Psi_{\beta,\vec{m}},\nonumber
\end{equation}
which is relatively weak since most of the properties of $f_{M}$ no longer exist after passing thought the limit. The following proposition, as the main result of Section~\ref{sec:general}, defines the RKHS, $\mathcal{H}_{\beta}\subset L^{2}(\mathbb{R}^{d}, \bm{W}^{1-2\beta})\cap C_{0}(\mathbb{R}^{d})$, corresponding to the kernel $k_{\beta}$ in \eqref{eq:kernel_d}, such that $\forall f \in \mathcal{H}_{\beta}$, the expansion \eqref{eq:expansion} holds pointwisely and $f_{M}$ converges uniformly to $f$ in $C_{0}(\mathbb{R}^{d})$.

\begin{proposition}\label{prop:RKHS_d}
Let $W(x)$, $p_{n}(x)$, and $\{\lambda_{n}\}$ be as in Lemma~\ref{lem:kernel}. Then, for any fixed $\beta \geq \frac{1}{2}$, we have the following results.
\begin{enumerate}[\roman{enumi}]
\item For any sequence $\left\{ \hat{f}_{\vec{m}} \right\} \in \ell_{2}$ satisfying 
\begin{equation}\label{eq:RKHS_cond1}
         \sum_{\vec{m}\geq 0}\frac{\hat{f}^2_{\vec m}}{\lambda_{\vec m}} < \infty,
\end{equation}
where $\lambda_{\vec{m}}$ is defined in \eqref{eq:kernel_d}, the sequence of functions
\begin{equation*}
f_{n}: = \sum_{\|\vec{m}\|_{1}\leq n} \hat{f}_{\vec{m}} \Psi_{\beta,\vec{m}}, \quad n\geq 0,
\end{equation*}
converge uniformly in $C_{0}(\BR^{d})$. Moreover, the limit, denoted as $f^{*}$, satisfies $f^{*} \in L^{2}(\mathbb{R}^{d}, \bm{W}^{1-2\beta})\cap C_{0}(\mathbb{R}^{d})$.
\item The function space
\begin{equation}\label{eq:RKHS}
\mathcal{H}_{\beta}: = \left\{ f = \sum_{\vec{m}\geq 0} \hat{f}_{\vec{m}} \Psi_{\beta,\vec{m}} \;\Big| \;   \sum_{\vec{m}\geq 0}\frac{\hat{f}^2_{\vec m}}{\lambda_{\vec m}} < \infty \right\},
\end{equation}
is a well-defined subspace of $L^{2}(\mathbb{R}^{d}, \bm{W}^{1-2\beta})\cap C_{0}(\mathbb{R}^{d})$. Further, define the map $\left \langle \cdot, \cdot \right\rangle : \mathcal{H}_{\beta} \times \mathcal{H}_{\beta} \rightarrow \mathbb{R}$ as
\begin{equation*}
       \left \langle f, g  \right\rangle : = \sum_{\vec{m}\geq 0}\frac{\hat{f}_{\vec m}\hat{g}_{\vec{m}}}{\lambda_{\vec m}}, \quad f = \sum_{\vec{m}\geq 0} \hat{f}_{\vec{m}} \Psi_{\beta, \vec{m}}, \;g= \sum_{\vec{m}\geq 0} \hat{g}_{\vec{m}} \Psi_{\beta, \vec{m}} \in \mathcal{H}_{\beta}.
\end{equation*}
Then $\left \langle \cdot, \cdot \right\rangle $ defines an inner product, and $\mathcal{H}_{\beta}$, equipped with the inner product $\left \langle \cdot, \cdot \right\rangle$, is a Hilbert space.
\item $\mathcal{H}_{\beta}$ is the RKHS with reproducing kernel $k_{\beta}$ in \eqref{eq:kernel_d}.
\end{enumerate}
\end{proposition}

See Appendix~\ref{app:A} for the proofs. With $k_{\beta}$ being bounded, Lemma~\ref{lem:uniform} implies that the norm convergence in $\mathcal{H}_{\beta}$ is a sufficient condition of the uniform convergence in $C_{0}(\BR^{d})$. In particular, $\forall f \in \mathcal{H}_{\beta}$, the expansion formula \eqref{eq:RKHS} converges uniformly, that is, there is a one-to-one correspondence between the \textit{function} $f \in \mathcal{H}_{\beta}$ and the sequence of expansion coefficients $\left(\hat{f}_{\vec{m}}\right) \in \ell_{2}$ satisfying \eqref{eq:RKHS_cond1}. Since the condition \eqref{eq:RKHS_cond1} is independent of $\beta$, then the class of RKHS $\{\mathcal{H}_{\beta}\}$ defined in Proposition~\ref{prop:RKHS_d} are isometrically isomorphic to each other. In particular, $\forall \beta_{1}, \beta_{2} >\frac{1}{2}$, the linear map $I_{\beta_{1}, \beta_{2}}: \mathcal{H}_{\beta_{1}} \rightarrow \mathcal{H}_{\beta_{2}}$ given by
\begin{equation*}
 I_{\beta_{1}, \beta_{2}}f =  f \cdot  \bm{W}^{\beta_{2} - \beta_{1}}, \quad \forall f\in \mathcal{H}_{\beta_{1}},
 \end{equation*}
defines the isometrical isomorphism between $\mathcal{H}_{\beta_{1}}$ and $\mathcal{H}_{\beta_{2}}$ with $ I_{\beta_{1}, \beta_{2}}^{-1} = I_{\beta_{2}, \beta_{1}}$.

To connect the kernel $k_{\beta}$ in \eqref{eq:kernel_d} with Mercer's theorem, we define an integral operator $T_{k_{\beta}}: L^{2}(\BR^{d}, \bm{W}^{1-2\beta}) \rightarrow L^{2}(\BR^{d}, \bm{W}^{1-2\beta})$ as follows
\begin{equation}\label{eq:int_oper}
(T_{k_{\beta}}f)(\bm{x}): = \int_{\BR^{d}} k_{\beta}(\bm{y}, \bm{x}) f(\bm{y}) \bm{W}^{1- 2\beta}(\bm{y}) \td \bm{y}.
\end{equation}
Set $f = \Psi_{\beta,\vec{n}}$ in \eqref{eq:int_oper}, and we obtain
\begin{equation*}
\left( T_{k_{\beta}}\Psi_{\beta, \vec{n}} \right)(\bm{x}) = \sum_{\vec{m} \geq 0} \lambda_{\vec{m}} \Psi_{\beta, \vec{m}}(\bm{x})  \int_{\BR^{d}} \Psi_{\beta, \vec{m}}(\bm{y}) \Psi_{\beta, \vec{n}}(\bm{y}) \bm{W}^{1-2\beta}(\bm{y})\td \bm{y} =  \lambda_{\vec{n}}\Psi_{\beta, \vec{n}}(\bm{x}),
\end{equation*}
where we have used the fact that the convergence in \eqref{eq:kernel_d} is uniform to switch the order of integration and summation. Thus, $\lambda_{\vec{m}}$, and $\Psi_{\vec{m}}$ correspond to the eigenvalue and eigenfunction of $T_{k_{\beta}}$, respectively. 

We have introduced a framework to construct RKHS as a subspace of $L^{2}(\BR^{d}, \bm{W}^{1-2\beta}) \cap C_{0}(\BR^{d})$. The resulting RKHS extracts features of both $L^{2}(\BR^{d}, \bm{W}^{1-2\beta})$ and $C_{0}(\BR^{d})$, e.g., the expansion formula \eqref{eq:expansion} makes sense not only in $L^{2}(\BR^{d}, \bm{W}^{1-2\beta})$ but also pointwise.  The main advantage with the representation in \eqref{eq:RKHS} over the radial-type kernel is as follows. In practice, given data $\left\{ X_{n} \right\}_{n=1}^{N}$ sampled from the target density $f$, we will choose a function in $\mathcal{H}_\beta$ with a finite sum, $\|\vec{m}\|_1\leq M,$ where $M \ll N$, as an estimator for $f$. While the choice of $M$ allows us to specify the theoretical ``bias'' or ``approximation error'', thanks to the orthogonal representation (see Section~\ref{sec:KELR}), the resulting hypothesis function is parametric and the evaluation of $f$ on a new $\bm{x} \in \BR^{d}$ amounts to evaluating ${M+d}\choose{M}$ components of $\left\{\Psi_{\beta,\vec{m}}(\bm{x}) \;|\; \|\vec{m}\|_1\leq M \right\}$. This is computationally much cheaper than evaluating  $f(\bm{x}) = \langle f,k(\cdot,\bm{x})\rangle_\mathcal{H}$, with radial-type kernels, such as $k\left(\bm{x}, \bm{y}\right) = h\left(\|\bm{x}-\bm{y}\|\right)$ for some positive function $h$, since the computation of the inner product requires evaluating $h\left(\|X_n-\bm{x}\|\right)$, for all $n=1,\ldots, N$.

It is worthwhile to mention that although Proposition~\ref{prop:RKHS_d} relies on the weight function $W$ satisfying the condition in Lemma~\ref{lem:bound}, there are kernels and RKHS that are generated from different choices of $W$. For example, when $W$ has a compact support, the resulting kernel and RKHS is an immediate consequence of Mercer's theorem. In the remaining of this subsection, we discuss two examples, the classes of Mercer-type kernels based on the Hermite and Laguerre polynomials. In particular, for Laguerre polynomials, the weight function does not satisfy the assumption in Lemma~\ref{lem:bound} and we use a different approach to justify the boundedness of the kernel.

\subsubsection{Hermite polynomials}\label{sec:Her}

Let $\bm{W}$ be the $d$-dimensional standard Gaussian distribution, that is, $\bm{W}(\bm{x}) = (2\pi)^{-\frac{d}{2}} \exp\left(-\frac{1}{2} \| \bm{x}\|^{2}  \right)$. As a result, the corresponding orthonormal polynomials are the normalized Hermite polynomials, denoted by $\{\psi_{n}\}$. To satisfy the condition in Lemma~\ref{lem:kernel}, we take $\lambda_{n} = \rho^{n}$ with $\rho\in (0,1)$. Following \eqref{eq:kernel_d}, we define the kernel
\begin{equation} \label{eq:kernel}
k_{\beta, \rho}(\bm{x}, \bm{y}): = \sum_{\vec{m}\geq 0} \rho^{\|\vec{m}\|_{1}} \Psi_{\beta, \vec{m}}(\bm{x}) \Psi_{\beta, \vec{m}}(\bm{y}), \quad \forall \bm{x}, \bm{y} \in \BR^{d}, \; \Psi_{\beta, \vec{m}} = \psi_{\vec{m}} \bm{W}^{\beta}.
\end{equation}
For this special case, we do have an explicit expression for $k_{\beta, \rho}$. When $d=\beta = 1$, the kernel $k_{1,\rho}$ in \eqref{eq:kernel} is known as the Mehler kernel \cite{mehler1866ueber} with
\begin{equation}\label{eq:mehler}
k_{1,\rho}(x,y) = \sum_{m=0}^{\infty} \rho^{m} \Psi_{m}(x) \Psi_{m}(y) = \frac{1}{2\pi \sqrt{1-\rho^2}} \exp \left( -\frac{x^2 - 2\rho xy + y^2}{2(1-\rho^2)}  \right).
\end{equation}
For general $d$-dimensional problems, we have  
\begin{eqnarray}\label{eq:kernel_exp}
k_{\beta,\rho}(\bm{x}, \bm{y}) = \left(2\pi \sqrt{1-\rho^2}\right)^{-d} \times \exp \left[ -\frac{1}{2(1-\rho^2)}\left(\|\bm x\|^2 + \|\bm y\|^2 - 2\rho \sum_{i=1}^{d} x_{i}y_{i} \right)   \right] W^{\beta - 1} (\bm{x}) W^{\beta - 1}(\bm{y}),
\end{eqnarray}
and
\begin{equation}
\|k_{\beta, \rho}\|_{\infty} = (2\pi)^{-\frac{\beta d}{2}}\left(1-\rho^2\right)^{-\frac{d}{4}}.\nonumber
\end{equation}
For $\beta \in [\frac{1}{2}, \infty)$ and $\rho \in (0,1)$, we will call the kernel $k_{\beta, \rho}$ in \eqref{eq:kernel_exp} and the corresponding RKHS
\begin{equation}\label{eq:RKHS_exp}
\mathcal{H}_{\beta, \rho}: = \left\{ f= \sum_{\vec{m}\geq 0} \hat{f}_{\vec{m}} \Psi_{\beta, \vec{m}}  \; \Big|  \;  \sum_{\vec{m}\geq 0}\frac{\hat{f}^2_{\vec m}}{\rho^{\|\vec{m}\|_{1}}} < \infty \right\}, 
\end{equation}
following Proposition~\ref{prop:RKHS_d}, the $d$-dimensional Mehler kernel and Mehler RKHS, respectively. The explicit formula of the Mehler kernel \eqref{eq:kernel_exp} is convenient for verifying the properties of the Mehler RKHS. For example, since the kernel $k_{\beta, \rho}$ \eqref{eq:kernel_exp} is smooth, by Lemma~\ref{lem:diff}, we have the following regularity result.

\begin{corollary}\label{lem:smooth}
Every function $f \in \mathcal{H}_{\beta, \rho}$ is smooth with derivatives satisfying,
\begin{equation*}
\quad \left | \partial^{\vec{m}}f(\bm{x})  \right| \leq  \left\|f \right\|_{\mathcal{H}_{\beta, \rho}} \cdot \left ( \partial^{\vec{m}, \vec{m}}k_{\beta, \rho}(\bm{x},\bm{x})\right)^{\frac{1}{2}}, \quad \forall \bm{x}\in \mathbb{R}^{d}, \quad \forall \vec{m}\geq 0.
\end{equation*}
In particular,
\begin{equation}\label{eq:nablaf}
\left| \frac{\partial }{\partial x_{i}}f(\bm{x}) \right | \leq  \left\|f\right\|_{\mathcal{H}_{\beta, \rho}}k^{\frac{1}{2}}_{\beta, \rho}(\bm{x},\bm{x})  \cdot   \sqrt{\frac{\rho}{1-\rho^{2}} +  \left(\frac{1}{1+\rho} + \beta -1 \right)^{2} x_{i} }, \quad i = 1,2,\dots,d.
\end{equation}
\end{corollary}

When $\vec{m} = \vec{0}$, Corollary~\ref{lem:smooth} reduces to \eqref{eq:decay_rate}, which means the functions in $\mathcal{H}_{\beta, \rho}$ yields the same decay rate as $k^{\frac{1}{2}}_{\beta, \rho}(\bm{x}, \bm{x})$. In particular, $\forall f \in \mathcal{H}_{\beta, \rho}$,
\begin{equation} \label{eq:de_rate}
\left|f(\bm{x}) \right| \leq (2\pi)^{-\frac{d}{2}} \left(1-\rho^2\right)^{-\frac{d}{4}} \left\| f \right \|_{ \mathcal{H}_{\beta, \rho}} \exp\left (- \frac{1}{2(1+\rho)} \| \bm{x} \|^{2} \right) \bm{W}^{\beta-1}(\bm{x}).
\end{equation}

\subsubsection{Laguerrre polynomials} \label{sec:Lag}

The normalized Laguerre polynomials, denoted by $\ell^{(\theta)}_{n}(x)$, are the orthonormal polynomials with respect to the Gamma distribution $G(x; \theta) \propto x^{\theta} e^{-x}$, for $\theta =0,1,\dots$, and $x\geq 0$. Following \eqref{eq:kernel_d}, we introduce the $d$-dimensional Gamma distribution and the corresponding normalized Laguerre polynomials
\begin{equation}\label{eq:d_Lag}
G(\bm{x}; \vec{\theta})= \prod_{i=1}^{d} \frac{1}{\Gamma(\theta_{i}+1)} x_{i}^{\theta_{i}}e^{-x_{i}}, \quad    \ell^{(\vec{\theta})}_{\vec{m}}(\bm{x}): = \prod_{i=1}^{d} \ell^{(\theta_{i})}_{m_{i}}(x_{i}), \quad \bm{x} \in [0, \infty)^{d},  \quad \vec{\theta}\in\{0,1,\dots\}^{d},
\end{equation}
respectively. It is clear that the eigenfunctions $\Psi_{\beta, \vec{\theta}, \vec{m}}$ defined as
\begin{equation}\label{eq:eigf_Lag}
\Psi_{\beta, \vec{\theta}, \vec{m}}(\bm{x}): =  \ell^{(\vec{\theta})}_{\vec{m}}(\bm{x})G^{\beta}(\bm{x}; \vec{\theta}), \quad \vec{m}\geq 0
\end{equation} 
form an orthonormal basis on $L^{2}\left([0,\infty)^{d}, G^{1-2\beta}( \;\cdot \;; \vec{\theta}) \right)$. 

Notice that the Laguerre weight function $G(x,\theta)$ does not satisfy the hypothesis of Lemma~\ref{lem:bound}, so we cannot use Lemma~\ref{lem:kernel} to verify the boundedness of the Mercer-type kernel as in Eq. \eqref{eq:kernel}. Nevertheless, we can still verify the boundedness of the kernel using existing identities. We summarize the result as follows.

\begin{lemma} \label{lem:HHF}(Hille-Hardy kernel)
Let $G(\bm{x}; \vec{\theta})$, $\ell^{(\vec{\theta})}_{\vec{m}}(\bm{x})$, and $\Psi_{\beta, \vec{\theta}, \vec{m}}(\bm{x})$ be defined as in Eq. \eqref{eq:d_Lag} and \eqref{eq:eigf_Lag}. For $\rho\in(0,1)$, and $\beta \in [\frac{1}{2}, \infty)$, 
\begin{equation}\label{eq:Mercer_Lag}
k_{\beta, \rho, \vec{\theta}}(\bm{x},\bm{y}) := \sum_{\vec{m} \geq 0} \rho^{\|\vec{m}\|_{1}} \Psi_{\beta, \vec{\theta}, \vec{m}}(\bm{x})\Psi_{\beta, \vec{\theta}, \vec{m}}(\bm{y}), \quad \bm{x}, \bm{y} \in [0, \infty)^{d},
\end{equation}
defines a bounded kernel. In particular, we have the following $d$-dimensional generalized Hille-Hardy formula \cite{watson1933notes}
\begin{eqnarray}\label{eq:HHF}
k_{\beta, \rho, \vec{\theta}}(\bm{x},\bm{y}) = \frac{\rho^{-\frac{1}{2}\|\vec{\theta}\|_{1}}}{(1-\rho)^{d}}\exp \left(-\frac{1+\rho}{2(1-\rho)}\|\bm{x}+\bm{y}\|_{1}\right) G^{\beta-\frac{1}{2}}(\bm{x}; \vec{\theta})G^{\beta-\frac{1}{2}}(\bm{y}; \vec{\theta}) \prod_{i=1}^{d} I_{\theta_{i}} \left( \frac{2\sqrt{x_{i}y_{i}\rho}}{1-\rho}  \right),
\end{eqnarray}
where $I_{\theta_{i}}$ denotes the modified Bessel function \cite{abramowitz1948handbook}.
\end{lemma}
See Appendix~\ref{app:C} for the proof. The Hille-Hardy formula (Eq. \eqref{eq:HHF} with $d=1$) can be interpreted as a generalization of the Mehler kernel \eqref{eq:mehler}, since Hermite polynomials can be entirely deduced from the Laguerre polynomials \cite{szeg1939orthogonal}.

With Lemma~\ref{lem:HHF}, following Proposition~\ref{prop:RKHS_d}, one can construct the RKHS corresponding to the Hille-Hardy kernel in \eqref{eq:HHF} as an analogue of the Melher RKHS in \eqref{eq:RKHS_exp}. One can also specify the regularity of this function space as in  Lemma~\ref{lem:smooth} by checking the boundedness of the derivatives of the Hille-Hardy kernel using some results in \cite{abramowitz1948handbook}, which we omitted here for brevity.

\subsection{Universality of the Mehler RKHS}\label{sec:richness}

As mentioned in the introduction, for a reliable estimation, we would like to construct a hypothesis space (RKHS) that is ``rich'' enough to capture particular behavior of the target function. Specifically, we will use the notion of weighted $c_0$-universality described in Lemma~\ref{lem:weight_c0} to specify the appropriate Mehler RKHS that can capture the decaying property (Gaussian or faster) of the target density function. The key idea is to understand the ``richness'' of the RKHS space $\mathcal{H}_{\beta, \rho}$ as we summarize now.

\begin{corollary} \label{corollary:rich}
For Mehler kernel $k_{\beta, \rho}$ in \eqref{eq:kernel_exp} and Mehler RKHS $\mathcal{H}_{\beta,\rho}$ in \eqref{eq:RKHS_exp}, let 
\begin{equation} \label{eq:qalpha}
q_{\beta, \rho}(\bm{x}) : = \exp\left [-  \left(\frac{1}{2(1+\rho)}+ \frac{\beta-1}{2} \right)\| \bm{x} \|^{2}  \right], \quad \bm{x} \in \BR^{d},
\end{equation}
then $\mathcal{H}_{\beta, \rho}$ is dense in $C_{0}\left(\BR^{d}, q_{\beta, \rho}^{-1}\right)$.
\end{corollary}

\begin{proof}
By Lemma~\ref{lem:weight_c0}, we need to check the following two conditions.
\begin{enumerate}
\item $k_{\beta, \rho}(\cdot, \bm{x})\in C_{0}\left(\BR^{d},q_{\beta, \rho}^{-1}\right)$, $\forall \bm{x} \in \BR^{d}$,
\item $\tilde{k}(\bm{x},\bm{y}):= q_{\beta, \rho}^{-1}(\bm{x}) k_{\beta,\rho}(\bm{x}, \bm{y}) q_{\beta, \rho}^{-1}(\bm{y})$ is $c_{0}$-universal.
\end{enumerate}

Notice that for any fixed $\bm{y} \in \BR^{d}$, as a function of $\bm{x}$, 
\begin{equation*}
k_{\beta, \rho}(\bm{x}, \bm{y} )q^{-1}_{\beta, \rho}(\bm{x}) \propto \exp \left[ -\frac{1}{2(1-\rho^2)}\left( \rho\|\bm x\|^2 - 2\rho \sum_{i=1}^{d} x_{i}y_{i} \right)   \right],
\end{equation*}
that is, $k_{\rho}(\cdot, \bm{y})q_{\rho}^{-1}(\cdot) \in C_{0}(\BR^{d})$ since $\rho\in (0,1)$. By the definition of $C_{0}\left(\BR^{d}, q_{\beta, \rho}^{-1}\right)$ (see Lemma~\ref{lem:weighted}), the first condition holds.

For the second condition, simply notice that 
\begin{equation*}
\tilde{k}(\bm{x},\bm{y}) \propto \exp \left[ -\frac{1}{2(1-\rho^2)}\left(\rho \|\bm x\|^2 + \rho\|\bm y\|^2 - 2\rho \sum_{i=1}^{d} x_{i}y_{i} \right)   \right] = \exp\left[ -\frac{\rho}{2(1-\rho^2)}\|\bm{x} - \bm{y}\|^{2}  \right],
\end{equation*}
and $\frac{\rho}{2(1-\rho^2)}>0$, that is, $\tilde{k}(\bm{x},\bm{y})$ is a Gaussian kernel which is $c_0$-universal.

\end{proof}

\begin{remark}\label{remark_rich}
Corollary~\ref{corollary:rich} suggests that the Mehler RKHS $\mathcal{H}_{\beta, \rho}$ is rich enough to approximate the functions in $C_{0}\left(\BR^{d}, q_{\beta, \rho}^{-1} \right)$. Moreover, notice that in the formula of $q_{\beta, \rho}$ in \eqref{eq:qalpha}, we have
\begin{equation*}
\left\{ \frac{1}{2(1+\rho)}+ \frac{\beta-1}{2}  \; \Big| \; \beta \in (\frac{1}{2}, \infty), \; \rho \in (0,1) \right \} = (0, \infty).
\end{equation*}
Thus, for any continuous target density function $f$ with a Gaussian decay rate, that is, there exist constants $\theta_{f}, C_{f}>0$ such that
\begin{equation}\label{gau_decay}
  \left| f(\bm{x}) \right|  \leq C_{f} \exp\left( -\theta_{f} \| \bm{x} \|^2 \right), \quad \forall \bm{x} \in \BR^{d},
\end{equation}
we can find $\beta^{*} > \frac{1}{2}$ and $\rho^{*} \in (0,1)$ such that 
\begin{equation}\label{eq:tuning}
       \frac{1}{2(1+\rho^{*})}+ \frac{\beta^{*}-1}{2} < \theta_{f},
\end{equation}
and the decay rate \eqref{gau_decay} implies $f \in C_{0}\left(\BR^{d},q_{\beta^{*}, \rho^{*}}^{-1}\right)$. This means that one can approximate any continuous function that has Gaussian (or faster) decaying rate as in \eqref{gau_decay}
up to any desirable accuracy using an estimator that belongs to the Mehler RKHS $\mathcal{H}_{\beta^{*}, \rho^{*}}$, which is a space of functions with Gaussian decay rate slower than \eqref{gau_decay}. In practice, one can use the sample excess kurtosis to approximate the decay rate of the target density from the given samples on each coordinate to identify the value of $\theta_{f}$ in \eqref{gau_decay}. In terms of tuning parameters $\rho$ and $\beta$ via \eqref{eq:tuning}, since the left-hand side of Eq.\eqref{eq:tuning} depends on $(\rho, \beta)$ with $\frac{\beta-1}{2}$ being the leading-order term, we fix $\rho = \rho^{*}\in (0,1)$ and seek for $\beta =\beta^{*}>\frac{1}{2}$ that satisfies the inequality \eqref{eq:tuning}. Theoretically, we can tune parameter $\rho$ as well by further posing decay rate condition to the partial derivative of the density function $f$ (e.g., Eq.\eqref{eq:nablaf}). However, such partial derivatives data, in practice, may not be available.
\end{remark}

We should point out that Corollary~\ref{corollary:rich} relies on the analytical expression of the Mehler kernel, which allows us to find an isometrically isomorphic map that transforms the Mehler kernel into a radial kernel that is known to be $c_0$-universal. Unfortunately, such a technique is not trivially applicable for the Hille-Hardy kernel since its analytical expression, as shown in \eqref{eq:HHF}, is more complicated. One possible approach to justify the ``richness'' of the Hille-Hardy RKHS is developing a weighted-$c_0$ universality characterization based on the Proposition~12 in \cite{sriperumbudur2011universality}, and we leave it as an open problem. Practically, one can choose the appropriate basis guided by the shape parameter $\theta_i$ in \eqref{eq:d_Lag} that can be determined from the sample data by the standard MLE method on each coordinate.

\section{Kernel embedding linear response} \label{sec:KELR}

In this section, we first define the kernel embedding estimate for density functions. Subsequently, we introduce the kernel embedding linear response as the main application. As mentioned in the introduction, finite orthogonal polynomials will admit negative-values on the density estimation. Here, we will address this issue and discuss to which extent the consistency of the kernel embedding linear response can be achieved with the proposed algorithm. We will close this section with an error bound for the empirical kernel embedding estimates based on non-i.i.d. data. For simplicity, we only discuss the linear response estimate based on the Mehler RKHS hypothesis space. Similar techniques can be applied to the Hille-Hardy kernel.

We assume the target $d$-dimensional density function $f$ lives in the Mehler RKHS $\mathcal{H}_{\beta, \rho}$ for some $\rho \in (0,1)$ and $\beta \geq \frac{1}{2}$. Recall that, by the reproducing property and the expansion formula \eqref{eq:expansion}, we have
\begin{equation}\label{eq:rep_for}
       f   = \sum_{\vec{m}\geq 0} \hat{f}_{\vec m} \Psi_{\beta, \vec{m}}, \quad \hat{f}_{\vec{m}}  =  \int_{\mathbb{R}^{d}} f(\bm x) \psi_{\vec m}(\bm x) \bm{W}^{1- \beta}(\bm{x})  \td \bm x, \quad \psi_{\vec m} = \prod_{i = 1}^{d} \psi_{m_i},
\end{equation}
where $\bm{W}(\bm{x}) \propto \exp\left(-\frac{1}{2} \| \bm{x}\|^{2}\right)$, and $\{ \psi_{m_{i}} \}$ are the normalized Hermite polynomials. We define the \textit{order-$M$ kernel embedding estimates} of $f$, denoted by $f_{M}$, as the order-$M$ truncation of Eq. \eqref{eq:rep_for}, that is,
\begin{equation} \label{eq:kme1}
 f_{M} := \sum_{\|\vec{m}\|_{1}\leq M} \hat{f}_{\vec m} \Psi_{\beta,\vec{m}}.
\end{equation}
We should point out that, with this choice of basis representation, we arrive at a polynomial chaos approximation\cite{xiu:2010} of $f$. But the convergence $f_{M}\rightarrow f$ as $M \rightarrow \infty$ is valid in both $L^{2}(\BR^{d}, \bm{W}^{1-2\beta})$ and $C_{0}(\BR^{d})$. Another remark is that, although the formula of the kernel embedding estimate \eqref{eq:kme1} is independent of the parameter $\rho$, this parameter does affect the rate of convergence of the residual error $ \left\| f - f_{M}\right\|_{\mathcal{H}_{\beta,\rho}}$(see Proposition~\ref{prop:mixing}).

In practice, the integral in \eqref{eq:rep_for} can be approximated by a Monte-Carlo average, that is,
\begin{equation}\label{eq:kme_MC}
\hat{f}_{\vec{m}} \approx \hat{f}_{\vec{m},N} :=  \frac{1}{N} \sum_{n=1}^{N} \psi_{\vec{m}}(X_{n}) \bm{W}^{1-\beta}(X_{n})  ,
\end{equation}
where $\{X_{n}\}_{n=1}^{N}$ are sampled from the target density function $f$. Here, the decay rate of $f$ in \eqref{eq:de_rate} ensures the convergence of $\hat{f}_{\vec{m},N}$ in \eqref{eq:kme_MC} as the sample size $N$ goes to $\infty$. One can define the order-$M$ \textit{empirical kernel embedding estimate} of $f$ as
\begin{equation}\label{eq:kme2}
f_{M,N} : = \sum_{\|\vec{m}\|_{1}\leq M} \hat{f}_{\vec{m},N} \Psi_{\beta,\vec{m}}.
\end{equation}
In our application, the target function is $f = p_{eq}$. Here, the sample $\{X_{n}\}$ corresponds to a time series, $\{X(t_n)\}$, of the unperturbed dynamics \eqref{eq:Ito} generated at the equilibrium. We should point out that the error bound, to be discussed in Proposition~\ref{prop:mixing}, will account for the fact that the empirical estimator in \eqref{eq:kme_MC} is based on non-i.i.d. data.

Let's replace the unknown $p_{eq}$ in \eqref{eq:lin_oper} by its order-$M$ kernel embedding estimates, denoted as $p_{M}$. We can naively define the order-$M$ kernel embedding linear response operator as
\begin{equation}\label{eq:KMLE1}
\hat{k}_{A}(t): = \mathbb{E}_{p_{eq}} \left [A(X(t)) \otimes \hat{B} (X(0)) \right], \quad \hat{B}_{i}(X):= -\frac{\partial_{X_i}\left[c_{i}(X)p_{M}(X)\right]}{p_{M}(X)},
\end{equation}
for integrable observable $A(X)$ and external forcing $c(X)\delta f(t)$. Unfortunately, Eq. \eqref{eq:KMLE1} is not a well-defined estimator for the linear response operator $k_{A}(t)$, since the basis function $\Psi_{\beta,\vec{m}}(\bm{x})$ in \eqref{eq:kernel}, as a product of normalized Hermite polynomials and the standard Gaussian distribution, can be negative.

To overcome this issue, we introduce the following restriction to the order-$M$ density estimate, 
\begin{equation*}
    p_{M,\delta}:= \max \{p_{M}, \delta\} = p_{M}\cdot \chi_{D_{M,\delta}}, \quad \delta > 0,
\end{equation*}
where $D_{M,\delta}: = \{\bm{x}\in \mathbb{R}^{d}\; |\; p_{M}(\bm{x})\geq \delta\}$ and $\chi_{D_{M,\delta}}$ denotes the characteristic function with respect to $D_{M,\delta}$. For the analysis below, it is helpful to see that 
\begin{equation}\label{eq:D_Mdelta}
   \left\{\bm{x}\; |\; p_{eq}(\bm{x})\geq \delta -\|p_{eq} - p_{M}\|_{\infty} \right\} \supset D_{M,\delta} \supset \left\{\bm{x}\; |\; p_{eq}(\bm{x})\geq \delta +\|p_{eq} - p_{M}\|_{\infty} \right\}.
\end{equation}
Under the decaying rate assumption of $p_{eq}$, the set $\{\bm{x}\in \mathbb{R}^{d}\; |\; p_{eq}(\bm{x})\geq \delta \}$ is compact $\forall \delta>0$. Indeed, with a positive lower bound on $p_{M}$, the estimator $\hat{B}_{i}$ in \eqref{eq:KMLE1} is well-defined on $D_{M,\delta}$. Here, the choice of $\delta>0$ depends on the desired precision of the estimator for a fixed $M$, as we shall explain next.

Similar to the correction of the density estimate, we define the restricted linear response operator
\begin{equation}\label{eq:KMLE2}
k_{A,D_{M,\delta}}(t): = \mathbb{E}_{p_{eq}} \left [A(X(t)) \otimes B_{D_{M,\delta}} (X(0)) \right], \quad B_{D_{M,\delta}}: = B \cdot \chi_{D_{M,\delta}}, \; \forall t \geq 0.
\end{equation}
Further, we assume that both $A(\cdot)$ and $B(\cdot)$ are fixed and have finite second moments with respect to $p_{eq}$, then, $B_{D_{M,\delta}}$ in \eqref{eq:KMLE2} also has a finite second moment. Thus, according to Lemma~\ref{lem:bound_kA}, $k_{A,D_{M,\delta}}(t)$ is well-defined  $\forall t \geq 0$, and,  replacing $B$ by $B(1-\chi_{D_{M,\delta}})$ in \eqref{eq:CS_1}, we reach
\begin{equation}\label{eq:KMLE3}
\left | k_{A}(t) - k_{A,D_{M,\delta}}(t) \right|  \leq  \left(\int_{\mathbb{R}^{d}} A^{2}(\bm{x}) p_{eq}(\bm{x}) \td \bm{x} \right)^{\frac{1}{2}} \otimes  \left(\int_{D_{M,\delta}^{c}}  B ^{2}(\bm{x}) p_{eq}(\bm{x}) \td \bm{x}\right)^{\frac{1}{2}}, \quad \forall t \geq 0.
\end{equation}
We claim that $\forall \epsilon > 0 $, there exists a $\delta^{*} >0$ such that $\forall \delta\in(0, \delta^{*})$ and for $M$ large enough,
\begin{equation}\label{eq:deps}
    |k_{A}(t) - k_{A,D_{M,\delta}}(t)|<\epsilon, \quad \forall t \geq 0,
\end{equation}
where the inequality is defined componentwise and uniformly in $t$. To clarify, without loss of generality, we assume that both $A$ and $B$ are scalars. With ${\mathbb{E}_{p_{eq}}[B^2(X)]}<\infty$, by the continuity of the integral, there exists a compact set $D\subset \operatorname{supp}\{p_{eq}\}$ such that
\begin{equation*}
\int_{D^{c}}  B ^{2}(\bm{x}) p_{eq}(\bm{x}) \td \bm{x} < \frac{\epsilon^2}{\mathbb{E}_{p_{eq}}[A^2(X)]},
\end{equation*}
where $D^c=\mathbb{R}^d \backslash D$ denotes the complementary set. Take $\delta^{*}:= \min\{p_{eq}(\bm{x})\; | \; \bm{x} \in D \}$. As a result of the uniform convergence of $p_{M}$ to $p_{eq}$, we have $\forall \delta \in (0, \delta^{*})$ and $M$ large enough, $\delta + \|p_{eq} - p_{M}\|_{\infty} \leq \delta^{*}$, that is,
\begin{equation*}
D_{M,\delta} \supset \{\bm{x}\in \mathbb{R}^{d}\; |\; p_{eq}(\bm{x})\geq \delta +\|p_{eq} - p_{M}\|_{\infty} \} \supset  \{\bm{x}\in \mathbb{R}^{d}\; |\; p_{eq}(\bm{x})\geq \delta^{*} \} \supset D,
\end{equation*}
where we have used the relation in \eqref{eq:D_Mdelta}. Thus,
\begin{equation*}
\int_{D_{M,\delta}^{c}}  B ^{2}(\bm{x}) p_{eq}(\bm{x}) \td \bm{x}<\int_{D^{c}}  B ^{2}(\bm{x}) p_{eq}(\bm{x}) \td \bm{x}  < \frac{\epsilon^2}{\mathbb{E}_{p_{eq}}[A^2(X)]},
\end{equation*}
and, following Eq.\eqref{eq:KMLE3}, Eq.\eqref{eq:deps} holds for all $t \geq 0$. 

The implication of \eqref{eq:deps} is the existence of the domain $D_{M,\delta}$ such that for large enough $M>0$, the restricted linear response estimator in \eqref{eq:KMLE2} is consistent with $k_A(t)$ up to the desirable precision $\epsilon>0$. 
Numerically, one can consistently estimate the restricted linear response operator in \eqref{eq:KMLE2} with the following well-defined estimator of $k_{A,D_{M,\delta}}$,
\begin{equation}\label{eq:deps_hat}
\quad \hat{k}_{A,D_{M,\delta}}(t) := \mathbb{E}_{p_{eq}} \left [A(X(t)) \otimes \hat{B}_{D_{M,\delta}} (X(0)) \right],  \quad \hat{B}_{D_{M,\delta}} = \hat{B}\cdot \chi_{D_{M,\delta}}.
\end{equation}
The following proposition characterizes the consistency of the estimator $\hat{k}_{A,D_{M,\delta}}(t)$ as $M \rightarrow \infty$.

\begin{proposition} \label{prop:consis}
Consider a $d$-dimensional It\^o diffusion \eqref{eq:Ito} with a positive equilibrium density function $p_{eq} \in \mathcal{H}_{\beta, \rho}$, and its perturbed dynamics \eqref{eq:Ito_per}. Assume that, both, the observable $A(\cdot)$ and the conjugate variable $B(\cdot)$ in \eqref{eq:lin_oper} have finite second-moments with respect to $p_{eq}$. For the restricted linear response operator $k_{A,D_{M,\delta}}(t)$ in \eqref{eq:KMLE2} and its estimator $\hat{k}_{A,D_{M,\delta}}(t)$ in \eqref{eq:deps_hat}, we have
\begin{equation} \label{eq:kA_con}
 \lim_{M\rightarrow +\infty} \|k_{A,D_{M,\delta}}(t) - \hat{k}_{A,D_{M,\delta}}(t)\|_{F} = 0,
\end{equation}
uniformly in $t$. Here, $\|\cdot\|_{F}$ denotes the Frobenius norm of matrices.
\end{proposition}

\begin{proof}
To simplify the discussion, without loss of generality, we assume that $c(\bm{x})$ in the perturbed dynamics \eqref{eq:Ito_per} is a scalar function. As a result, following Eq. \eqref{eq:lin_oper}, we have
\begin{equation} \label{eq:B_Bhat}
B(\bm{x}) = \left(\nabla p_{eq}(\bm{x}) c(\bm{x}) + \nabla c(\bm{x}) p_{eq}(\bm{x})\right) p_{eq}^{-1}(\bm{x}), \hat{B}(\bm{x}) = \left( \nabla p_{M}(\bm{x}) c(\bm{x}) + \nabla c(\bm{x}) p_{M}(\bm{x})\right) p_{M}^{-1}(\bm{x}).
\end{equation}
By the integrability assumptions of $A$ and $B$, both $k_{A,D_{M,\delta}}(t)$ and its estimator $\hat{k}_{A,D_{M,\delta}}(t)$ are well-defined, and we have
\begin{equation} \label{eq:term_zero}
\begin{split}
\left | k_{A,D_{M,\delta}}(t) - \hat{k}_{A,D_{M,\delta}}(t) \right| & = \mathbb{E}_{p_{eq}} \left [A(X(t)) \otimes \left( B_{D_{M,\delta}} - \hat{B}_{D_{M,\delta}}\right) (X(0)) \right]  \\
& \leq \mathbb{E}_{p_{eq}}\left[A(X)^{2}\right] \otimes \mathbb{E}_{p_{eq}}\left[\left( B_{D_{M,\delta}}(X) - \hat{B}_{D_{M,\delta}}(X)\right)^{2}\right].
\end{split}
\end{equation}
Furthermore we have
\begin{equation*}
\begin{split}
 \mathbb{E}_{p_{eq}} & \left[\left( B_{D_{M,\delta}}(X) - \hat{B}_{D_{M,\delta}}(X)\right)^{2}\right] 
= \int_{D_{M,\delta}}  \left( B - \hat{B}\right)^{2}(\bm{x}) p_{eq} (\bm{x}) \td \bm{x} \\ 
& = \int_{D_{M,\delta}}  \left[ \left(B - Bp_{eq}p^{-1}_{M} \right) +  \left(Bp_{eq}p^{-1}_{M} - \hat{B}\right)\right]^{2}(\bm{x}) p_{eq} (\bm{x}) \td \bm{x}  \\
& = \int_{D_{M,\delta}}  \left[B(1-p_{eq} p_{M}^{-1})+ \left( \nabla (p_{eq}-p_{M}) c +\nabla c (p_{eq} - p_{M})\right) p_{M}^{-1} \right]^{2}(\bm{x}) p_{eq} (\bm{x}) \td \bm{x} \\
&  \leq  2\int_{D_{M,\delta}} B^{2}(1-p_{eq} p_{M}^{-1})^{2}(\bm{x}) p_{eq}(\bm{x}) \td \bm{x} +2\int_{D_{M,\delta}}  p_{M}^{-2}\left[ \nabla (p_{eq}-p_{M}) c +\nabla c (p_{eq} - p_{M})\right]^{2}(\bm{x}) p_{eq}(\bm{x}) \td \bm{x}\\
& \leq \frac{2}{\delta^2} \left\|p_{M} -  p_{eq} \right\|_{\infty}^{2}\int_{D_{M,\delta}} B^{2}(\bm{x}) p_{eq}(\bm{x}) \td \bm{x} +\frac{2}{\delta^2} \int_{D_{M,\delta}} \left[ \nabla (p_{eq}-p_{M}) c +\nabla c (p_{eq} - p_{M})\right]^{2}(\bm{x})p_{eq}(\bm{x}) \td \bm{x},
\end{split}
\end{equation*}
where we have used the fact that $\delta$ is the lower bound of $p_{M}$ on $D_{M,\delta}$. By Lemma~\ref{lem:uniform} and Corollary~\ref{lem:smooth}, we have the uniform convergence of $p_{M} \rightarrow p_{eq}$ and $\nabla p_{M} \rightarrow \nabla p_{eq}$ as $M\rightarrow \infty$, that is,
\begin{equation} \label{eq:term_one}
\begin{split}
 \lim\sup_{M\rightarrow \infty}  &\int_{D_{M,\delta}}  \left[ \nabla (p_{eq}-p_{M}) c +\nabla c (p_{eq} - p_{M})\right]^{2}(\bm{x})p_{eq}(\bm{x}) \td \bm{x} \\
\leq  \lim_{M\rightarrow \infty} &\int_{\{\bm{x}\in \mathbb{R}^{d}\; |\; p_{eq}(\bm{x})\geq \frac{\delta}{2} \}}  \left[ \nabla (p_{eq}-p_{M}) c +\nabla c (p_{eq} - p_{M})\right]^{2}(\bm{x})p_{eq}(\bm{x}) \td \bm{x} =0,
\end{split}
\end{equation}
where we used the compactness of $\{\bm{x}\in \mathbb{R}^{d}\; |\; p_{eq}(\bm{x})\geq \frac{\delta}{2} \}$ to pass the limit. The first term, 
\begin{equation}\label{eq:term_two}
\left\|p_{M} -  p_{eq} \right\|_{\infty}^{2}\int_{D_{M,\delta}} B^{2}(\bm{x}) p_{eq}(\bm{x}) \td \bm{x} \leq \left\|p_{M} -  p_{eq} \right\|_{\infty}^{2}\int_{\mathbb{R}^d} B^{2}(\bm{x}) p_{eq}(\bm{x}) \td \bm{x} \rightarrow 0,
\end{equation}
as $M \rightarrow \infty$. Here, we used the assumption that $B$ has a finite second moment with respect to $p_{eq}$. Combining \eqref{eq:term_one} with \eqref{eq:term_two}, we reach the convergence
\begin{equation*}
\lim_{M\rightarrow \infty} \mathbb{E}_{p_{eq}}\left[\left( B_{D_{M,\delta}}(X) - \hat{B}_{D_{M,\delta}}(X)\right)^{2}\right] = 0.
\end{equation*}
Together, Eq.~\eqref{eq:term_zero} and the square integrability of $A$, Eq. \eqref{eq:kA_con} holds and the convergence is uniform in $t$.

\end{proof}

\begin{remark}\label{rem:cond_forcing}
It is worthwhile to mention that, for $p_{eq} \in \mathcal{H}_{\beta, \rho}$, by the decay rate of $p_{eq}$ \eqref{eq:de_rate} and $\nabla p_{eq}$ \eqref{eq:nablaf}, a sufficient condition for $B$ in \eqref{eq:B_Bhat} to have a finite second moment is that $c(\cdot)$ is continuously differentiable and 
\begin{equation*}
\exists\; \theta \in \left( 0, \frac{1}{4(1+\rho)} + \frac{\beta-1}{4} \right) \text{ s.t. }\sup_{\bm{x} \in \BR^{d}} \left( |c(\bm{x})| + \|\nabla c(\bm{x}) \|  \right)  \exp\left(-\theta \|x\|^{2} \right) < \infty.
\end{equation*}
This is a concrete condition that can be used in practice to check whether the FDT response is well-defined under external forces of the form $f(x,t)=c(x)\delta f(t)$ for unperturbed dynamics with equilibrium density having Gaussian (or faster) decay rate of arbitrary variance. 
\end{remark}

In practice, we do not have to determine the domain $D_{M,\delta}$. When approximating $\hat{k}_{A, D_{M,\delta}}$ in \eqref{eq:deps_hat} via the Monte-Carlo method (Eq. \eqref{eq:MC_approx}), it is enough to truncate those samples $X_{n}$ where  $p_{M,N}(X_n)< \delta$ with $p_{M,N}$ given by \eqref{eq:kme2}. Numerically, given a training data set of length $N$, there are two parameters to be determined: $M$, as the order of the estimator, and $\delta$, as the threshold of the truncation. For each $M$, to make full use of the data, we set
\BEA
\delta = \delta_{M}:= \min_i \{p_{M,N}(X_{i})>0\}.\label{eq:deltaM}
\EEA
Besides $\delta_M$ \eqref{eq:deltaM}, we further introduce:
\BEA
\mathcal{R}_M &:=& \frac{\mbox{\{Number of $X_i$ where $p_{M,N}(X_i) < \delta_M$\}}}{N}, \label{crit1} \\
\eta_M &:=&  \| p_{M+1,N} - p_{M,N} \|_{L^{2}(\bm{W}^{1-2\beta})} = \left(\sum_{\|\vec{m}\|_1 = M+1} \hat{p}^2_{\vec{m},N}\right)^{\frac{1}{2}}, \label{crit2}
\EEA
where $p_{M,N}$ and $\hat{p}_{\vec{m},N}$ are given by Eqs.~\eqref{eq:kme2} and \eqref{eq:kme_MC} with $f = p_{eq}$, respectively. The parameter $\mathcal{R}_M$ in \eqref{crit1} characterizes the number of training samples whose estimated density function values are smaller than the threshold $\delta_M$; we refer to this parameter as the rejection probability. Small value of $\mathcal{R}_M$ corresponds to small truncation error ($\epsilon$ in \eqref{eq:deps}). The parameter $\eta_M$ in \eqref{crit2} provides an indicator whether the estimator converges as $M\to\infty$. In practice, the computation of $\eta_M$ can be done in an iterative manner by increasing the parameter $M$ since the orthogonality condition allows one to identify the error by the expansion coefficients $\hat{p}_{\vec{m},N}$, which can be computed consecutively. For accurate estimations, we will choose $M$ such that both $\mathcal{R}_M$ and $\eta_M$ are small.

We summarize the numerical scheme into a pseudo-algorithm. Given a time series $\{X_i\}$ sampled from the target unknown $d$-dimensional density function $p_{eq}$, and a target two-point satisfies $k_{A}(t)$ in \eqref{eq:lin_oper}, we proceed our estimation following these steps:
\begin{enumerate}
    \item For each dimension, run a statistical test to identify the tail of the marginal distribution of $p_{eq}$. Subsequently, we choose an appropriate RKHS basis based on the tail information (e.g., tuning the parameter $(\beta, \rho)$ in the Mehler RKHS $\mathcal{H}_{\beta, \rho}$ in \eqref{eq:RKHS_exp}, see Remark~\ref{remark_rich} for the details);
    \item Use tensor product to construct the basis on $\mathbb{R}^d$, and compute the empirical kernel embedding estimates of $p_{eq}$, $p_{M,N}$, following \eqref{eq:kme_MC} and \eqref{eq:kme2} with $f = p_{eq}$;
    \item Use $\delta$ \eqref{eq:deltaM}, $\mathcal{R}_{M}$ \eqref{crit1}, and $\eta_{M}$ \eqref{crit2} to identify the proper order of the estimates and truncate sample points with estimate value, $p_{M,N}(X_i)$, less than $\delta$;
    \item Replace the unknown function $B$ in the two-point statistics $k_{A}(t)$ \eqref{eq:lin_oper} with $\hat{B}$ in \eqref{eq:KMLE1} and estimate the value via Monte-Carlo based on the truncated data set.
\end{enumerate}

\begin{figure}[ht!]
    \center
    \includegraphics[width=0.48\textwidth]{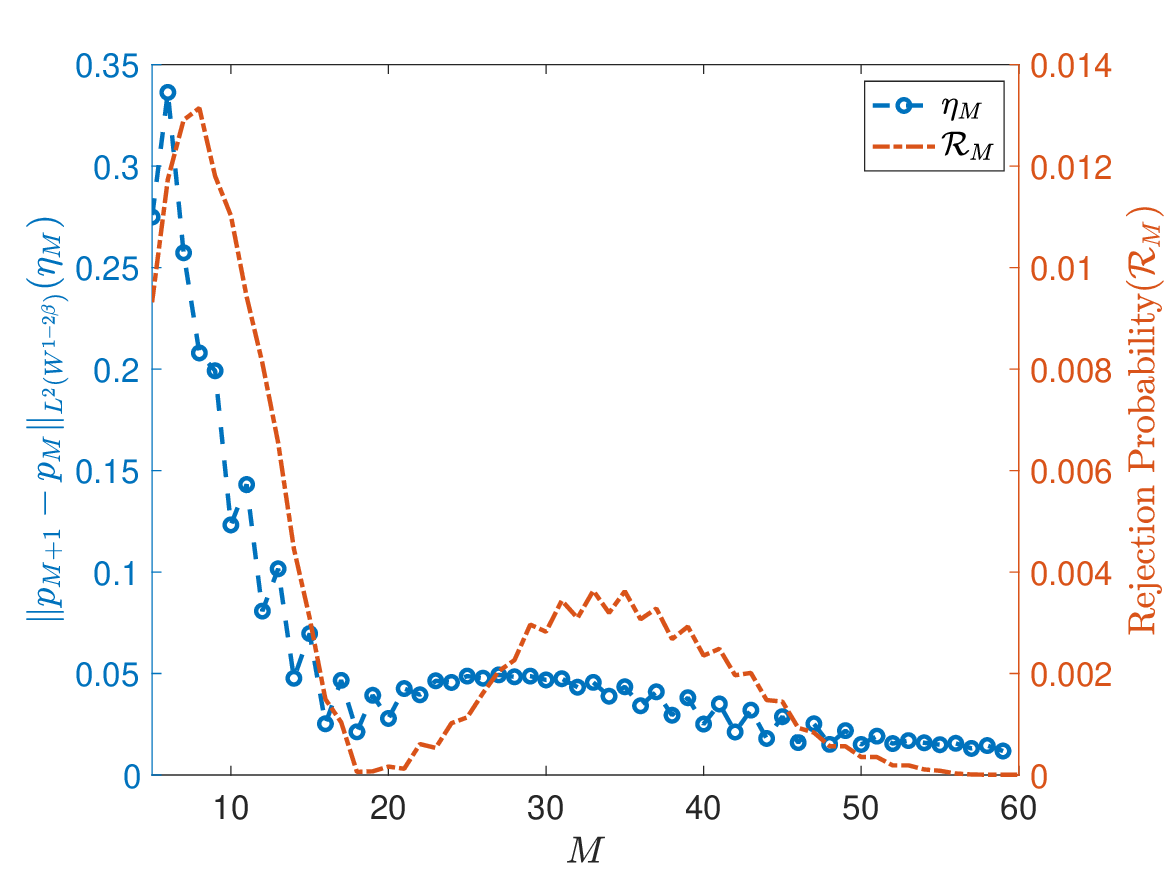}
    \includegraphics[width=0.48\textwidth]{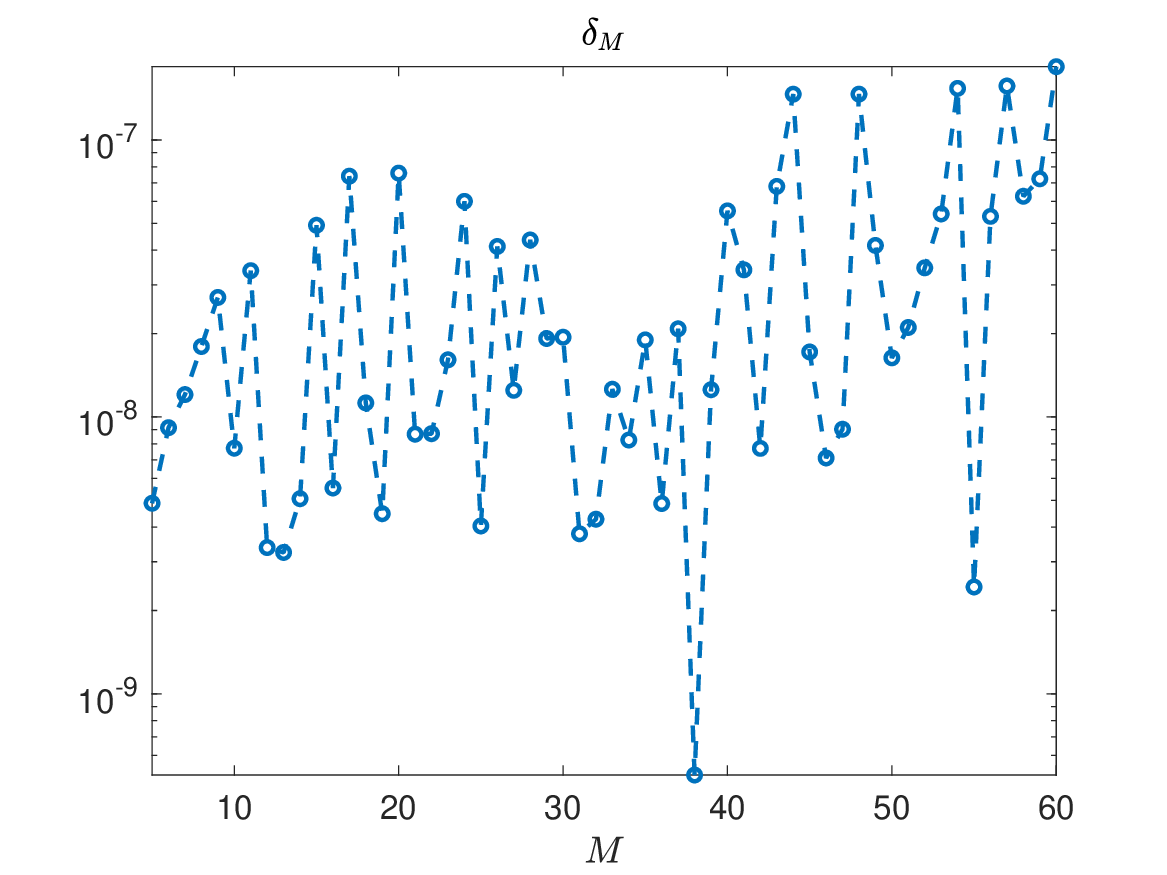}
    \caption{On the left panel, we show an example of how $\eta_M$ in \eqref{crit2} behaves as a function of $M$. On the same panel, we also show the rejection probability $\mathcal{R}_M$ as a function of $M$. One can see that as the algorithm converges (with small $\eta_M$), the rejection probability converges to a relatively small value. On the right panel, we show $\delta_M$ as a function of $M$. Here, there is no pattern for $\delta_M$. In practice, since $\delta_M$ can be arbitrarily small, one can set $\delta$ in \eqref{eq:KMLE2} to be slightly larger than the floating-point single precision to guarantee a well-posed estimator. The results in this figure are based on the gradient system to be discussed in Section~\ref{sec:Triple}.}
\label{fig:numerics}
\end{figure}

In Figure~\ref{fig:numerics}, we provide a numerical example as an illustration. In this example, the data points are the time series of the solutions of the dynamical system in Section~\ref{sec:Triple}. Here, we illustrate the effectiveness of the criterion in \eqref{crit1}-\eqref{crit2} for choosing the parameter $M$. First, notice that $\delta_M$ specified as in \eqref{eq:deltaM} is roughly constant as a function of $M$ with values fluctuating around the floating-point single precision, $10^{-8}$. To guarantee a well-defined estimator, one may need to set $\delta$ to be slightly larger than the floating-point single precision since \eqref{eq:deltaM} can be arbitrarily small. In this example, the choice of $M=60$ will correspond to $\delta_M>10^{-7}$ that is larger than the floating-point single precision. Notice that while the value of $\eta_M$ fluctuates, it eventually converges as $M>50$. On the same panel, we also plot the rejection probability $\mathcal{R}_M$, which is smaller than .4\% as $M$ increases. For $20<M<50$, the large value of $\mathcal{R}_M$ is due to Gibbs phenomenon (as the estimate becomes oscillatory near zeros). As we shall see in Section~\ref{sec:Triple}, we will obtain a reasonably accurate estimation of $k_A$ with $M=60$. We will also use the same criterion to choose $M$ in the other numerical example in Section~\ref{sec:Lan}.

The following proposition addresses the Monte-Carlo error of the order-$M$ empirical kernel embedding estimate $f_{M,N}$ \eqref{eq:kme2}.

\begin{proposition} \label{prop:mixing}

Given a $d$-dimensional stationary (time-homogeneous) Markov process $X(t)$ with a discretization $\left \{X_{n} := X((n-1)\Delta t) \right\}_{n=1}^{N}$ that satisfying
\begin{enumerate}
\item $\{X_{n}\}$ yields a stationary distribution $f \in \mathcal{H}_{\beta, \rho}$ for some $\rho \in (0,1)$ and $\beta \in [\frac{1}{2}, \frac{1}{1+\rho}]$.

\item $\{X_{n}\}$ is an $\alpha$-mixing process \cite{Davydov:68} with mixing coefficient $\alpha(k)$ satisfies
\begin{equation}\label{eq:cond2}
C_{\epsilon} := \sum_{k=1}^{\infty} \alpha(k)^{\frac{\epsilon}{2+\epsilon}} < \infty,
\end{equation}
for some $\epsilon \in (0,2]$.
\end{enumerate}

Then, for $d\leq \frac{3}{2}M +1$, we have the following error bound for the empirical kernel embedding estimate $f_{M,N}$ in \eqref{eq:kme2}
\begin{eqnarray}\label{eq:error_bound}
\mathbb{E}\left [ \left\|f - f_{M,N} \right\|_{L^{2}(\bm{W}^{1-2\beta})}^{2} \right] < \frac{1}{N}\left[ C_{f}(M+1)^{d} + 24 C_{\epsilon} C^{\frac{1}{2}}_{f} M^{d-1} \left( \frac{5}{2} \right) ^{M+3}\right] + \rho^{M+1} \left\| f \right \|^{2}_{\mathcal{H}_{\beta, \rho}},
\end{eqnarray}
where $C_{f} := (2\pi)^{(\beta-1)\frac{d}{2}} \left(1-\rho^2\right)^{-\frac{d}{4}}\left\| f \right \|_{\mathcal{H}_{\beta, \rho}}$, and the expectation in \eqref{eq:error_bound} is defined over random variables $\left \{X_{n} \right\}_{n=1}^{N}$.
\end{proposition}

See Appendix~\ref{app:B} for the proof. The first part of the upper bound in \eqref{eq:error_bound}, with factor $\frac{1}{N}$, corresponds to the estimation (or Monte-Carlo) error, which consists of the variance and covariance terms. The covariance term accounts for the error due to averaging over correlated or non-i.i.d. samples. The second part, depending on the parameter $\rho\in (0,1)$, corresponds to the approximation error or bias due to finite truncation, $M$. Asymptotically, it is clear that the error vanishes as $M\to \infty$ after $N\to\infty$.

If $\{X_{n}\}_{n=1}^{N}$ are i.i.d. samples of $f$, the error bound \eqref{eq:error_bound} reduces to
\begin{equation}\label{eq:error_iid}
\mathbb{E}\left [ \left\|f - f_{M,N} \right\|_{L^{2}(\bm{W}^{1-2\beta})}^{2} \right] < \frac{C_{f}}{N} (M+1)^{d}  + \rho^{M+1} \left\| f \right \|^{2}_{\mathcal{H}_{\beta, \rho}},
\end{equation}
where we are left with the variance and bias terms in the upper bound. We should point out that, in this case, the assumption of $d\leq \frac{3}{2}M+1$, which is used to obtain the upper bound for the non-i.i.d. estimation error  \eqref{eq:error_bound}, can be neglected. First, we note that the error estimate shows that the approximation error (or bias) is independent of the dimension. In fact, the approximation error in \eqref{eq:error_iid} decays exponentially as a function of the model parameters, $M$, since $0<\rho<1$. This fact is not so surprising since the target function, $f \in \mathcal{H}_{\beta,\rho}$, is bounded, smooth, and has Gaussian decay. However, the estimation error depends exponentially on the dimension. This curse of dimensionality is due to the use of the tensor product on an orthonormal basis. While the constants in error bound may be different, the estimation with the Hille-Hardy RKHS should also suffer the same curse of dimensionality issue, and the method of proof is not different from the one presented in Appendix~\ref{app:B}.

For data points that lie on an $m$-dimensional compact manifold $\mathcal{M}$ embedded in $\mathbb{R}^d$, it was shown in \cite{JH:18} that the estimation error (or variance) can be improved to be exactly $MN^{-1}$ using the Mercer kernel constructed by the eigenvalues correspond to the orthonormal eigenfunctions of the Laplace-Beltrami (or a weighted Laplacian) operator, defined on functions that take values on $\mathcal{M}$. In fact, if the target function is a Sobolev class $H^\ell(\mathcal{M})$, one can show that the approximation error is of order $M^{-\frac{2\ell}{m}}$, where the Weyl asymptotic for the eigenvalue of the Laplace-Beltrami (or a weighted Laplacian) operator on compact manifolds \cite{colbois2013eigenvalues} is used. Balancing these errors, we obtain the famous optimal bias-variance trade-off rate of order $N^{-\frac{2\ell}{2\ell+m}}$ for linear estimators \cite{stone1982optimal}. When the target function is smooth such as $\ell=m$, the optimal error rate of order $N^{-\frac{2}{3}}$ suggests that this approach overcomes the curse of dimension. While this is appealing, the numerical method in approximating the orthonormal eigenfunctions of the Laplace-Beltrami operator, such as the diffusion maps \cite{cl:06}, suffers from the curse of dimension with a spectral convergence rate of order-$\left(\frac{\log N}{N}\right)^{\frac{1}{2m}}$ \cite{trillos2019error}. Thus, the curse of dimensionality is not completely avoided. The only advantage of such nonparametric estimation is when the target function is supported on a low-dimensional intrinsic manifold $\mathcal{M}$ that is embedded in a very high-dimensional ambient space $\mathbb{R}^d$. 

In the next section, we will numerically verify the kernel embedding linear response estimators, constructed using the Mehler RKHS and the Hille-Hardy RKHS. For applications with target functions defined on smooth manifolds embedded in $\mathbb{R}^d$, see \cite{berry2017correcting, JH:18}.

\section{Numerical Examples} \label{sec:num}

In this section, we will test the proposed method on two models,  a stochastic gradient system with a triple-well potential (Section~\ref{sec:Triple}) and a Langevin equation with the Morse potential (Section~\ref{sec:Lan}). In \cite{HLZ:19}, the authors have explored these two examples.  Here, our goal is to demonstrate the numerical implementation of the kernel embedding linear response approach introduced in Section~\ref{sec:KELR} based on different choices of orthogonal polynomials. We also stress that we picked examples with explicit equilibrium distributions so that we can verify the estimates by directly comparing to the exact linear response via the FDT formula as reviewed in Section~\ref{sec:FDT}. 

As for a comparison, we consider the classical kernel density estimator (KDE). Formally, given $\{X_n\}_{n=1}^{N}$ sampled from a target density function $\peq$, the KDE estimator is given by,
\begin{equation}\label{eq:KDE}
\tilde{p}_{N}(\bm{x}) = \frac{1}{N} \sum_{n=1}^{N} \frac{1}{h} K\left(\frac{\bm{x} - X_i}{h} \right),
\end{equation}
where $h>0$ is the bandwidth parameter, and $K$ is a smooth density function with zero mean and finite second moment \cite{Wass_book:06}.  The uniform convergence of the KDE estimates and its derivatives have been discussed in \cite{nadaraya1965non} and \cite{schuster1969estimation}, respectively. Similar to the kernel embedding linear response in Section~\ref{sec:KELR}, one may replace $p_{eq}$ in \eqref{eq:lin_oper} with the corresponding KDE estimate, $\tilde{p}_{N}$, and obtain the KDE-based approximation of the linear response statistics. As opposed to the kernel embedding formulation in Section~\ref{sec:KELR}, the KDE-based estimator possesses the positivity but does not have the orthogonality of the basis functions. As a result, we do not need to propose any truncation to the sample points in the numerical implementation. However, each evaluation of $\tilde{p}_{N}$ or its derivatives requires visiting $K$ a total of $N$ times, which makes the computation expensive.

As with the kernel regression in \eqref{eq:KDE}, the choice of kernel $K$ is not crucial, but the choice of bandwidth $h$ is important. For simplicity, throughout the section, we fix $K$ to be the Gaussian kernel, $K(\bm{x}) = (2\pi)^{-\frac{d}{2}}\exp\left(-\frac{1}{2}\|\bm{x}\|^2\right)$, and set $h$ to be the Silverman's bandwidth \cite{silverman1986density},
\begin{equation*}
h = \left(\frac{4}{(d+2)N}\right)^{\frac{1}{d+4}}\sigma( \{X_n\} ),
\end{equation*}
where $\sigma(\{X_n\})$ denotes the standard deviation of the given sample $\{X_n\}_{n=1}^{N}$. While more elaborate bandwidth selection methods can be done (such as the cross validation), we do not pursue them due to the very slow computational time in computing the linear response statistics (see Table~\ref{tab:run_time} for a wall-clock time in both examples below with $N=10^7$).

For the triple-well model, the potential function contains a quadratic retaining potential, which introduces a Gaussian tail to the density function. Thus, we will consider a two-dimensional Mehler kernel in the kernel embedding linear response. For the Langevin equation, the marginal distribution of the velocity $v$ is Gaussian, while the marginal distribution of the displacement $x$, governed by the Morse potential, is asymmetric. In computing the kernel embedding linear response, to obtain the best result, the kernel will be derived based on a tensor product of Hermite (for the variable $v$) and Laguerre (for the variable $x$) polynomials.

\subsection{A gradient system with a triple-well potential} \label{sec:Triple}

\begin{figure}[ht!] 
	\centering
	\includegraphics[scale=0.41]{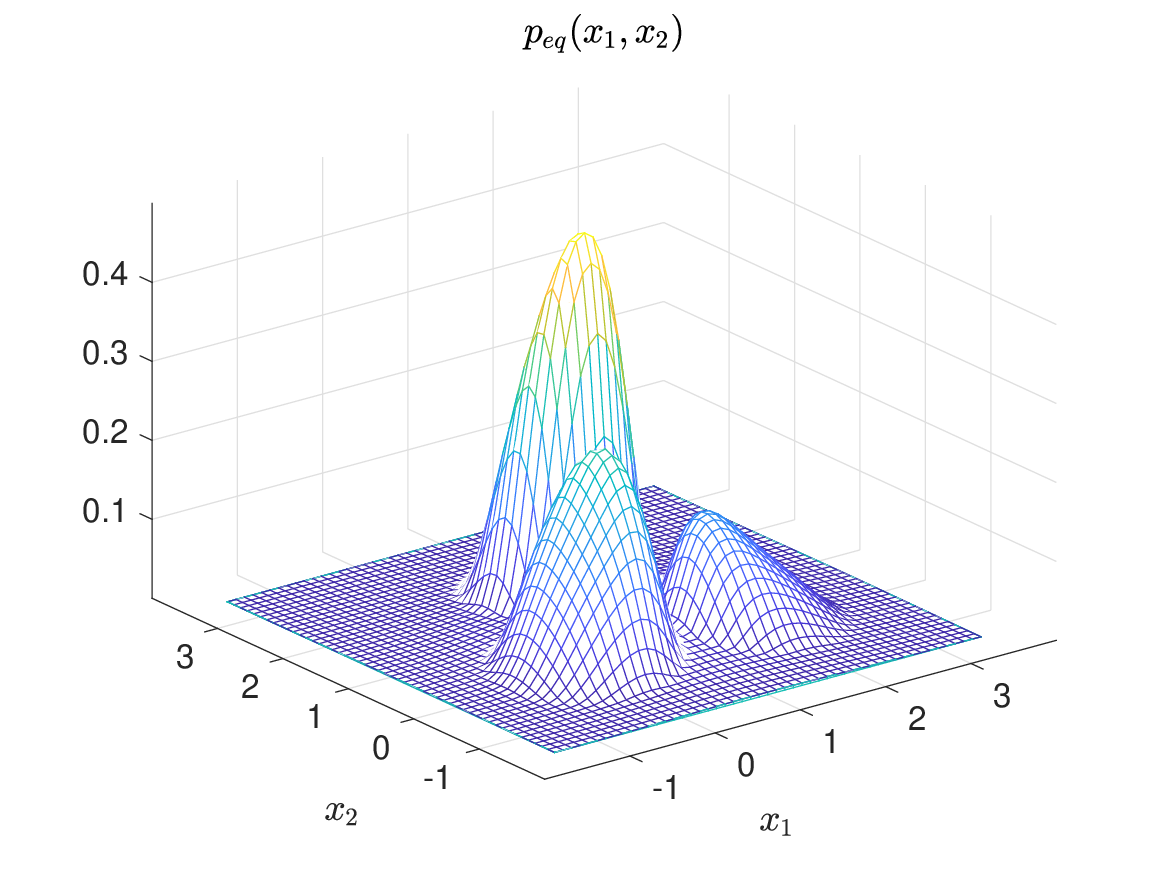}
    \includegraphics[scale=0.41]{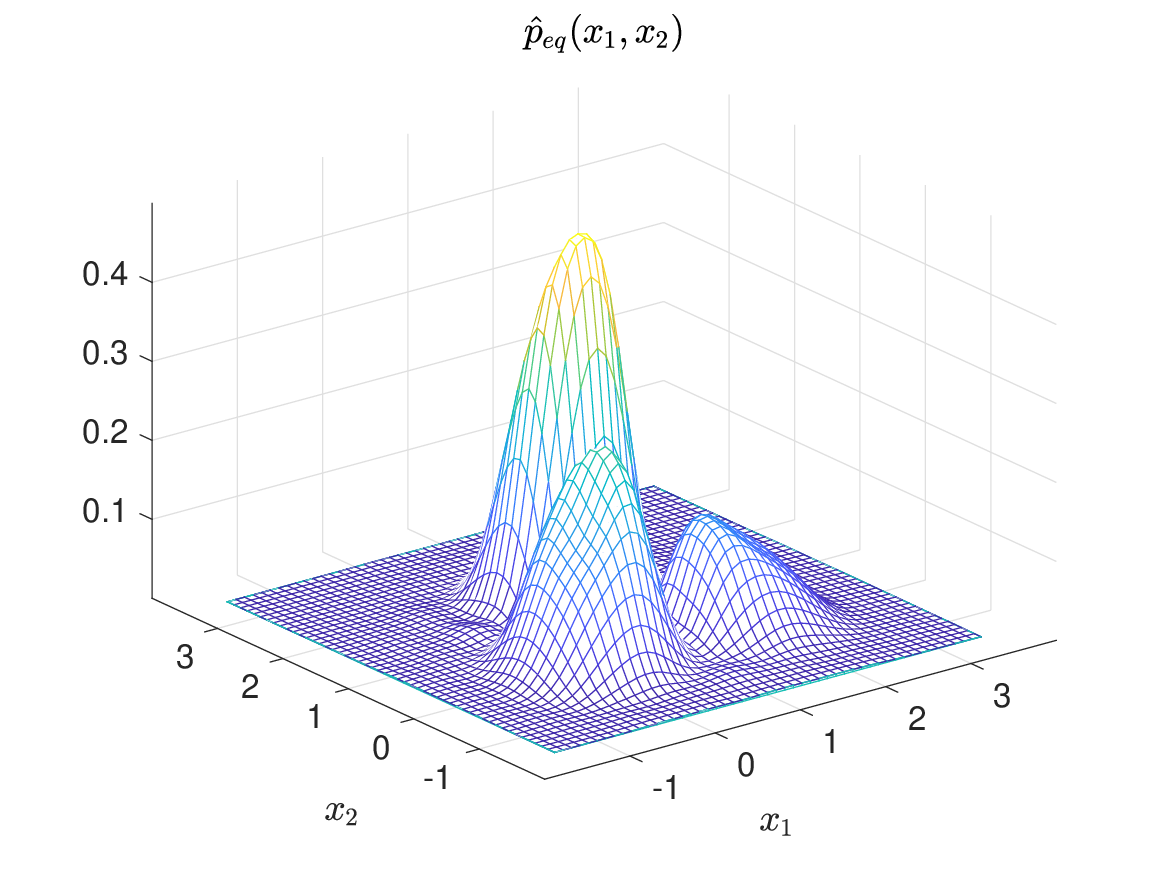}
    \includegraphics[scale=0.41]{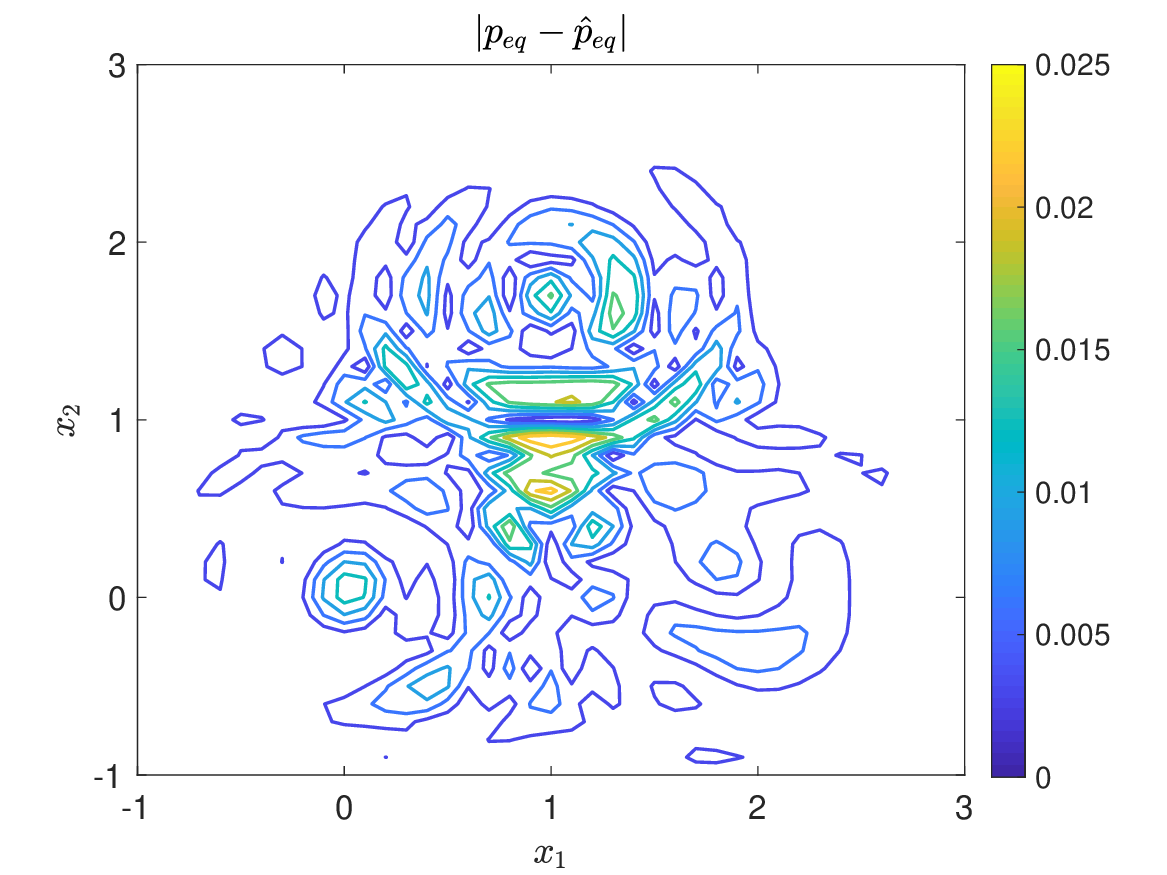}
    \includegraphics[scale=0.41]{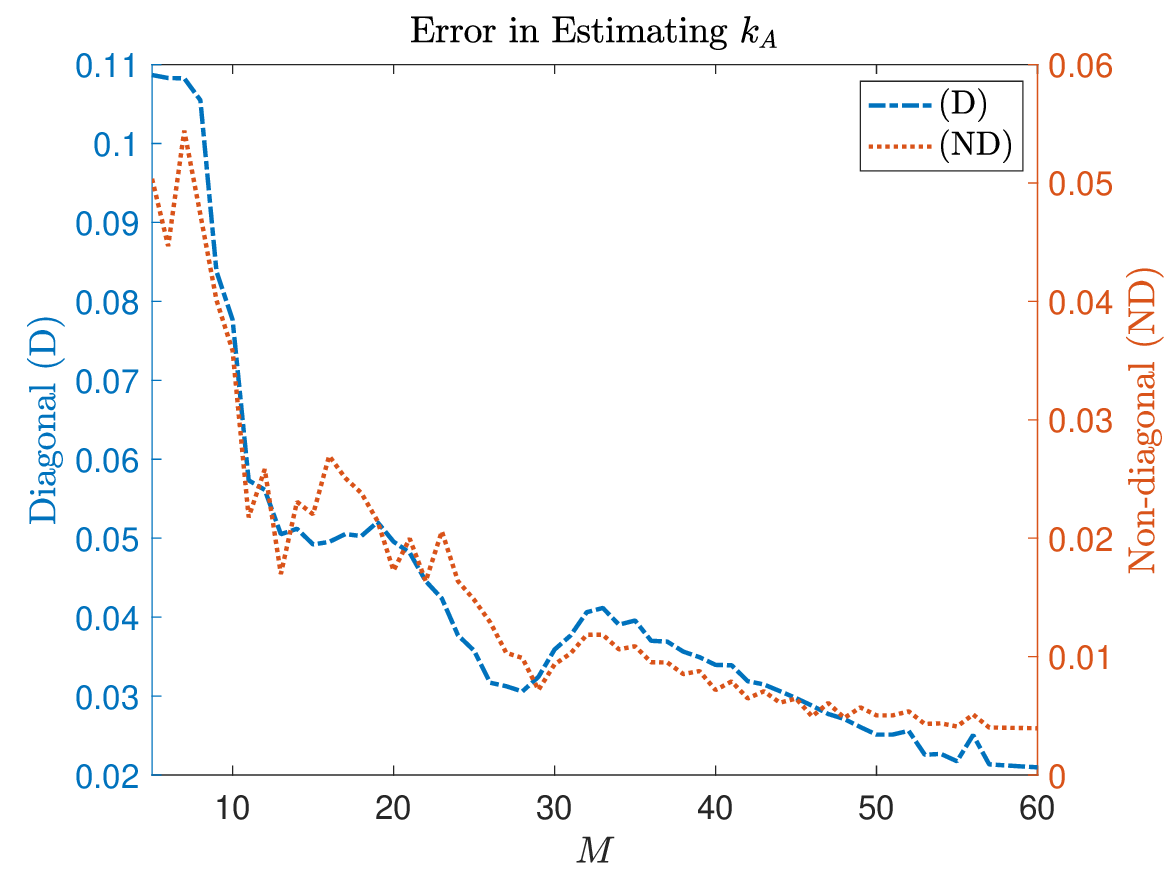}
    \caption{The equilibrium distribution of the triple-well model \eqref{triple} (upper left panel) and its kernel embedding estimate (upper right panel) based on a total of $1\times 10^{7}$ sample. The contour plot (lower left panel) displays the absolute error of the estimate. The error plot (lower right panel) shows the $\ell_{\infty}$-error of the estimates $\hat{k}_{A}$ via kernel embedding linear response.  We separate the diagonal entries (D) from the non-diagonal entries (ND) due to their scale difference.}
	\label{fig:p_tri}
\end{figure}

We first consider a two-dimensional stochastic gradient system as follows, 
\begin{equation}\label{triple}
\dot{\bm{x}} = -\nabla V(\bm{x}) + \sqrt{2k_{B}T} \dot{\bm{W}}_{t}
\end{equation}
where $\bm{W}_t$ is a two-dimensional Wiener process, and $V$, similar to the model in \cite{Hannachi:01}, is a triple-well potential,
\begin{equation}\label{potential_well}
\begin{split}
V(\bm{x})  = & -v\left(x_{1}^{2}+x_{2}^{2}\right)-(1-\gamma)v \left((x_{1}-2a)^{2}+x_{2}^{2} \right)-(1+\gamma)v \left((x_{1}-a)^{2}+(x_{2}-\sqrt{3}a)^{2}\right) \\ &+0.8\left[(x_{1}-a)^{2}+(x_{2}- {a}/{\sqrt{3}})^{2}\right],
\end{split}
\end{equation}
with
\begin{equation}
v(z)=10 \exp\left(\frac{1}{z^{2}-a^{2}}\right)\cdot \chi_{(-a,a)}(z),\nonumber
\end{equation}
where $\chi_{(-a,a)}$ denotes the characteristic function over the interval $(-a,a)$. The parameter $\gamma \in (0,1)$ in \eqref{potential_well} indicates the depth of the three wells. The additional quadratic term $0.8[(x_{1}-a)^{2}+(x_{2}-a/\sqrt{3})^{2}]$ in the triple-well potential \eqref{potential_well} is a smooth retaining potential (also known as a confining potential). The triple-well model \eqref{triple}, as a stochastic gradient system \cite{Pav_book:14}, yields an equilibrium distribution given by
\begin{equation} \label{triple_peq}
p_{eq}(\bm{x}) \propto \exp\left({-\frac{V(\bm{x})}{k_{B}T}}\right), \quad \bm{x}\in \BR^{2}.
\end{equation}
To derive a linear response operator, we consider an external forcing that is constant in $\bm{x}$. Subsequently, the perturbed dynamics is given by,
\begin{equation*}
\dot{\bm{x}} = -\nabla V(\bm{x})+ f(t) \delta  + \sqrt{2k_{B}T} \dot{\bm{W}}_{t},
\end{equation*}
where $|\delta| \ll 1$. If we select $A(\bm{x}):=\bm{x}$ as the observable, the corresponding linear response operator, as a $2$-by-$2$ matrix-valued function given by \eqref{eq:lin_oper}, reads
\begin{equation}\label{eq:lin_oper_tri}
k_{A}(t)= -\mathbb{E}_{p_{eq}}\left[\bm{x}(t)\otimes \nabla \log(p_{eq}(\bm{x}(0)))\right]= \frac{1}{k_{B}T}\mathbb{E}_{p_{eq}}\left[\bm{x}(t)\otimes \nabla V(\bm{x}(0))\right].
\end{equation}
For the quadratic retaining potential in \eqref{potential_well}, $p_{eq}$ in \eqref{triple_peq} yields a Gaussian tail. Therefore, we apply the kernel embedding linear response to the linear response operator $k_{A}(t)$ in \eqref{eq:lin_oper_tri} based on the Hermite polynomials. Let $\hat{p}_{eq}$ denote the kernel embedding estimate of $p_{eq}$ in \eqref{triple_peq}, and by Eq. \eqref{eq:KMLE1}, the kernel embedding linear response operator, as an estimate of $k_{A}(t)$ in \eqref{eq:lin_oper_tri}, is defined as
\begin{equation}\label{eq:KMLE_tri}
\hat{k}_{A}(t)= -\mathbb{E}_{p_{eq}}\left[\bm{x}(t)\otimes \nabla \log\big(\hat{p}_{eq}(\bm{x}(0))\big)\right].
\end{equation}
In the numerical experiment, we set $(a,k_B T,\gamma)=(1,1.5,0.25)$. To generate the data from \eqref{triple}, we apply the weak trapezoidal method  \cite{Anderson:09} with step size $h= 1\times 10 ^{-3}$, followed by  a $1/5$-subsample. Figure~\ref{fig:p_tri} shows the order-$60$ (we take $\beta=1$ in \eqref{eq:rep_for}, and $M = 60$ in \eqref{eq:kme1}) kernel embedding estimates $\hat{p}_{eq}$ compared with the $p_{eq}$ in \eqref{triple_peq}.
The estimate $\hat{p}_{eq}$ captures the well structure of the model with a decent accuracy. In this figure, we also see that the error in estimating $k_A$ decreases as a function of $M$.

\begin{figure}[ht!]
	\centering
	\includegraphics[scale=0.4]{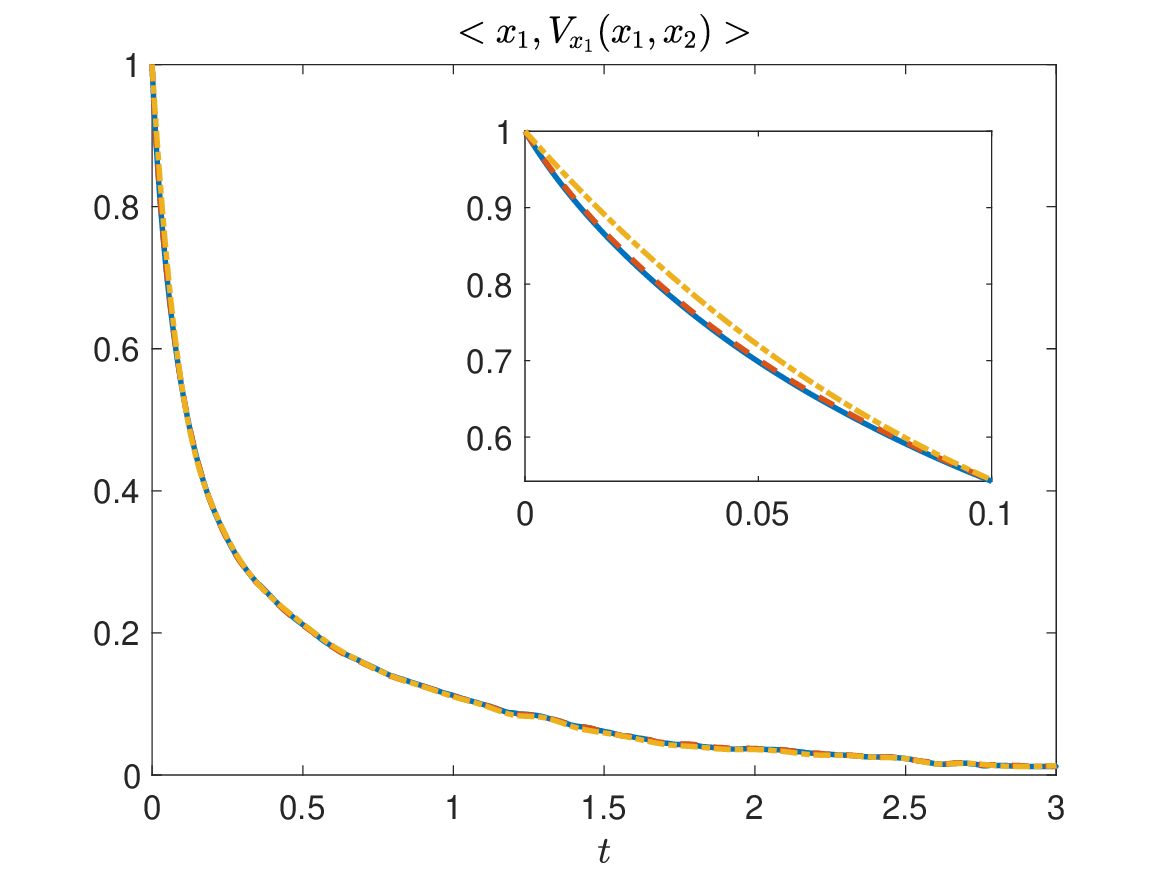}
    \includegraphics[scale=0.4]{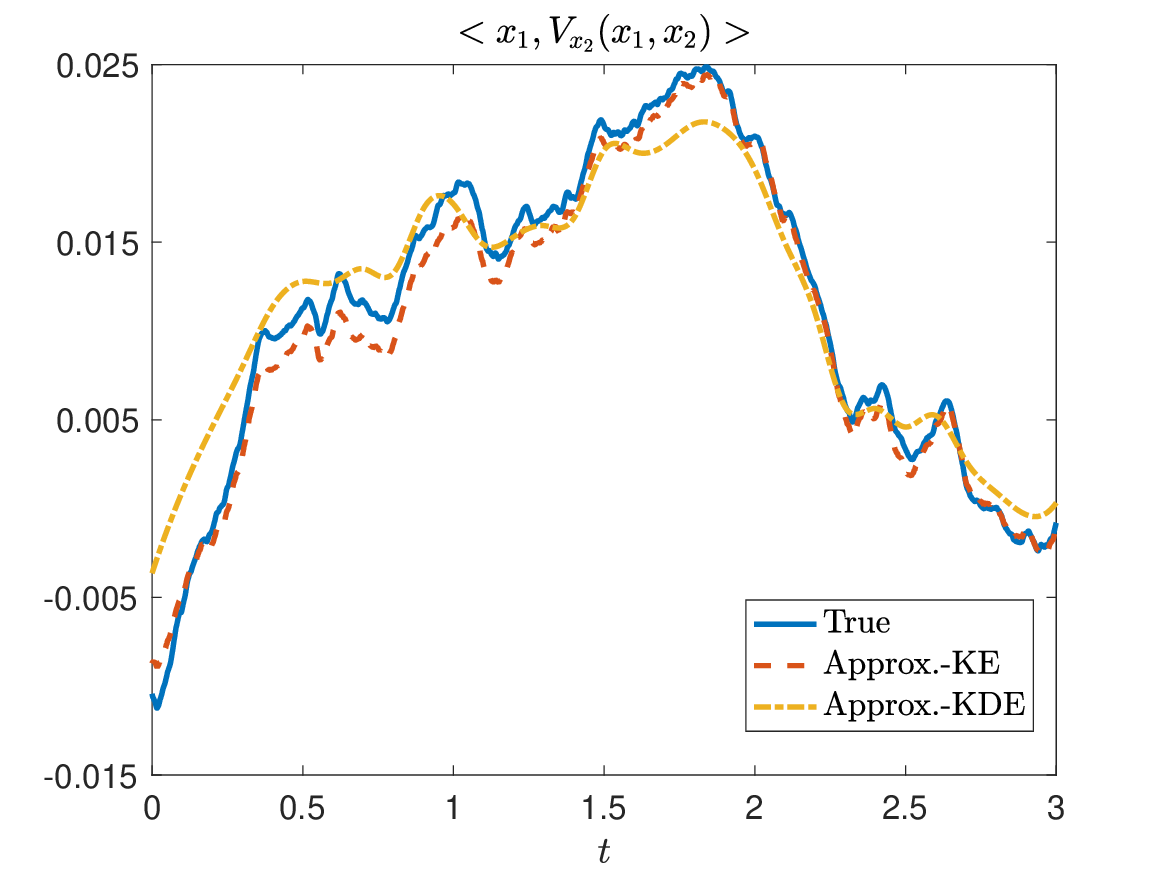}
    \includegraphics[scale=0.4]{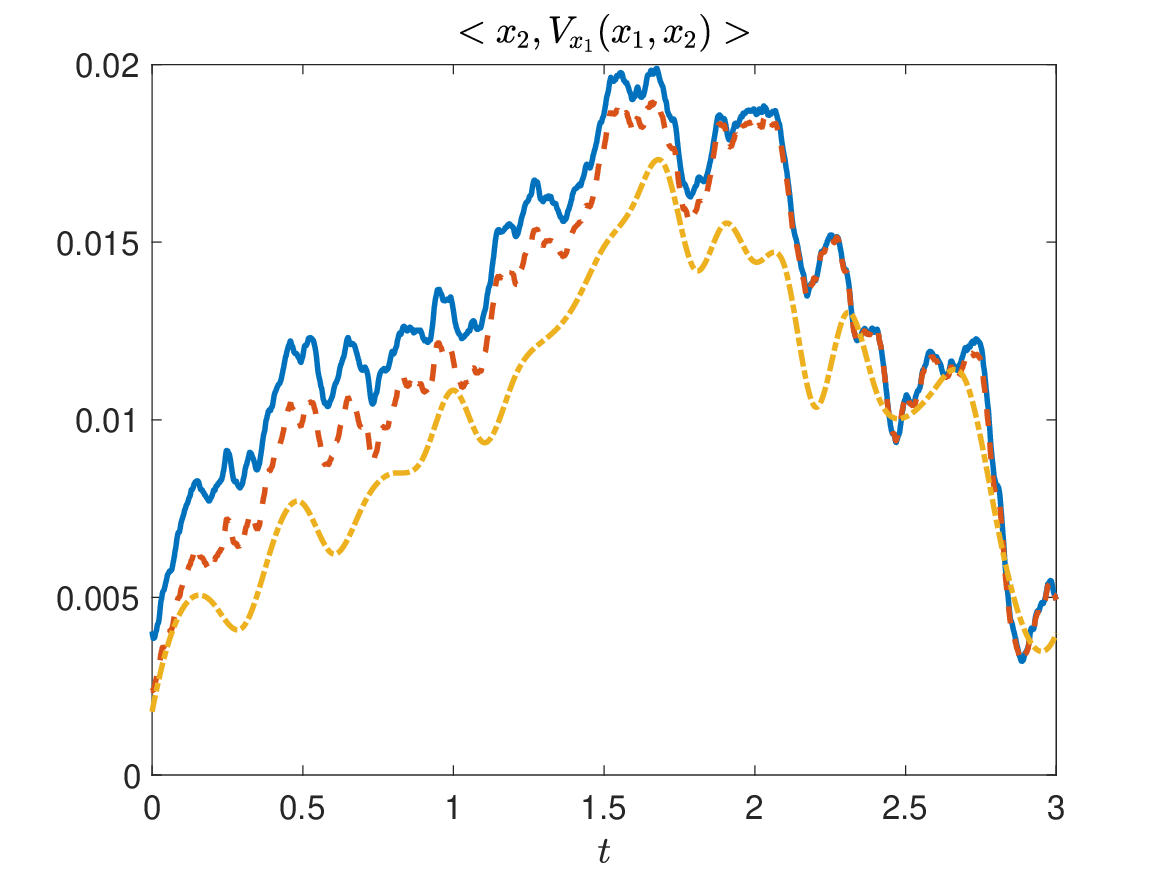}
    \includegraphics[scale=0.4]{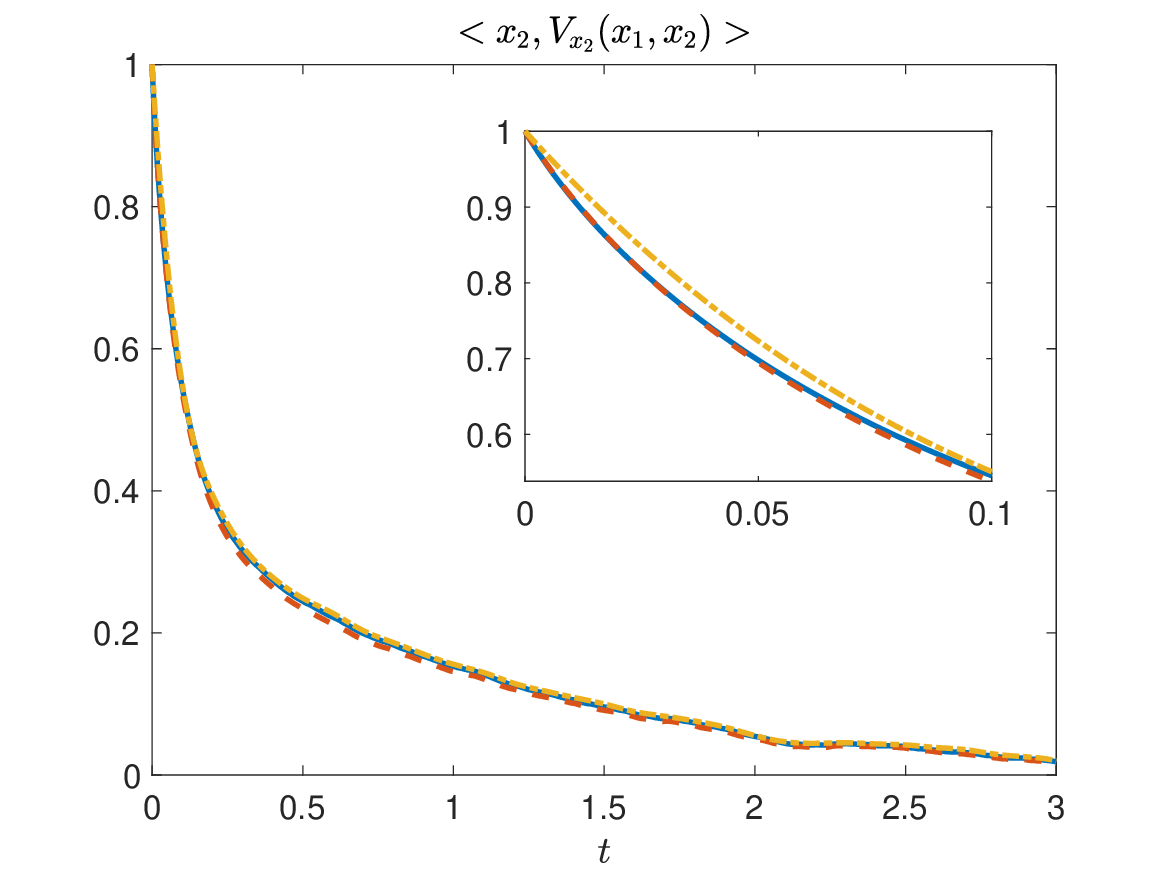}
	\caption{The linear response operator $k_{A}$ in \eqref{eq:lin_oper_tri} (blue solid) and the corresponding estimates $\hat{k}_{A}$ in \eqref{eq:KMLE_tri} via kernel embedding linear response (red dash) and KDE (yellow dot-dash).  For the two-point statistics,  both $k_{A}$ and $\hat{k}_{A}$ are computed via Monte-Carlo.  The diagonal entries of $k_{A}$ and $\hat{k}_{A}$ are normalized so that they share the same initial value $1$. Two insert figures are added to the diagonal entries to show the details of the estimates at the initial stage.}
	\label{fig:KELR_tri}
\end{figure}

Figure~\ref{fig:KELR_tri} compares the linear response estimates in \eqref{eq:lin_oper_tri}, obtained from the kernel embedding formulation and KDE, with the true linear response. Notice that kernel embedding with $M=60$ produces more accurate estimates compared to the KDE-based estimator.

Since the unperturbed dynamics \eqref{triple} is a stochastic gradient system, the off-diagonal entries of the linear response operator $k_{A}(t)$ in \eqref{eq:lin_oper_tri} should be zero at $t=0$ due to the equipartition of the energy in statistical mechanics \cite{HLZ:19}. However, the empirical estimates of $k_{A}$ reported in Figure~\ref{fig:KELR_tri} shows a slight deviation from the exact value. This error can be attributed to the fact that our data is generated from a numerical discretization, which introduces an error in the equilibrium density and the time correlation \cite{leimkuhler2016computation}. Nevertheless, it is quite remarkable that the kernel embedding estimate still captures the true off-diagonal component profile, considering that the true response is relatively small itself. Here, the difference between the KDE-based estimates and the true response is more profound.

\subsection{A Langevin equation with Morse potential} \label{sec:Lan}

\begin{figure}[ht]
	\centering
	\includegraphics[scale=0.41]{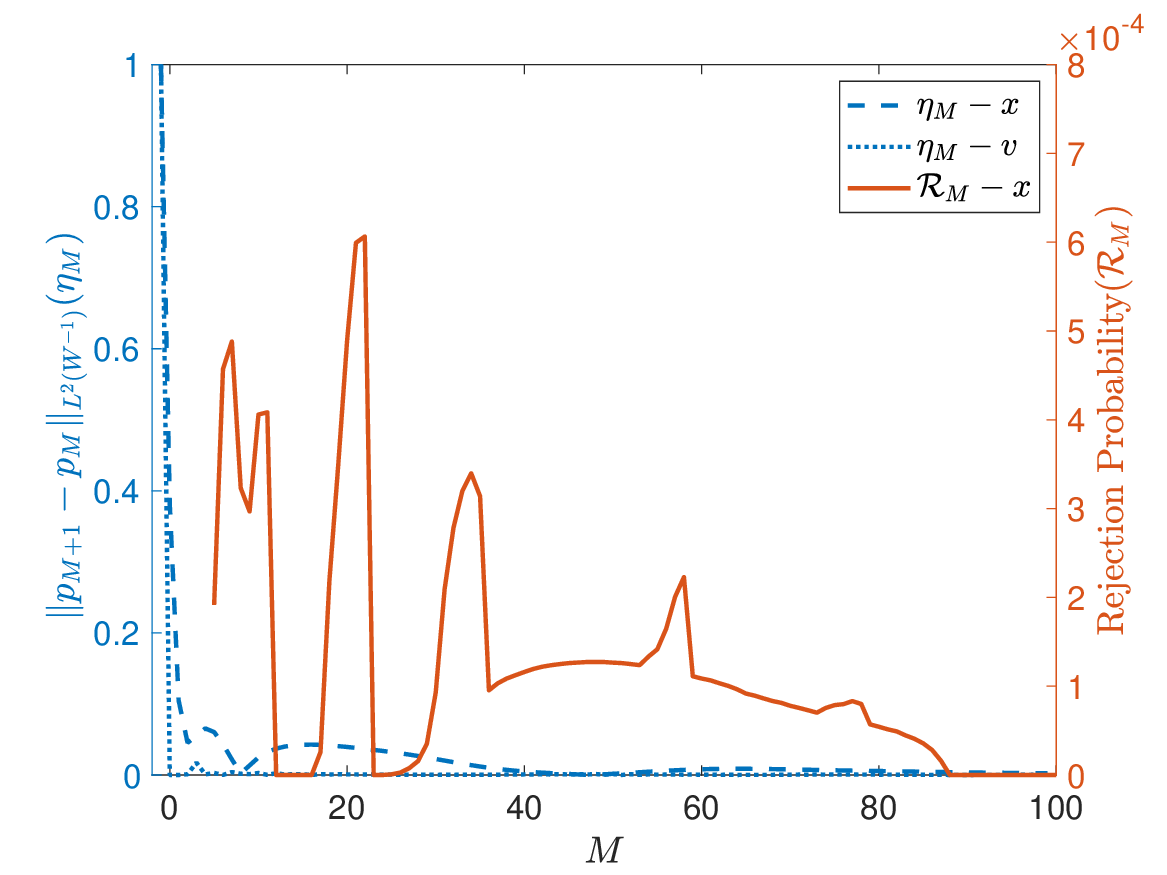}
	\includegraphics[scale=0.41]{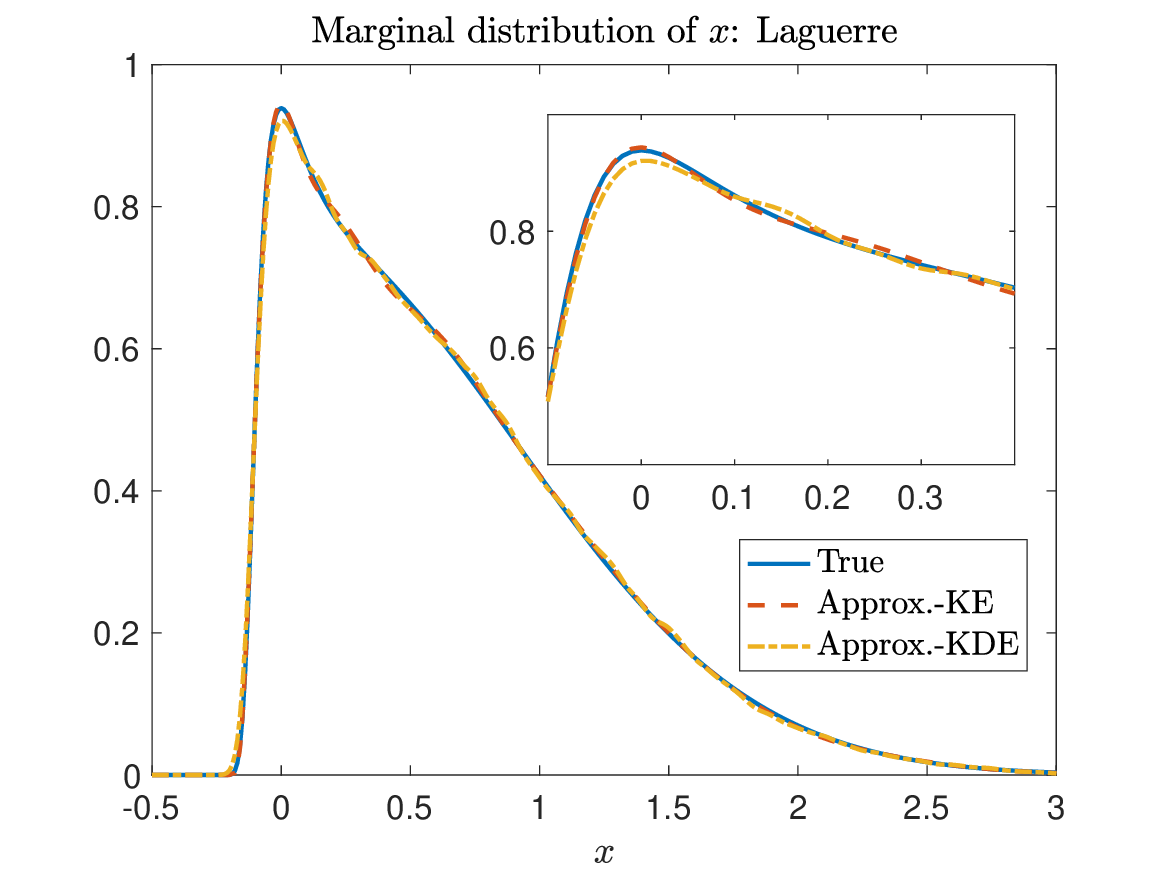}
	\caption{Left panel: $\eta_M$ for the Gaussian marginal density of variable $v$ as a function of $M=M_2$ (dotted blue line). In the same panel, we also show $\eta_M$ for the marginal density of variable $x$ (dashes blue) and the rejection probability $\mathcal{R}_M$ (solid red) as functions of $M=M_1$ for a fixed $M_2=0$. Right panel: The marginal distribution of $x$ (left) of the Langevin dynamics \eqref{Lan_sys} at equilibrium. The kernel embedding estimate uses Laguerre polynomials with $M=90$. All of the results in this picture is based on a total of $N=10^{7}$ samples.}
	\label{fig:phat_Lan}
\end{figure}

For the second example, we consider a Langevin dynamics in statistical mechanics, describing the dynamics of a particle driven by a conservative force, a damping force, and a stochastic force. In particular, we choose the conservative force based on the Morse potential
\begin{equation}
U(x)= U_0(a(x-x_0)), \quad U_0(x)= \epsilon(e^{-2x}-2e^{-x} + 0.03 x^2), \nonumber
\end{equation}
where the last quadratic term in $U_{0}$, similar to the last term in the triple-well potential  \eqref{potential_well}, acts as a retaining potential, preventing the particle from moving to infinity. For this one-dimensional model, we rescale the mass to unity, and write the dynamics as follows
\begin{equation}\label{Lan_sys}
\begin{cases}
\dot{x}=v,  \\
\dot{v}=-U'(x)-\gamma v+\sqrt{2\gamma k_{B}T} \dot{W}_t, \\
\end{cases}
\end{equation}
where $\dot{W}_t$ is a white noise. The smooth retaining potential $U(x)$ guarantees the ergodicity of the Langevin system in \eqref{Lan_sys} (see Appendix C of \cite{HLZ:19} for the proof, which is an application of the result in \cite{Mattingly:02}). Namely, there is an equilibrium distribution (also known as the Gibbs measure) $p_{eq}(x,v)$, given by
\begin{equation}\label{eqdis} 
p_{eq}(x,v)\propto \exp\left[-\frac{1}{k_{B}T}\left(U(x)+\frac{1}{2}v^{2}\right)\right].
\end{equation}
To derive a linear response operator in \eqref{eq:lin_oper}, we introduce a constant external forcing together with a constant external velocity field, and the corresponding perturbed system is given by
\begin{equation*}
\begin{cases}
\dot{x}=v +\delta f_1 \\
\dot{v}=-U'(x)-\gamma v+\delta f_2 +  \sqrt{2\gamma k_BT} \dot{W}_t. \\
\end{cases}
\end{equation*}
By selecting the observable $A= (x , v)^{\top}$, the resulting linear response operator is given by
\begin{equation}\label{esst_lan}
k_{A}(t) =-\mathbb{E}_{\peq}\left[(x(t),v(t))^{\top} \otimes \nabla \log(p_{eq}(x(0),v(0))) \right] = \frac{1}{k_B T}\mathbb{E}_{\peq} \left[ 
\begin{pmatrix}  x(t)U'(x(0)) & x(t)v(0) \\ v(t)U'(x(0)) & v(t)v(0)  \end{pmatrix}
\right].
\end{equation}
Here, to derive the $2$-by-$2$ two-point statistics matrix in \eqref{esst_lan}, we have used the fact that the variables $x$ and  $v$ are independent at the equilibrium. However, in our estimation problem, we do not assume that the underlying density, $\peq$, is a product of two marginal densities. In particular, adopting the notation in Section~\ref{sec:Her} and \ref{sec:Lag}, we define the order-$(M_1,M_2)$ kernel embedding estimate of $p_{eq}$ in \eqref{eqdis} as
\begin{equation}\label{eq:kme_Lan}
\hat{p}_{eq}(x,v) = \sum_{m_1=0}^{M_1}\sum_{m_2=0}^{M_2} \hat{p}_{m_1, m_2} \ell^{(1)}_{m_1}(x) \psi_{m_2}(v) G(x;1)W(v), \quad \hat{p}_{m_1, m_2} := \int\int \ell^{(1)}_{m_1}(x) \psi_{m_2}(v)p_{eq}(x,v) \td x \td v,
\end{equation}
where $\left\{\ell^{(1)}_{m_1}\right\}$ denotes the normalized Laguerre polynomials with respect to the Gamma distribution $G(x;1) \propto x e^{-x}$; while $\left\{\psi_{m_2}\right\}$ denotes the normalized Hermite polynomials with respect to the Gaussian distribution $W$. The representation in \eqref{eq:kme_Lan} is motivated by the asymmetric structure of the marginal distribution of $x$ (See Figure~\ref{fig:phat_Lan}). To determine $M_1$ and $M_2$, consider the estimation of the marginal distribution of $x$ given by
\begin{equation*}
\hat{\rho}(x) = \int \hat{p}_{eq}(x,v) \td v = \sum_{m_1=0}^{M_1}\left(\sum_{m_2=0}^{M_2}\int \hat{p}_{m_1, m_2} \psi_{m_2}(v) W(v) \td v \right)
\ell^{(1)}_{m_1}(x) G(x;1),
\end{equation*}
which is a linear superposition of the orthogonal basis functions $\{\ell^{(1)}_{m_1}(x) G(x;1)\}$ of order $M_1$. The estimates of the marginal distribution of $v$ follows the order-$M_2$ kernel embedding estimates Eq.\eqref{eq:rep_for} with $\beta =1$. Thus, we may introduce two sub-problems of learning the marginal distributions of $x$ and $v$ based on the Laguerre and Hermite polynomials, respectively.

Figure~\ref{fig:phat_Lan} (left) shows $\eta_{M}$ for the Gaussian marginal density of $v$ (blue dotted line) as a function of $M=M_2$. One can see that even with $M_2=0$, the error already converges, which is not surprising since $v$ is Gaussian at the equilibrium. In particular, using Hermite polynomials, we reach a perfect fit at order-$0$, that is, $M_2 = 0$. With $M_{2} = 0$, suggested by \eqref{eq:kme_Lan}, our density estimation problem reduces to a one-dimensional problem in learning the marginal distribution of $x$. In the same plot, we show the error (in blue dashes) in the estimation of the marginal density of $x$, $\eta_M$, as a function of $M=M_1$ for a fixed $M_2=0$. Notice that the error decreases as a function of $M$. In the same panel, we also plot the rejection probability, $\mathcal{R}_M$ in \eqref{crit1}, as a function of $M=M_1$ for a fixed $M_2=0$. Notice that  while the rejection probability fluctuates, it eventually becomes very small as $M>90$. In the right hand panel, we compare the marginal distribution of $x$ estimated using the Laguerre polynomials via the kernel embedding and the KDE. A close inspection suggests that the kernel embedding produces a more accurate estimate. As for the Gaussian distributed variable, $v$, both estimators are comparably accurate (results are not shown).

With the estimate $\hat{p}_{eq}$ of $p_{eq}$, we define the corresponding kernel embedding linear response, 
\begin{equation}\label{eq:KELR_Lan}
\hat{k}_{A}(t)=-\mathbb{E}_{\peq}\left[(x(t),v(t))^{\top} \otimes \nabla \log(\hat{p}_{eq}(x(0),v(0))) \right],
\end{equation}
as the estimates of $k_{A}(t)$ in \eqref{esst_lan}, which does not rely on the independence of $x$ and $v$ at the equilibrium.

\begin{figure}[ht!]
	\centering
	\includegraphics[scale=0.4]{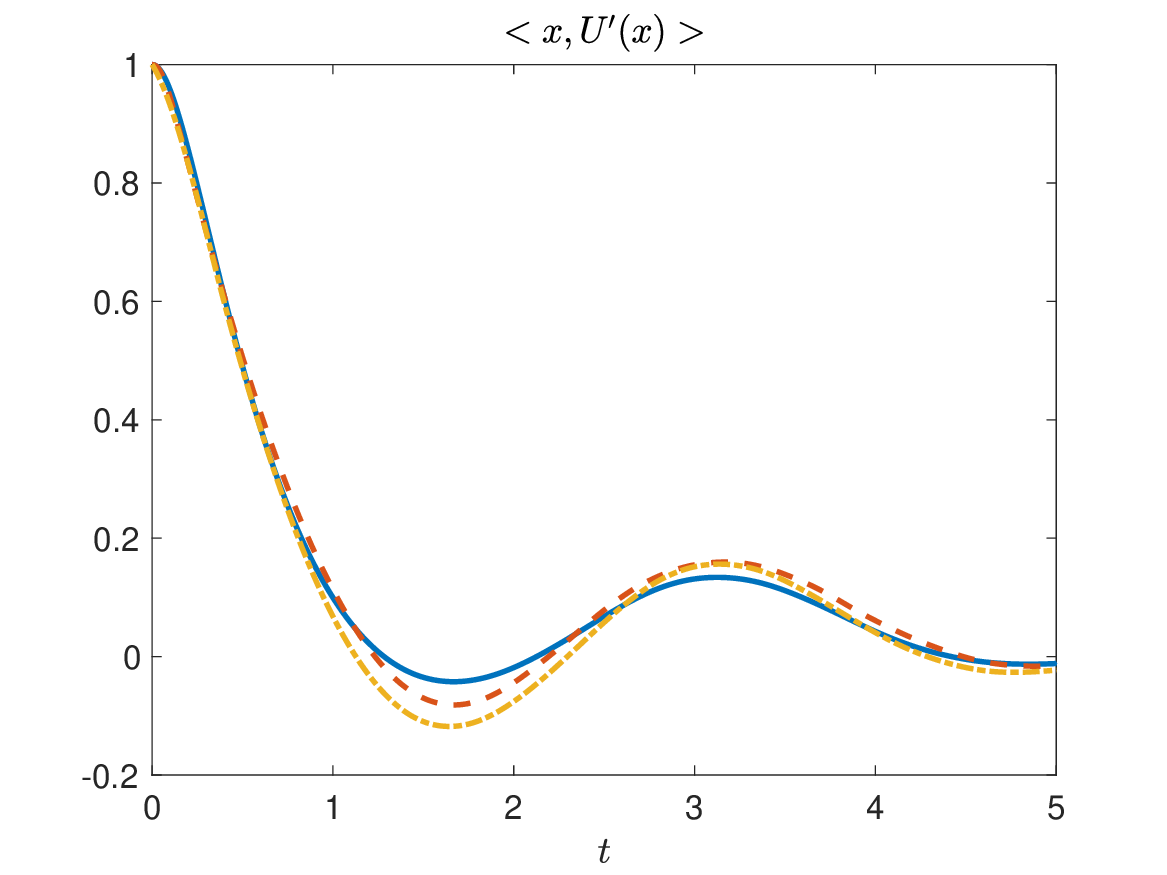}
    \includegraphics[scale=0.4]{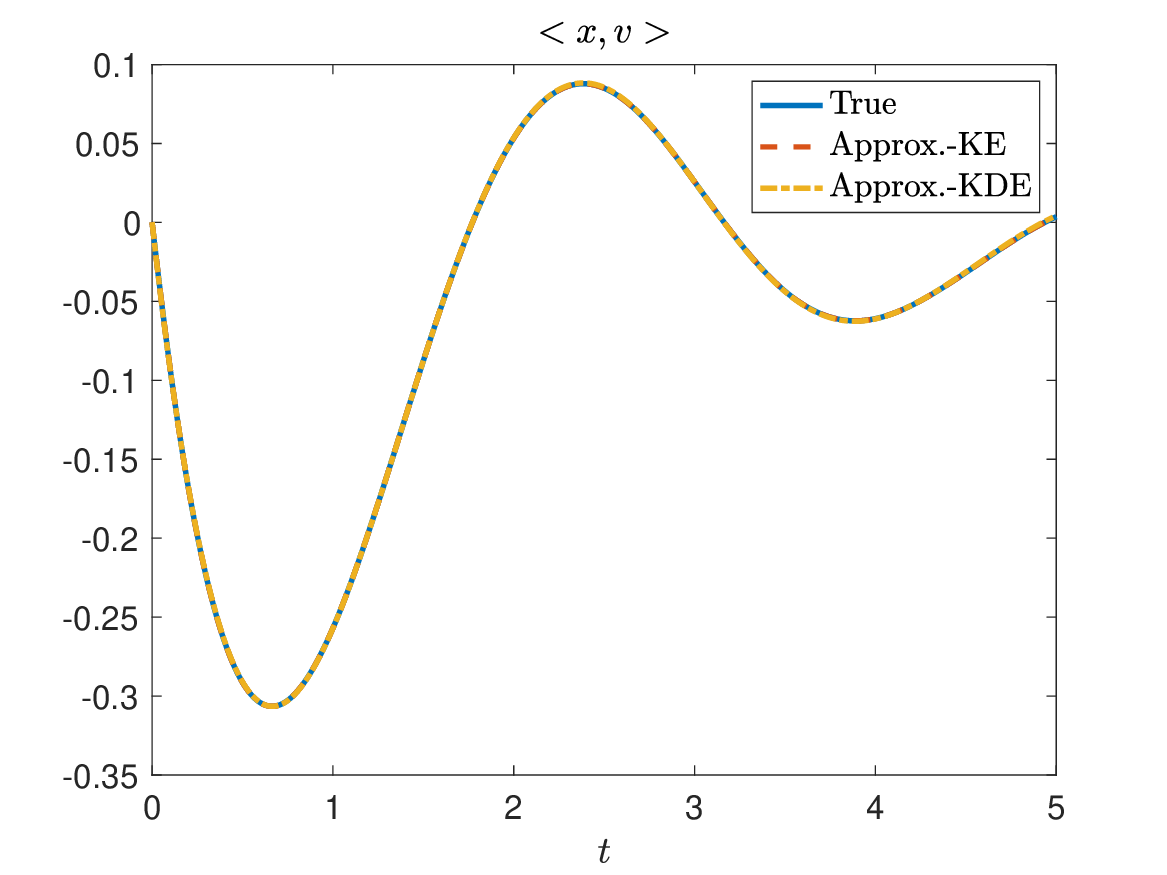}
    \includegraphics[scale=0.4]{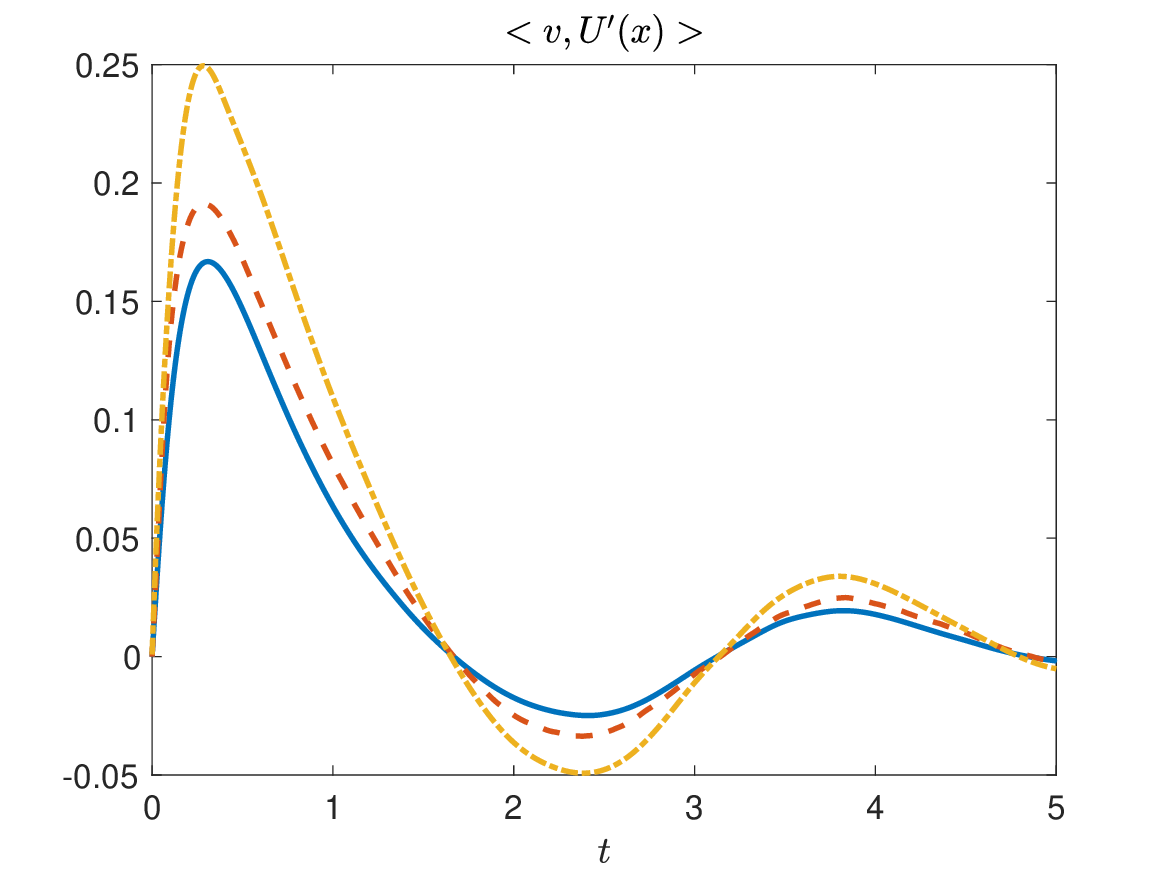}
    \includegraphics[scale=0.4]{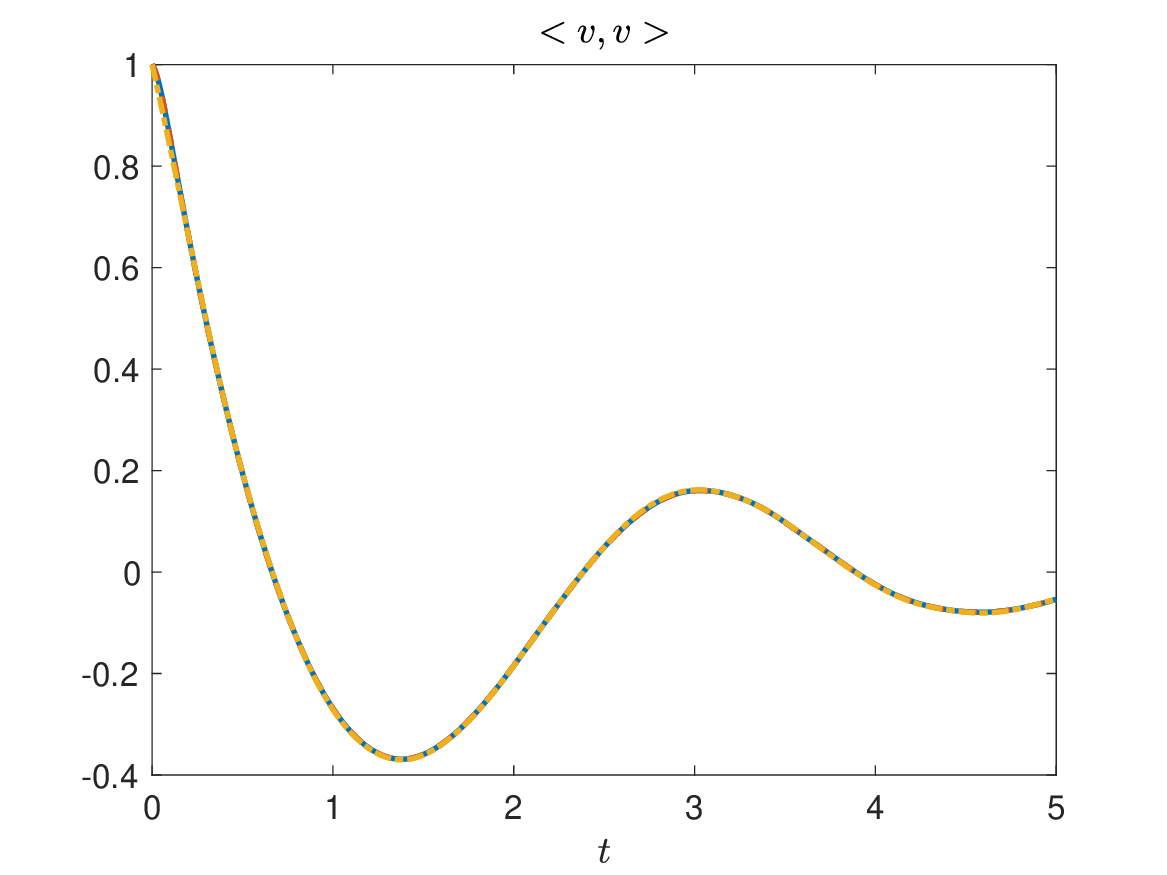}
	\caption{The linear response operator $k_{A}$ in \eqref{esst_lan} (blue solid) and the corresponding estimates $\hat{k}_{A}$ in \eqref{eq:KELR_Lan} via kernel embedding linear response (red dash) and KDE (yellow dot-dash).  All the statistics are computed via Monte-Carlo.  Similar to Figure~\ref{fig:KELR_tri}, the diagonal entries of $k_{A}$ and $\hat{k}_{A}$ are normalized so that they share the same initial value $1$. The $(1,2)$ and $(2,2)$ components reach perfect fits for both methods since $v$ is Gaussian at the equilibrium.}
	\label{fig:KELR_Lan}
\end{figure}

In the numerical test, we set $\gamma = 0.5, k_{B}T = 1.0, \epsilon = 0.2, a = 10$, and $x_0 = 0$.  The data, in the form of a time series, are generated from the model \eqref{Lan_sys} using an operator-splitting method \cite{telatovich2017strong} with step size $h = 1\times 10^{-3}$, followed by a $1/10$-subsample. Figure~\ref{fig:phat_Lan} presents the marginals of the kernel embedding estimates $\hat{p}_{eq}$ in \eqref{eq:kme_Lan} with $M_1=90$ and $M_2=0$. Notice that the domain of the Laguerre polynomials $\ell^{(1)}_{m}$ in \eqref{eq:kme_Lan} is $[0, \infty)$. 
In the computation of the coefficients $\hat{p}_{mn}$ in \eqref{eq:kme_Lan} via Monte-Carlo, we have shifted the observations of $x$ so that they are all positive. Figure~\ref{fig:KELR_Lan} compares the linear response operator $k_{A}$ in \eqref{esst_lan} with its estimates $\hat{k}_{A}$ in \eqref{eq:KELR_Lan} based on the kernel embedding linear response {    and KDE.  We reach perfect fits in $(1,2)$ and $(2,2)$ components of $k_{A}$ since $v$ is Gaussian at the equilibrium.} In our Langevin dynamics, the damping coefficient $\gamma = 0.5$ in \eqref{Lan_sys} is used, which   is in the under-damped regime \cite{HLZ:19}, in which case we have a relatively strong memory effect and the two-point statistics $k_{A}(t)$ in \eqref{esst_lan} exhibits a more complicated behavior. The numerical result suggests that the kernel embedding linear response estimate is more accurate compare to the KDE-based estimate.

We report the elapsed wall-clock time of running the simulations based on MATLAB\textregistered using a desktop computer (equipped with a 3.2GHz quad-core Intel Core i5 processor with 32Gb RAM) in Table~\ref{tab:run_time}. Notice that since the complexity of the KDE method only depends on the sample size $N=10^7$ and $d=2$, the wall-clock time for the KDE estimates in both examples above are similar.
With much fewer basis functions, thanks to the orthogonality, the kernel embedding approach is much faster.  We want to point out that, in practice, there are many ways to reduce the computational cost of KDE. Here, we followed the formulation in \eqref{eq:KDE} for easy implementations. From these two examples, we can see that the kernel embedding implemented with appropriate orthogonal polynomials is more efficient and produces more accurate estimates with a relatively simple tuning procedure as proposed in \eqref{crit1}-\eqref{crit2}.

\begin{table}[ht!]
    \centering
    \begin{tabular}{c|c|c}\hline
    Method     & Number of Basis &  Elapsed Time (s)\\
    \hline \hline
    KDE  (Triple-Well, $N = 1\times 10^7$)    &  $1\times 10^7$  & $1.99 \times 10^4$\\
    \hline
    Kernel Embedding (Triple-Well, $M = 60$) & $1891$ & $1.54 \times 10^3$  \\
    \hline
    Kernel Embedding (Langevin, $M_1 = 90$, $M_2 = 0$) & $91$    & $8.21$\\
    \hline
    \end{tabular}
    \caption{   Elapsed time (based on a desktop computer, equipped with a 3.2GHz quad-core Intel Core i5 processor with 32Gb RAM) of the KDE approach and the kernel embedding approach in computing the linear response statistics.}
    \label{tab:run_time}
\end{table}

\section{Summary and discussion} \label{sec:sum}

In this paper, we considered estimating the equilibrium density of unperturbed dynamics from non-i.i.d observed time series. This problem arises from the estimation of linear response statistics. In particular, we employed a nonparametric density estimator formulated by the kernel embedding of distributions. We chose the corresponding hypothesis space (model) as an RKHS so that the properties of the kernel can be carried over to the functions in the RKHS. To avoid the computational expense that arises using radial type kernels, we considered the ``Mercer-type'' kernels constructed based on the classical orthogonal bases defined on non-compact domains, e.g., $\BR^d$. Here, the orthogonality corresponds to a weighted-$L^2$ space, which naturally provides the coefficients' integral definition in the estimates. In practice, the number of bases involved, thanks to the orthogonality, are much fewer than the sample size. For example, in the test models, the sample size is of order $10^7$, while the number of bases is of order $10^3$ or less.

We used the orthogonal polynomials, assigned with a specific power of the corresponding weight, to build a Mercer-type kernel corresponding to an RKHS. To overcome the difficulty caused by the non-compact domain, we showed that the boundedness of the kernel could be achieved either using the asymptotic behavior of a class of orthogonal polynomials (e.g., Hermite polynomials) or existing identities (e.g., Laguerre polynomials). By studying the orthogonal polynomial approximation in the RKHS setting, we established the estimator's uniform convergence, which justifies using the estimator for interpolation. With this choice of basis representation, we arrived at the well-known polynomial chaos approximation. However, the convergence is valid in both $L^2$-norm sense and uniform-norm sense. An important issue that we addressed is the practical problem of choosing hypothesis space (or polynomials), which is critical for accurate estimations. Using the RKHS formulation, we generalized the theory of $c_0$-universality to guarantee a consistent estimator with a hypothesis space constructed to respect the decaying property of the target function that can be empirically quantified from the available data.

In terms of linear response estimation, we defined the kernel embedding linear response based on the kernel embedding estimate of the equilibrium density. Our study provides practical conditions for the estimator's well-posedness and the well-posedness of the underlying response statistics. Given a well-posed estimator, supported by the theory of RKHS, we provided a theoretical guarantee for the convergence of the estimator to the underlying actual linear response statistics. Finally, since we approximate the coefficients in the estimates using a Monte-Carlo average, we derived a statistical error bound for the density estimation that accounts for the Monte-Carlo averaging over non-i.i.d time series with $\alpha$-mixing property and biases due to the finite basis truncation. This error bound provides a means to understand the kernel's feasibility and limitation with the ``Mercer-type'' kernels (and specifically, polynomial chaos expansion). 

Numerically, we validated the kernel embedding linear response estimator on two stochastic dynamics with known but non-trivial equilibrium densities. In the triple-well model, we explored the effectiveness of the kernel embedding estimate of a two-dimensional target density with Gaussian decay. In the Langevin model, the marginal distribution of the displacement, due to the Morse potential, is asymmetric (the decaying properties are different on two sides). We considered a hypothesis space (RKHS) based on a tensor product of Hermite (for the velocity) and Laguerre (for the displacement) polynomials in constructing the kernel. In both examples, we found that the proposed estimator is computationally more efficient and more accurate compared to the KDE-based linear response estimator.

Overall, the kernel embedding linear response provides a systematic and data-driven approach in computing the linear response statistics without knowing the explicit form of the underlying density function. However, this approach is still subjected to the curse of dimensionality due to the usage of the tensor product.

\section*{Acknowledgments} 
The research of XL was supported under the NSF grant DMS-1819011 and JH was supported under the NSF grant DMS-1854299. The authors thank Bharath Sriperumbudur for helpful discussions on the RKHS.

\appendix

\section{Proof of Proposition~\ref{prop:RKHS_d}} \label{app:A}
In this Appendix, we discuss the proof of Proposition~\ref{prop:RKHS_d} which specifies the RKHS $\mathcal{H}_\beta$ constructed using  the Mercer kernel in \eqref{eq:kernel_d}. 
\begin{enumerate}[\roman{enumi}]
\item 
To begin with, notice that $\forall n \geq 0$, and $\forall \bm{x} \in \BR^{d}$,
\begin{equation}
\left|   f_{n}(\bm{x}) \right|   \leq  \sum_{ \|\vec{m}\|_{1}\leq n } \left| \hat{f}_{\vec m} \Psi_{\beta,\vec{m}}( \bm{x})\right|  \leq \left( \sum_{\vec{m}\geq 0}  \frac{\hat{f}^{2}_{\vec m}}{\lambda_{\vec{m}}} \right)^{\frac{1}{2}}\cdot \left( \sum_{ \vec{m} \geq 0 } \lambda_{\vec{m}} \Psi^{2}_{\beta,\vec{m}}( \bm{x})\right)^{\frac{1}{2}} = \left( \sum_{\vec{m}\geq 0}  \frac{\hat{f}^{2}_{\vec m}}{\lambda_{\vec{m}}} \right)^{\frac{1}{2}} k_{\beta}^{\frac{1}{2}}(\bm{x}, \bm{x}),\nonumber
\end{equation}
and together with the decay rate \eqref{eq:decay_d}, we have, by Eq. \eqref{eq:RKHS_cond1}, $\left \{f_{n} \right\} \subset C_{0}(\BR^{d})$. Following the same idea, for positive integers $n_{2}>n_{1}$, and $\forall \bm{x} \in \BR^{d}$,
\begin{equation*}
\left| f_{n_2}(\bm{x}) - f_{n_1}(\bm{x}) \right| \leq \sum_{ n_{1}<\|\vec{m}\|_{1}\leq n_{2}} \left| \hat{f}_{\vec m} \Psi_{\beta,\vec{m}}( \bm{x})\right|  \leq \left( \sum_{ n_{1}<\|\vec{m}\|_{1}\leq n_{2}}  \frac{\hat{f}^{2}_{\vec m}}{\lambda_{\vec{m}}} \right)^{\frac{1}{2}}k_{\beta}^{\frac{1}{2}}(\bm{x}, \bm{x}) \leq  C_{3}^{\frac{d}{2}}\left( \sum_{ n_{1}<\|\vec{m}\|_{1}}  \frac{\hat{f}^{2}_{\vec m}}{\lambda_{\vec{m}}} \right)^{\frac{1}{2}},
\end{equation*}
which, as a result of Eq. \eqref{eq:RKHS_cond1}, implies that $\left\{f_{n} \right\}$ is a Cauchy sequence in $C_{0}(\BR^{d})$, and $f_{n}$ uniformly converges to a function$f^{*}\in C_{0}(\BR^{d})$ satisfying
\begin{equation*}
f^{*}(\bm{x}) = \sum_{\vec{m}\geq 0}\hat{f}_{\vec{m}} \Psi_{\beta,\vec{m}}(\bm{x}), \quad \forall \bm{x} \in \BR^{d}.
\end{equation*}
In particular,
\begin{equation*}
 \int_{\BR^{d}} \left(f^{*}(\bm{x}) \right)^{2} \bm{W}^{1-2\beta}(\bm{x}) \td \bm{x} = \sum_{\vec{m}\geq 0} \hat{f}^{2}_{\vec{m}} < \infty,
\end{equation*}
that is, $f^{*} \in L^{2}(\mathbb{R}^{d}, \bm{W}^{1-2\beta})\cap C_{0}(\mathbb{R}^{d})$.

\item 
From the first result, we have learned that $\mathcal{H}_{\beta}$ \eqref{eq:RKHS} is a well-defined subspace of $L^{2}(\mathbb{R}^{d}, \bm{W}^{1-2\beta})\cap C_{0}(\mathbb{R}^{d})$. With the positivity of $\{\lambda_{n}\}$, it is straightforward to show that $\langle \cdot,\cdot \rangle$ defines an inner product on $\mathcal{H}_{\beta}$. We now prove $\mathcal{H}_{\beta}$ is closed with respect to the topology induced by the inner product $\langle \cdot, \cdot \rangle$.

Take a Cauchy sequence $\left\{f^{(n)} \right\} \subset \mathcal{H}_{\beta}$ with
\begin{equation*}
         f^{(n)} = \sum_{\vec{m}\geq 0} \hat{f}^{(n)}_{\vec{m}} \Psi_{\beta,\vec{m}}, \quad n =1,2,\dots,
\end{equation*}
that is, $\forall \epsilon >0$, there exists $N>0$ such that $\forall n_{1}, n_{2} >N$, we have $ \left \| f^{(n_{1})} - f^{ (n_{2})} \right\|_{\mathcal{H}_{\beta}}<\epsilon$. Here, $\|\cdot\|_{\mathcal{H}_{\beta}}$ denotes the norm induced by the inner product $\langle \cdot , \cdot\rangle$. Since $\lambda_{n}\rightarrow 0$ as $n\rightarrow \infty$, we are able to show that $\left \{f^{(n)} \right\}$ is also Cauchy in $L^{2}(\mathbb{R}^{d}, \bm{W}^{1-2\beta})$. In particular, let $f^{(n)} \rightarrow f^{\dagger}$ in $L^{2}(\mathbb{R}^{d}, \bm{W}^{1-2\beta})$ with 
\begin{equation} \label{eq:f_dagger}
       f^{\dagger} = \sum_{\vec{m}} f^{\dagger}_{\vec{m}} \Psi_{\beta,\vec{m}}, \quad      \lim_{n\rightarrow \infty} \hat{f}^{(n)}_{\vec{m}} =\hat{f}^{\dagger}_{\vec{m}},  \quad \forall \vec{m} \geq 0.
\end{equation}
By Fatou's Lemma, we have
\begin{equation*}
    \sum_{\vec{m} \geq 0 } \frac{\left(\hat{f}^{\dagger}_{\vec{m}}\right)^{2}}{\lambda_{\vec{m}}} =  \sum_{\vec{m} \geq 0 } \lim_{n\rightarrow \infty} \frac{\left(\hat{f}^{(n)}_{\vec{m}}\right)^{2}}{\lambda_{\vec{m}}} \leq \liminf_{n\rightarrow \infty} \sum_{\vec{m} \geq 0 } \frac{\left(\hat{f}^{(n)}_{\vec{m}}\right)^{2}}{\lambda_{\vec{m}}} = \liminf_{n\rightarrow \infty} \left \|f^{(n)}  \right\|_{\mathcal{H}_{\beta}}^{2} < \infty,
\end{equation*}
which means $f^{\dagger} \in \mathcal{H}_{\beta}$ with
\begin{equation*}
   f^{\dagger}(\bm{x}) = \sum_{\vec{m}} f^{\dagger}_{\vec{m}} \Psi_{\beta,\vec{m}}(\bm{x}), \quad \forall \bm{x} \in \mathbb{R}^{d},
\end{equation*}
is a bounded continuous function as a representative of $f^{\dagger} \in L^{2}(\mathbb{R}^{d}, \bm{W}^{1-2\beta})$.

Finally, we need to show that $f^{(n)} \rightarrow f^{\dagger}$ in $\mathcal{H}_{\beta}$. We first claim that $\forall \epsilon >0$, there exist positive integers $N_{0}$ and $M_{0}$ such that 
\begin{equation}\label{eq:eps_1}
\sum_{\| \vec{m}\|_{1}> M_{0}} \frac{\left( \hat{f}^{(n)}_{\vec{m}}\right)^{2}}{\lambda_{\vec{m}}} \leq \frac{\epsilon}{8}, \quad \forall n > N_{0}.
\end{equation}
Simply notice that
\begin{equation*}
\sum_{\| \vec{m}\|_{1}> M_{0}} \frac{\left( \hat{f}^{(n)}_{\vec{m}}\right)^{2}}{\lambda_{\vec{m}}} = \sum_{\| \vec{m}\|_{1}> M_{0}} \frac{\left( \hat{f}^{(n)}_{\vec{m}}-\hat{f}^{(N_{0})}_{\vec{m}} + \hat{f}^{(N_{0})}_{\vec{m}} \right)^{2}}{\lambda_{\vec{m}}} \leq 2 \left\| f^{(n)}- f^{(N_{0})} \right\|^{2}_{\mathcal{H}_{\beta}} + 2 \sum_{\| \vec{m}\|_{1}> M_{0}} \frac{\left( \hat{f}^{(N_{0})}_{\vec{m}}\right)^{2}}{\lambda_{\vec{m}}},
\end{equation*}
and, by the fact that $\left\{f^{(n)}\right\}$ is Cauchy in $\mathcal{H}_{\beta}$, we first pick $N_{0}$ large enough such that 
\begin{equation*}
2 \left\| f^{(n)}- f^{(N_{0})}\right\|^{2}_{\mathcal{H}_{\beta}}  \leq \frac{\epsilon}{16}, \quad \forall n > N_0.
\end{equation*}
Then, from $f^{(N_{0})}\in \mathcal{H}_{\beta}$, we pick $M_{0}$ large enough such that
\begin{equation*}
 2 \sum_{\| \vec{m}\|_{1}> M_{0}} \frac{\left( \hat{f}^{(N_{0})}_{\vec{m}}\right)^{2}}{\lambda_{\vec{m}}} \leq \frac{\epsilon}{16}.
\end{equation*}
Meanwhile, since $f^{\dagger}\in \mathcal{H}_{\beta}$, there exists a positive integer $M_{1}$ such that
\begin{equation}\label{eq:eps_2}
\sum_{\| \vec{m}\|_{1}> M_{1}} \frac{\left( \hat{f}^{\dagger}_{\vec{m}}\right)^{2}}{\lambda_{\vec{m}}} \leq \frac{\epsilon}{8}.
\end{equation}
Fix $M = \max\{M_{0}, M_{1} \}$, we decompose $\left \| f^{\dagger} - f^{(n)} \right\|^{2}_{\mathcal{H}_{\beta}}$ as follows,
\begin{equation*}
\left \| f^{\dagger}  - f^{(n)}  \right\|^{2}_{\mathcal{H}_{\beta}} =  \sum_{\|\vec{m}\|_{1}\leq M}  \frac{  \left(\hat{f}^{\dagger}_{\vec m} -\hat{f}^{(n)}_{\vec m}   \right)^{2}}{\lambda_{\vec{m}}} + \sum_{\|\vec{m}\|_{1} > M}  \frac{  \left(\hat{f}^{\dagger}_{\vec m} -\hat{f}^{(n)}_{\vec m}   \right)^{2}}{\lambda_{\vec{m}}},
\end{equation*}
and, by the convergence of $\hat{f}^{(n)}_{\vec{m}}$ \eqref{eq:f_dagger}, there exists a positive integer $N_{1}$ such that,
\begin{equation}\label{eq:eps_3}
\sum_{\|\vec{m}\|_{1}\leq M}  \frac{  \left(\hat{f}^{\dagger}_{\vec m} -\hat{f}^{(n)}_{\vec m}   \right)^{2}}{\lambda_{\vec{m}}} < \frac{\epsilon}{2}, \quad  \forall n > N_{1}.
\end{equation}
Finally, take $N = \max\{N_{0}, N_{1}\}$, combining \eqref{eq:eps_1}, \eqref{eq:eps_2}, and \eqref{eq:eps_3}, and $\forall n > N$ we have
\begin{equation*}
\left \| f^{\dagger} - f^{(n)} \right\|^{2}_{\mathcal{H}_{\beta}} < \frac{\epsilon}{2} +  \sum_{\|\vec{m}\|_{1} > M}  \frac{  \left(\hat{f}^{\dagger}_{\vec m} -\hat{f}^{(n)}_{\vec m}   \right)^{2}}{\lambda_{\vec{m}}} \leq \frac{\epsilon}{2} + 2 \sum_{\| \vec{m}\|_{1}> M} \frac{\left( \hat{f}^{\dagger}_{\vec{m}}\right)^{2}}{\lambda_{\vec{m}}} + 2 \sum_{\| \vec{m}\|_{1}> M} \frac{\left( \hat{f}^{(n)}_{\vec{m}}\right)^{2}}{\lambda_{\vec{m}}} < \epsilon,
\end{equation*}
that is, $f^{(n)} \rightarrow f^{\dagger}$ in $\mathcal{H}_{\beta}$.

\item
With all the preparation work we have done so far, the reproducing property \eqref{eq:rep_prop} of $\mathcal{H}_{\beta}$ becomes almost trivial.

For any fixed $ \bm{x}_{0} \in \mathbb{R}^{d}$, $k_{\beta}(\cdot, \bm{x}_{0}) \in \mathcal{H}_{\beta}$. From  $\eqref{eq:kernel_Mercer}$,
\begin{equation*}
\left \langle k_{\beta}(\cdot,\bm{x}_{0}) ,k_{\beta}(\cdot ,\bm{x}_{0}) \right \rangle = \sum_{\vec{m}\geq 0} \lambda_{\vec{m}} \Psi^{2}_{\beta, \vec{m}}(\bm{x}_{0}) = k_{\beta}(\bm{x}_{0},\bm{x}_{0}) < \infty,
\end{equation*}
which means, $\forall f \in \mathcal{H}_{\beta}$, the inner product $\left\langle f ,k_{\beta}(\cdot,\bm{x}_{0}) \right\rangle$ is well-defined $\forall \bm{x}_{0} \in \BR^{d}$. In particular, we have
\begin{equation*}
\left\langle f ,k_\beta(\cdot,\bm{x}_{0}) \right\rangle= \sum_{\vec{m}\geq 0} \hat{f}_{\vec{m}} \Psi_{\beta, \vec{m}}(\bm{x}_{0}) = f(\bm{x}_{0}),
\end{equation*}
which leads to the reproducing property.
\end{enumerate}

\section{Proof of Lemma~\ref{lem:HHF}} \label{app:C}

In this Appendix we discuss the proof of Lemma~\ref{lem:HHF}, which provides a class of bounded kernel based on Laguerre polynomials. As prerequisites, we need the notions of the Bessel and modified Bessel functions.

\begin{figure}[ht]
	\centering
	\includegraphics[scale=0.45]{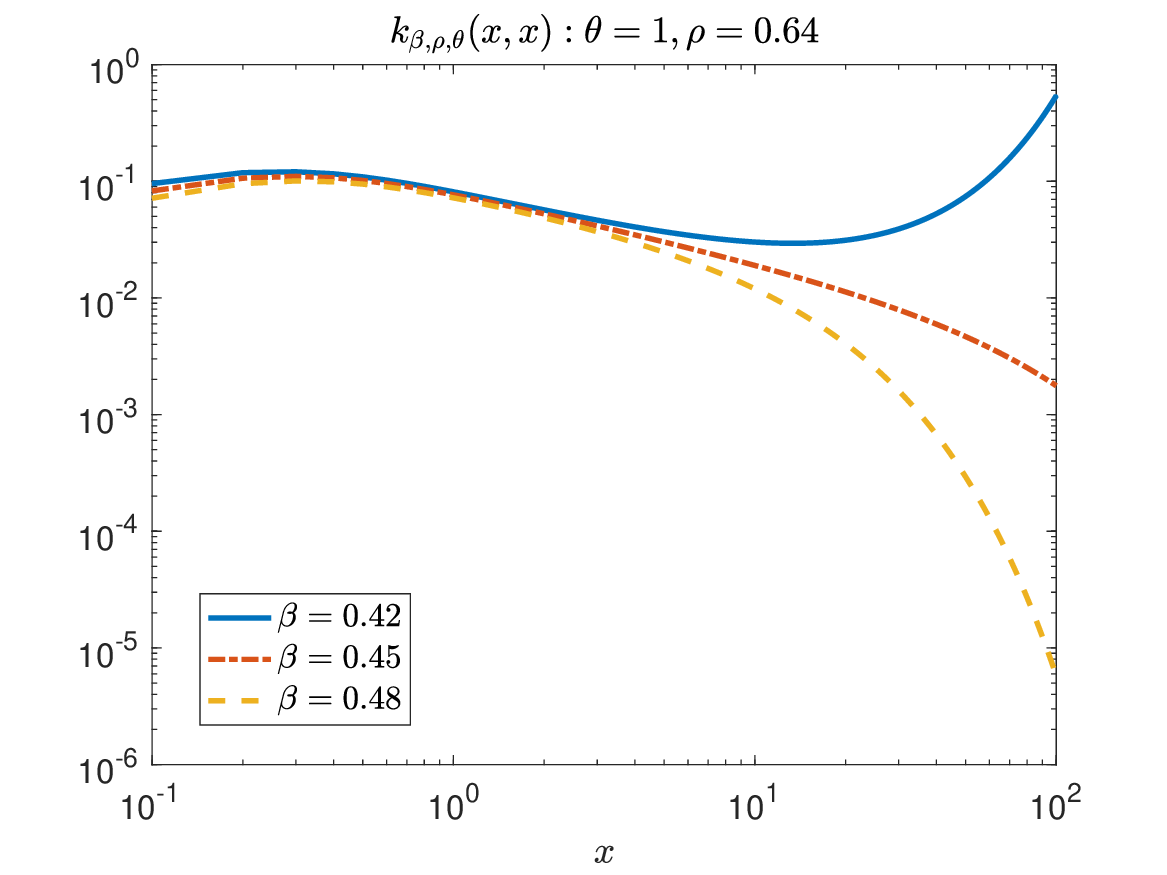}
	\caption{Graph of $k_{\beta, 0.64, 1}(x,x)$ in \eqref{eq:Mercer_Lag} for $\beta = 0.42$ (blue-solid), $\beta = 0.45$ (red-dot-dash), and $\beta = 0.48$ (yellow-dash). Notice that for $\rho = 0.64$, to ensure the boundedness, we need $\beta \geq \frac{0.8}{1+0.8}\approx 0.44$, which is consistent with the numerical results.}
	\label{fig:appc}
\end{figure}

We focus on the (first kind) Bessel functions of integer order \cite{abramowitz1948handbook, szeg1939orthogonal}, which can be defined as
\begin{equation*}
J_{\theta}(z) = \sum_{m = 0}^{\infty} \frac{(-1)^{m}(z/2)^{\theta+2m} }{m! \Gamma(m+\theta+1)}, \quad \theta = 0,1,\dots, \quad z \in \mathbb{C},
\end{equation*}
where $\Gamma(\cdot)$ denotes the Gamma function.  For $\arg z\in \left(-\frac{\pi}{2}, \frac{3\pi}{2}  \right)$, we have the following asymptotic expansion (Eq. (1.71.9) of \cite{szeg1939orthogonal})
\begin{equation} \label{eq:appc1}
i^{\theta} J_{\theta}(-i z) = (2\pi z)^{-\frac{1}{2}} \left[ e^{z} + (-1)^{\theta}ie^{-z}  \right] \left( 1 + O\left(|z|^{-1}\right)  \right), \quad |z|\rightarrow \infty,
\end{equation}
which plays a critical role in proving the boundedness of the kernel. The modified Bessel functions (of the first kind) \cite{abramowitz1948handbook} satisfies,
\begin{equation}\label{eq:appc2}
I_{\theta}(z) = i^{-\theta} J_{\theta}\left(iz \right), \quad \theta = 0,1,\dots
\end{equation}
To prove Lemma~\ref{lem:HHF}, we first recall the Hille-Hardy formula \cite{watson1933notes}, for $\rho\in(0,1)$,
\begin{equation}\label{eq:hhf_1}
\sum_{m=0}^{\infty} \rho^{m} \ell^{(\theta)}_{m}(x) \ell^{(\theta)}_{m}(y) G^{\frac{1}{2}}(x;\theta)G^{\frac{1}{2}}(y; \theta)= \frac{\rho^{-\frac{1}{2}\theta}}{1-\rho} \exp \left[- \frac{1}{2}(x+y)\frac{1+\rho}{1-\rho} \right] I_{\theta}\left( \frac{2\sqrt{xy\rho}}{1-\rho} \right), 
\end{equation}
where $\{\ell^{(\theta)}_{m}\}$ denote the normalized Laguerre polynomials with respect to the Gamma distribution $G(x;\theta) \propto x^{\theta} e^{-x}$ ($\theta = 0,1,\dots$). A modern proof of the Hille-Hardy formula can be found in \cite{al1964operational}. For general $d$-dimensional case, by Eq. \eqref{eq:hhf_1}, the Mercer-type kernel defined  in \eqref{eq:Mercer_Lag} satisfies
\begin{equation*}
\begin{split}
k_{\beta, \rho, \vec{\theta}}(\bm{x},\bm{y}) & = \prod_{i=1}^{d} \sum_{m=0}^{\infty} \rho^{m} \ell^{(\theta_{i})}_{m}(x_{i}) \ell^{(\theta_{i})}_{m}(y_{i}) G^{\beta}(x_{i};\theta_{i})G^{\beta}(y_{i}; \theta_{i}) \\
& = \prod_{i=1}^{d} \frac{\rho^{-\frac{1}{2}\theta_{i}}}{1-\rho} \exp \left[- \frac{1}{2}(x_{i}+y_{i})\frac{1+\rho}{1-\rho} \right] G^{\beta-\frac{1}{2}}(x_{i};\theta_{i})G^{\beta-\frac{1}{2}}(y_{i}; \theta_{i}) I_{\theta_{i}}\left( \frac{2\sqrt{x_{i}y_{i}\rho}}{1-\rho} \right) \\
& = \frac{\rho^{-\frac{1}{2}\|\vec{\theta}\|_{1}}}{(1-\rho)^{d}}\exp \left(-\frac{1+\rho}{2(1-\rho)}\|\bm{x}+\bm{y}\|_{1}\right)G^{\beta-\frac{1}{2}}(\bm{x}; \vec{\theta})G^{\beta-\frac{1}{2}}(\bm{y}; \vec{\theta}) \prod_{i=1}^{d} I_{\theta_{i}} \left( \frac{2\sqrt{x_{i}y_{i}\rho}}{1-\rho}  \right).
\end{split}
\end{equation*}
To prove the boundedness of the kernel $k_{\beta, \rho, \vec{\theta}}(\bm{x},\bm{y})$ in \eqref{eq:Mercer_Lag}, it is enough to show that the one-dimensional kernel $k_{\frac{1}{2}, \rho, \theta}(x,y)$ is bounded $\forall \rho\in(0,1)$ and $\forall \theta \in \{0,1,\dots\}$. Notice that by Eq. \eqref{eq:hhf_1}, \eqref{eq:appc1}, and \eqref{eq:appc2}, we have
\begin{equation*}
\begin{split}
\left|  k_{\frac{1}{2}, \rho, \theta}(x,x) \right| &\propto \exp \left(- x\frac{1+\rho}{1-\rho} \right) \left| I_{\theta}\left( \frac{2x\sqrt{\rho}}{1-\rho} \right)\right|  = \exp \left(- x\frac{1+\rho}{1-\rho} \right) \left| J_{\theta}\left( \frac{2x\sqrt{\rho}}{1-\rho}i \right)\right| \\
& \propto \exp \left(- x\frac{1+\rho}{1-\rho} \right) z^{-\frac{1}{2}}\sqrt{\exp(z)^{2}+ \exp(-z)^{2}} \left( 1 + O\left(|z|^{-1}\right)  \right), 
\end{split} 
\end{equation*}
where $z =  \frac{2x\sqrt{\rho}}{1-\rho}>0$. Thus, for $x$ large, $\left|  k_{\frac{1}{2}, \rho, \theta}(x,x) \right|$ yields the following asymptotic expansion,
\begin{equation}\label{eq:appc3}
\begin{split}
\left|  k_{\frac{1}{2}, \rho, \theta}(x,x) \right| & \propto \exp \left(- x\frac{1+\rho}{1-\rho} \right) \exp \left( \frac{2x\sqrt{\rho}}{1-\rho}\right) \left( O\left(x^{-\frac{1}{2}}\right) + O\left(x^{-\frac{3}{2}}\right)  \right), \\
& = \exp\left(-\frac{1-\sqrt{\rho}}{1+\sqrt{\rho}}x\right) \left( O\left(x^{-\frac{1}{2}}\right) + O\left(x^{-\frac{3}{2}}\right)  \right) \rightarrow 0,
\end{split}
\end{equation}
as $x\rightarrow +\infty$ since $\frac{1-\sqrt{\rho}}{1+\sqrt{\rho}}>0$. Thus, the kernel $k_{\beta, \rho, \vec{\theta}}(\bm{x},\bm{y})$ in \eqref{eq:Mercer_Lag} is bounded $\forall \rho\in (0,1)$, $\beta \geq \frac{1}{2}$, and $\theta \in \{0,1,2\dots\}$. In particular, notice that $\frac{1-\sqrt{\rho}}{1+\sqrt{\rho}}\rightarrow 0$ as $\rho \rightarrow 1^{-}$, the expansion \eqref{eq:appc3} suggests that $\frac{1}{2}$ is indeed the lower bound for $\beta$ to ensure the boundedness $\forall \rho\in (0,1)$. On the other hand, for a fix $\rho \in (0,1)$, the kernel $k_{\beta, \rho, \vec{\theta}}(\bm{x},\bm{y})$ is bounded for $\beta \geq \frac{\sqrt{\rho}}{1+\sqrt{\rho}}$. Figure~\ref{fig:appc} evaluates the kernel $k_{\beta, \rho, \theta}(x,x)$ under $\rho = 0.64$, $\theta= 1$, and three different values of $\beta$. The graph does support the asymptotic expansion \eqref{eq:appc3}.

\section{Proof of Proposition~\ref{prop:mixing}} \label{app:B}

Here we discuss the proof of Proposition~\ref{prop:mixing}, which provides an error bound for the empirical kernel embedding estimate $f_{M,N}$ in \eqref{eq:kme2}. To begin with, let's briefly review the concept of $\alpha$-mixing process and the corresponding Davydov's covariance inequality \cite{Davydov:68}. Given a stationary process $\{X_{n}\}_{n=1}^{\infty}$, let $\mathcal{A}^{b}_{a}$ ( $a<b \leq \infty$) denote the $\sigma$-algebra generated by $\{X_{n}, a\leq n \leq b\}$. We say $\{X_{n}\}$ is $\alpha$-mixing if
\begin{equation*}
\alpha(k): = \sup_{n} \sup_{A\in \mathcal{A}^{n}_{1}, B\in \mathcal{A}^{\infty}_{n+k}} |\mathbb{P}(AB) -\mathbb{P}(A)\mathbb{P}(B)|  \rightarrow 0, \quad k \rightarrow \infty,
\end{equation*}
where $\alpha(k)$ is known as the $\alpha$-mixing coefficient. Such condition imposes restrictions on $X_{n}$ amounting to a weak interdependence of the beginning and the end of the process. Davydov's covariance inequality states that,
\begin{equation}\label{eq:Dav}
    \left |\operatorname{cov}(h(X_{1}), h(X_{1+k})) \right| \leq 12 \alpha(k)^{1/r} \left(\mathbb{E}\left [|h(X_{1})|^{p} \right]  \right)^{1/p}\left(\mathbb{E} \left[|h(X_{1})|^{q} \right]  \right)^{1/q},
\end{equation}
where $1/r + 1/p + 1/q = 1$ given $\mathbb{E}\left[|h(X_{1})|^{p}\right]$ and $\mathbb{E}\left[|h(X_{1})|^{q}\right]$ exist.

To begin with, by the orthogonality of $\{\Psi_{\beta, \vec{m}}\}$ over $L^{2}(\BR^{d},\bm{W}^{1-2\beta})$, we split the error $\left \| f - f_{M,N}  \right\|^{2}_{L^{2}(\bm{W}^{1-2\beta})}$ into two parts,
\begin{equation*}
\left \| f - f_{M,N}  \right\|^{2}_{L^{2}(\bm{W}^{1-2\beta})} = \sum_{\|\vec{m}\|_{1}\leq M} \left ( \hat{f}_{\vec{m}} - \hat{f}_{\vec{m},N}\right)^{2} + \sum_{\|\vec{m}\|_{1} > M} \hat{f}_{\vec{m}}^2,
\end{equation*}
where the first term on the right-hand side can be interpreted as the ``estimation error'' introduced by approximating the coefficient $\hat{f}_{\vec{m}}$ using a Monte-Carlo sum \eqref{eq:kme_MC}; while the second term is the ``approximation error" (or bias) caused by the truncation.

For the bias term, with $f \in \mathcal{H}_{\beta, \rho}$ and Eq. \eqref{eq:RKHS_exp}, we have
\begin{equation}\label{eq:decay_bias}
 \sum_{\|\vec{m}\|_{1} > M} \hat{f}_{\vec{m}}^2 \leq \rho^{M+1} \sum_{\|\vec{m}\|_{1} > M} \frac{\hat{f}_{\vec{m}}^2}{\rho^{\| \vec{m}\|_{1}}} \leq \rho^{M+1} \left \| f  \right\|_{\mathcal{H}_{\beta, \rho}}.
\end{equation}

For the estimation error, notice that $\hat{f}_{\vec{m},N}$ in \eqref{eq:kme_MC} is an unbiased estimator of $\hat{f}_{\vec{m}}$, and we have
\begin{equation*}
\begin{split}
 \mathbb{E}\left[ \sum_{\|\vec{m}\|_{1}\leq M} \left ( \hat{f}_{\vec{m}} - \hat{f}_{\vec{m},N}\right)^{2}\right] & = \sum_{\|\vec m\|_{1}\leq M} \mathbb{E} \left[ \left( \frac{1}{N}\sum_{i = 1}^{N}\psi_{\vec m}(X_{i}) \bm{W}^{1-\beta}(X_{i})- \hat{f}_{\vec m} \right)^{2}\right]  \\
&= \frac{1}{N^{2}}\sum_{\|\vec m\|_{1}\leq M}\sum_{i,j=1}^{N} \operatorname{cov} \left(\psi_{\vec{m}}(X_{i})\bm{W}^{1-\beta}(X_{i}), \psi_{\vec{m}}(X_{j})\bm{W}^{1-\beta}(X_{j}) \right).
\end{split} 
\end{equation*}
Let $h_{\vec{m}}: = \psi_{\vec m}\bm{W}^{1-\beta} $. Since $\{X_{n}\}$ is stationary, the estimation error can be bounded by
\begin{equation}\label{eq:appd_3}
\mathbb{E}\left[ \sum_{\|\vec{m}\|_{1}\leq M} \left ( \hat{f}_{\vec{m}} - \hat{f}_{\vec{m},N}\right)^{2}\right]  \leq \frac{1}{N}\sum_{\|\vec m\|_{1}\leq M} \left[\operatorname{var} \left( h_{\vec{m}}(X_{1})  \right) + 2\sum_{k = 1}^{N-1} \operatorname{cov}\left( h_{\vec{m}}(X_{1}), h_{\vec{m}}(X_{1+k}) \right) \right],
\end{equation}
where the covariance terms represent the non-i.i.d. feature of $\{X_{i}\}$. To bound the variance term in \eqref{eq:appd_3}, consider
\begin{equation}\label{eq:appd_1}
\operatorname{var} \left( h_{\vec{m}}(X_{1})  \right) \leq \int_{\BR^{d}} h_{\vec m}^{2}(\bm{x}) f(\bm{x}) \td \bm{x} =  \int_{\BR^{d}} \psi_{\vec m}^{2}(\bm{x}) \bm{W}^{2-2\beta}(\bm{x}) f(\bm{x}) \td \bm{x}.
\end{equation}
Moreover, by the decay rate of $f$ \eqref{eq:de_rate} and the assumption $\beta \in [\frac{1}{2}, \frac{1}{1+\rho}]$, we have a new decay rate for $f \in \mathcal{H}_{\beta, \rho}$,
\begin{equation} \label{eq:appd_2}
\left|f(x)\right|  \leq (2\pi)^{(\beta-1)\frac{d}{2}} \left(1-\rho^2\right)^{-\frac{d}{4}}\left\| f \right \|_{\mathcal{H}_{\beta, \rho}} \bm{W}^{2\beta-1}(\bm{x}).
\end{equation}
Let $C_{f}: =  (2\pi)^{(\beta-1)\frac{d}{2}} \left(1-\rho^2\right)^{-\frac{d}{4}}\left\| f \right \|_{\mathcal{H}_{\beta, \rho}}$, and the variance in Eq. \eqref{eq:appd_1} can be further bounded by
\begin{equation*}
\operatorname{var} \left( h_{\vec{m}}(X_{1})  \right) \leq C_{f} \int_{\BR^{d}} \psi_{\vec m}^{2}(\bm{x}) \bm{W}(\bm{x}) \td \bm{x} =C_{f}.
\end{equation*}
To bound the covariance term in \eqref{eq:appd_3}, we apply the Davydov's covariance inequality \eqref{eq:Dav} with $h = h_{\vec{m}}$. Set $p = q = 2+\epsilon$ in \eqref{eq:Dav} for $\epsilon \in (0,2)$, we have
\begin{equation}\label{eq:Dav1}
       \left|\operatorname{cov}(h_{\vec{m}}(X_{1}), h_{\vec{m}}(X_{1+k}))\right| \leq  12 \alpha(k)^{\frac{\epsilon}{2+\epsilon}} \left(\mathbb{E}\left[|h_{\vec{m}}(X_{1})|^{2+\epsilon}\right]  \right)^{\frac{2}{2+\epsilon}} \leq 12 \alpha(k)^{\frac{\epsilon}{2+\epsilon}}\left(\mathbb{E}\left[h_{\vec{m}}^{4}(X_{1})\right]  \right)^{\frac{1}{2}},
\end{equation}
where we have used Jensen's inequality in the last step. 

To bound the fourth moment involved in \eqref{eq:Dav1}, we introduce the following expansion formula \cite{carlitz1962product} for the normalized Hermite polynomials.
\begin{equation*}
\psi_{m}^{2} = \sum_{r = 0}^{m} 2^{-r} {m \choose r}\frac{\sqrt{(2m-2r)!}}{(m-r)!} \psi_{2m-2r}.
\end{equation*}
Then, together with the decay rate of $f$ in \eqref{eq:appd_2}, the fourth moment in \eqref{eq:Dav1} can be bounded by
\begin{equation}\label{eq:appd_5}
\begin{split}
\mathbb{E}\left[h_{\vec{m}}^{4}(X_{1})\right] & = \int_{\BR^{d}} \psi_{\vec{m}}^{4}(\bm{x}) \bm{W}^{4- 4\beta}(\bm{x}) f(\bm{x}) \td \bm{x}  \leq C_{f}  \int_{\BR^{d}} \psi_{\vec{m}}^{4}(\bm{x}) \bm{W}^{3- 2\beta}(\bm{x}) \td \bm{x}  \\
& \leq C_{f}  \int_{\BR^{d}} \psi_{\vec{m}}^{4}(\bm{x}) \bm{W}(\bm{x}) \td \bm{x} \leq C_{f} \prod_{i=1}^{d} \sum_{r = 0}^{m_{i}} 2^{-2r} {m_{i} \choose r}^{2}\frac{(2m_{i}-2r)!}{[(m_{i}-r)!]^{2}}.
\end{split}
\end{equation}
To simplify the upper bound in \eqref{eq:appd_5}, recall the Stirling's approximation of the factorials
\begin{equation*}
    \sqrt{2\pi} n^{n + \frac{1}{2}}e^{-n} \leq n! \leq e n^{n+\frac{1}{2}}e^{-n},
\end{equation*}
which leads to
\begin{equation*}
    \frac{\sqrt{(2\ell)!}}{\ell !} \leq \sqrt{\frac{e}{2\pi}} \frac{(2\ell)^{\ell+\frac{1}{4}}}{\ell^{\ell+\frac{1}{2}}} < 2^{\ell}.
\end{equation*}
As a result, the summation in eq. \eqref{eq:appd_5} can be controlled by
\begin{equation} \label{eq:appd_6}
 \sum_{r = 0}^{m_{i}} 2^{-2r} {m_{i} \choose r}^{2}\frac{(2m_{i}-2r)!}{[(m_{i}-r)!]^{2}}<  \sum_{r=0}^{m_{i}} {m_{i} \choose r}^{2} 2^{2m_{i}-4r} <  2^{2m_{i}} \left[ \sum_{r=0}^{m_{i}} {m_{i} \choose r} \left( \frac{1}{4} \right)^{r}\right]^{2}= \left( \frac{5}{2} \right)^{2m_{i}}.
\end{equation}
Combining \eqref{eq:Dav1}, \eqref{eq:appd_5} and \eqref{eq:appd_6},  we have
\begin{equation*}
|\operatorname{cov}(h_{\vec{m}}(X_{1}), h_{\vec{m}}(X_{1+k}))| \leq  12 \alpha(k)^{\frac{\epsilon}{2+\epsilon}}C_{f}^{\frac{1}{2}} \left( \frac{5}{2} \right)^{\|\vec m\|_{1}}.
\end{equation*}
Finally, notice that for a positive integer $j$ the total number of different $d$-dimensional multi-indices such that $\|\vec{m}\|_{1}=j$ is ${j+d-1 \choose d-1}$, and 
\begin{equation*}
{j+d-1 \choose d-1}  =  \frac{(j+1)(j+2)\cdots (j+d-1)}{(d-1)!} < \frac{(j+d-1)^{d-1}}{(5/2)^{d-3}},
\end{equation*}
then the estimation error in \eqref{eq:appd_3} can be bounded by
\begin{equation} \label{eq:appd_7}
\begin{split}
\mathbb{E}\left[ \sum_{\|\vec{m}\|_{1}\leq M} \left ( \hat{f}_{\vec{m}} - \hat{f}_{\vec{m},N}\right)^{2}\right] & \leq \frac{1}{N} \left[ (M+1)^{d}C_{f} + 24 C_{f}^{\frac{1}{2}}\left( \sum_{k=1}^{N-1} \alpha(k)^{\frac{\epsilon}{2+\epsilon}}  \right) \left( \sum_{j=0}^{M} {j+d-1 \choose d-1}  \left( \frac{5}{2} \right)^{j}\right)\right] \\
& \leq  \frac{1}{N} \left[ (M+1)^{d}C_{f} + 24 C_{\beta,\epsilon} C_{f}^{\frac{1}{2}} \left( \sum_{j=0}^{M} (j+d-1) ^{d-1}\left( \frac{5}{2} \right) ^{j-d+3} \right)\right], \\
\end{split}
\end{equation}
where we have used the assumption \eqref{eq:cond2}. With $d\leq \frac{3}{2}M + 1$, the summation in \eqref{eq:appd_7} can be bounded by
\begin{equation*}
\sum_{j=0}^{M} (j+d-1) ^{d-1}\left( \frac{5}{2} \right) ^{j-d+3} < (M+d-1) ^{d-1}\sum_{j=0}^{M}\left( \frac{5}{2} \right) ^{j-d+3} \leq M^{d-1} \sum_{j=0}^{M}\left( \frac{5}{2} \right) ^{j+2} < M^{d-1} \left( \frac{5}{2} \right) ^{M+3},
\end{equation*}
which leads to
\begin{equation} 
\mathbb{E}\left[ \sum_{\|\vec{m}\|_{1}\leq M} \left ( \hat{f}_{\vec{m}} - \hat{f}_{\vec{m},N}\right)^{2}\right] < \frac{1}{N} \left[ (M+1)^{d}C_{f} + 24 C_{\beta,\epsilon} C^{\frac{1}{2}}_{f} M^{d-1} \left( \frac{5}{2} \right) ^{M+3}\right].\nonumber
\end{equation}
Together with \eqref{eq:decay_bias}, we obtain the error bound in \eqref{eq:error_bound}.



\end{document}